\documentclass[10pt]{amsart}
\usepackage{amsmath}
\usepackage{amsxtra}
\usepackage{amscd}
\usepackage{amsthm}
\usepackage{amsfonts}
\usepackage{amssymb}
\usepackage{eucal}

\newtheorem{thm}[subsection]{Theorem}
\newtheorem{cor}[subsection]{Corollary}
\newtheorem{lem}[subsection]{Lemma}
\newtheorem{prop}[subsection]{Proposition}
\newtheorem{propconstr}[subsection]{Proposition-Construction}

\theoremstyle{definition}

\theoremstyle{remark}

\newcommand{\propconstrref}[1]{Proposition-Construction~\ref{#1}}
\newcommand{\thmref}[1]{Theorem~\ref{#1}}
\newcommand{\secref}[1]{Sect.~\ref{#1}}
\newcommand{\lemref}[1]{Lemma~\ref{#1}}
\newcommand{\propref}[1]{Proposition~\ref{#1}}
\newcommand{\corref}[1]{Corollary~\ref{#1}}

\newcommand{\nc}{\newcommand}
\nc{\renc}{\renewcommand}
\nc{\ssec}{\subsection}
\nc{\sssec}{\subsubsection}
\nc{\on}{\operatorname}

\nc\ol{\overline}
\nc\wt{\widetilde}
\nc\tboxtimes{\wt{\boxtimes}}
\nc{\alp}{\alpha}

\nc{\ZZ}{{\mathbb Z}}
\nc{\NN}{{\mathbb N}}
\nc{\CC}{{\mathbb C}}
\nc{\OO}{{\mathbb O}}
\renc{\SS}{{\mathbb S}}
\nc{\DD}{{\mathbb D}}
\nc{\GG}{{\mathbb G}}

\nc{\Fq}{{\mathbb F}_q}
\nc{\Fqb}{\ol{{\mathbb F}_q}}
\nc{\Ql}{\ol{{\mathbb Q}_\ell}}
\nc{\id}{\text{id}}
\nc\X{\mathcal X}

\nc{\Hom}{\on{Hom}}
\nc{\Lie}{\on{Lie}}
\nc{\Loc}{\on{Loc}}
\nc{\Pic}{\on{Pic}}
\nc{\Bun}{\on{Bun}}
\nc{\IC}{\on{IC}}
\nc{\Aut}{\on{Aut}}
\nc{\rk}{\on{rk}}
\nc{\Sh}{\on{Sh}}
\nc{\Perv}{\on{Perv}}
\nc{\pos}{{\on{pos}}}
\nc{\Conv}{\on{Conv}}
\nc{\Sph}{\on{Sph}}
\nc{\Sym}{\on{Sym}}
\nc{\BunBb}{\overline{\Bun}_B}
\nc{\Buno}{\overset{o}{\Bun}}
\nc{\BunPb}{{\overline{\Bun}_P}}
\nc{\BunBM}{\overline{\Bun}_{B(M)}}
\nc{\BunPbw}{{\widetilde{\Bun}_P}}
\nc{\BunBP}{\widetilde{\Bun}_{B,P}}
\nc{\GUb}{\overline{G/U}}
\nc{\GUPb}{\overline{G/U(P)}}

\nc{\Hhom}{\underline{\on{Hom}}}
\nc\syminfty{\on{Sym}^{\infty}}
\nc\lal{\ol{\lambda}}
\nc\xl{\ol{x}}
\nc\thl{\ol{\theta}}
\nc\nul{\ol{\nu}}
\nc\mul{\ol{\mu}}
\nc\Sum\Sigma
\nc{\oX}{\overset{o}{X}{}}
\nc{\hl}{\overset{\leftarrow}h}
\nc{\hr}{\overset{\rightarrow}h}
\nc{\M}{{\mathcal M}}
\nc{\N}{{\mathcal N}}
\nc{\F}{{\mathcal F}}
\nc{\D}{{\mathcal D}}
\nc{\Q}{{\mathcal Q}}
\nc{\Y}{{\mathcal Y}}
\nc{\G}{{\mathcal G}}
\nc{\E}{{\mathcal E}}
\nc{\CalC}{{\mathcal C}}
\nc\Dh{\widehat{\D}}

\nc{\C}{{\mathcal C}}
\nc{\K}{{\mathcal K}}
\renewcommand{\H}{{\mathcal H}}

\nc{\T}{{\mathcal T}}
\nc{\V}{{\mathcal V}}
\renc{\P}{{\mathcal P}}
\nc{\A}{{\mathcal A}}
\nc{\B}{{\mathcal B}}
\nc{\U}{{\mathcal U}}
\renewcommand{\L}{{\mathcal L}}
\nc{\Gr}{\on{Gr}}

\nc{\frn}{{\check{\mathfrak u}(P)}}
\nc{\p}{\mathfrak p}
\nc{\q}{\mathfrak q}
\nc\f{{\mathfrak f}}

\nc{\qo}{{\mathfrak q}}
\nc{\po}{{\mathfrak p}}
\nc{\s}{{\mathfrak s}}
\nc\w{\text{w}}

\nc\mathi\iota
\nc\Spec{\on{Spec}}
\nc\Mod{\on{Mod}}
\nc{\tw}{\widetilde{\mathfrak t}}
\nc{\pw}{\widetilde{\mathfrak p}}
\nc{\qw}{\widetilde{\mathfrak q}}
\nc{\jw}{\widetilde j}

\nc{\grb}{\overline{\Gr}}
\nc{\I}{\mathcal I}

\nc{\lambdach}{{\check\lambda}}
\nc{\Lambdach}{{\check\Lambda}{}}
\nc{\much}{{\check\mu}}
\nc{\omegach}{{\check\omega}}
\nc{\nuch}{{\check\nu}}
\nc{\etach}{{\check\eta}}
\nc{\alphach}{{\check\alpha}}
\nc{\betach}{{\check\beta}}
\nc{\rhoch}{{\check\rho}}
\nc{\ch}{{\check h}}

\nc{\Hb}{\overline{\H}}

\nc{\BA}{{\mathbb{A}}}
\nc{\BC}{{\mathbb{C}}}
\nc{\BG}{{\mathbb{G}}}
\nc{\BM}{{\mathbb{M}}}
\nc{\BN}{{\mathbb{N}}}
\nc{\BP}{{\mathbb{P}}}
\nc{\BR}{{\mathbb{R}}}
\nc{\BZ}{{\mathbb{Z}}}
\nc{\BV}{{\mathbb{V}}}
\nc{\BW}{{\mathbb{W}}}
\nc{\BS}{{\mathbb{S}}}
\nc{\BQ}{{\mathbb{Q}}}

\nc{\CA}{{\mathcal{A}}}
\nc{\CB}{{\mathcal{B}}}

\nc{\CE}{{\mathcal{E}}}
\nc{\CF}{{\mathcal{F}}}
\nc{\CG}{{\mathcal{G}}}
\nc{\CL}{{\mathcal{L}}}
\nc{\CM}{{\mathcal{M}}}
\nc{\CN}{{\mathcal{N}}}
\nc{\CK}{{\mathcal{K}}}
\nc{\CO}{{\mathcal{O}}}
\nc{\CP}{{\mathcal{P}}}
\nc{\CJ}{{\mathcal{J}}}
\nc{\CI}{{\mathcal{I}}}
\nc{\CQ}{{\mathcal{Q}}}
\nc{\CR}{{\mathcal{R}}}
\nc{\CS}{{\mathcal{S}}}
\nc{\CT}{{\mathcal{T}}}
\nc{\CU}{{\mathcal{U}}}
\nc{\CV}{{\mathcal{V}}}
\nc{\CW}{{\mathcal{W}}}
\nc{\CZ}{{\mathcal{Z}}}

\nc{\cM}{{\check{\mathcal M}}{}}
\nc{\csM}{{\check{\mathcal A}}{}}
\nc{\oM}{{\overset{\circ}{\mathcal M}}{}}
\nc{\obM}{{\overset{\circ}{\mathbf M}}{}}
\nc{\oCA}{{\overset{\circ}{\mathcal A}}{}}
\nc{\obA}{{\overset{\circ}{\mathbf A}}{}}
\nc{\ooM}{{\overset{\circ}{M}}{}}
\nc{\osM}{{\overset{\circ}{\mathsf M}}{}}
\nc{\vM}{{\overset{\bullet}{\mathcal M}}{}}
\nc{\nM}{{\underset{\bullet}{\mathcal M}}{}}
\nc{\oD}{{\overset{\circ}{\mathcal D}}{}}
\nc{\obD}{{\overset{\circ}{\mathbf D}}{}}
\nc{\oA}{{\overset{\circ}{\mathbb A}}{}}
\nc{\op}{{\overset{\bullet}{\mathbf p}}{}}
\nc{\cp}{{\overset{\circ}{\mathbf p}}{}}
\nc{\oU}{{\overset{\bullet}{\mathcal U}}{}}
\nc{\oZ}{{\overset{\circ}{\mathcal Z}}{}}
\nc{\ofZ}{{\overset{\circ}{\mathfrak Z}}{}}
\nc{\oF}{{\overset{\circ}{\fF}}}

\nc{\fa}{{\mathfrak{a}}}
\nc{\fb}{{\mathfrak{b}}}
\nc{\fg}{{\mathfrak{g}}}
\nc{\fgl}{{\mathfrak{gl}}}
\nc{\fh}{{\mathfrak{h}}}
\nc{\fj}{{\mathfrak{j}}}
\nc{\fm}{{\mathfrak{m}}}

\nc{\fl}{{\mathfrak{l}}}
\nc{\fn}{{\mathfrak{n}}}
\nc{\fu}{{\mathfrak{u}}}
\nc{\fp}{{\mathfrak{p}}}
\nc{\fr}{{\mathfrak{r}}}
\nc{\fs}{{\mathfrak{s}}}
\nc{\fsl}{{\mathfrak{sl}}}
\nc{\hsl}{{\widehat{\mathfrak{sl}}}}
\nc{\hgl}{{\widehat{\mathfrak{gl}}}}
\nc{\hg}{{\widehat{\mathfrak{g}}}}
\nc{\chg}{{\widehat{\mathfrak{g}}}{}^\vee}
\nc{\hn}{{\widehat{\mathfrak{n}}}}
\nc{\chn}{{\widehat{\mathfrak{n}}}{}^\vee}

\nc{\fA}{{\mathfrak{A}}}
\nc{\fB}{{\mathfrak{B}}}
\nc{\fO}{{\mathfrak{O}}}
\nc{\fC}{{\mathfrak{C}}}
\nc{\fD}{{\mathfrak{D}}}
\nc{\fE}{{\mathfrak{E}}}
\nc{\fF}{{\mathfrak{F}}}
\nc{\fG}{{\mathfrak{G}}}
\nc{\fK}{{\mathfrak{K}}}
\nc{\fL}{{\mathfrak{L}}}
\nc{\fM}{{\mathfrak{M}}}
\nc{\fN}{{\mathfrak{N}}}
\nc{\fP}{{\mathfrak{P}}}
\nc{\fU}{{\mathfrak{U}}}
\nc{\fV}{{\mathfrak{V}}}
\nc{\fZ}{{\mathfrak{Z}}}
\nc{\fz}{{\mathfrak{z}}}

\nc{\bb}{{\mathbf{b}}}
\nc{\bc}{{\mathbf{c}}}
\nc{\bd}{{\mathbf{d}}}
\nc{\bg}{{\mathbf{g}}}
\nc{\be}{{\mathbf{e}}}
\nc{\bj}{{\mathbf{j}}}
\nc{\bn}{{\mathbf{n}}}
\nc{\bp}{{\mathbf{p}}}
\nc{\bq}{{\mathbf{q}}}
\nc{\bu}{{\mathbf{u}}}
\nc{\bv}{{\mathbf{v}}}
\nc{\bx}{{\mathbf{x}}}
\nc{\bs}{{\mathbf{s}}}
\nc{\by}{{\mathbf{y}}}
\nc{\bw}{{\mathbf{w}}}
\nc{\bA}{{\mathbf{A}}}
\nc{\bK}{{\mathbf{K}}}
\nc{\bB}{{\mathbf{B}}}
\nc{\bC}{{\mathbf{C}}}
\nc{\bD}{{\mathbf{D}}}
\nc{\bH}{{\mathbf{H}}}
\nc{\bM}{{\mathbf{M}}}
\nc{\bN}{{\mathbf{N}}}
\nc{\bV}{{\mathbf{V}}}
\nc{\bW}{{\mathbf{W}}}
\nc{\bX}{{\mathbf{X}}}
\nc{\bZ}{{\mathbf{Z}}}
\nc{\bS}{{\mathbf{S}}}

\nc{\sA}{{\mathsf{A}}}
\nc{\sB}{{\mathsf{B}}}
\nc{\sC}{{\mathsf{C}}}
\nc{\sD}{{\mathsf{D}}}
\nc{\sF}{{\mathsf{F}}}
\nc{\sG}{{\mathsf{G}}}
\nc{\sK}{{\mathsf{K}}}
\nc{\sM}{{\mathsf{M}}}
\nc{\sO}{{\mathsf{O}}}
\nc{\sQ}{{\mathsf{Q}}}
\nc{\sP}{{\mathsf{P}}}
\nc{\sZ}{{\mathsf{Z}}}
\nc{\sfp}{{\mathsf{p}}}
\nc{\sr}{{\mathsf{r}}}
\nc{\sg}{{\mathsf{g}}}
\nc{\ssf}{{\mathsf{f}}}
\nc{\sfb}{{\mathsf{b}}}
\nc{\sfc}{{\mathsf{c}}}
\nc{\sd}{{\mathsf{d}}}
\nc{\sk}{{\mathsf{k}}}

\nc{\BK}{{\bar{K}}}

\nc{\tA}{{\widetilde{\mathbf{A}}}}
\nc{\tB}{{\widetilde{\mathcal{B}}}}
\nc{\tG}{{\widetilde{G}}}
\nc{\TM}{{\widetilde{\mathbb{M}}}{}}
\nc{\tO}{{\widetilde{\mathsf{O}}}{}}
\nc{\tU}{{\widetilde{\mathfrak{U}}}{}}
\nc{\TZ}{{\tilde{Z}}}
\nc{\tx}{{\tilde{x}}}
\nc{\tbv}{{\tilde{\bv}}}
\nc{\tfP}{{\widetilde{\mathfrak{P}}}{}}
\nc{\tz}{{\tilde{\zeta}}}
\nc{\tmu}{{\tilde{\mu}}}

\nc{\urho}{\underline{\rho}}
\nc{\uB}{\underline{B}}
\nc{\uC}{{\underline{\mathbb{C}}}}
\nc{\ui}{\underline{i}}
\nc{\uj}{\underline{j}}
\nc{\ofP}{{\overline{\mathfrak{P}}}}
\nc{\oB}{{\overline{\mathcal{B}}}}
\nc{\og}{{\overline{\mathfrak{g}}}}
\nc{\oI}{{\overline{I}}}

\nc{\eps}{\varepsilon}
\nc{\hrho}{{\hat{\rho}}}

\nc{\one}{{\mathbf{1}}}
\nc{\two}{{\mathbf{t}}}

\nc{\Rep}{{\mathop{\operatorname{\rm Rep}}}}
\nc{\Tot}{{\mathop{\operatorname{\rm Tot}}}}
\nc{\Ker}{{\mathop{\operatorname{\rm Ker}}}}
\nc{\Hilb}{{\mathop{\operatorname{\rm Hilb}}}}
\nc{\End}{{\mathop{\operatorname{\rm End}}}}
\nc{\Ext}{{\mathop{\operatorname{\rm Ext}}}}
\nc{\CHom}{{\mathop{\operatorname{{\mathcal{H}}\it om}}}}
\nc{\GL}{{\mathop{\operatorname{\rm GL}}}}
\nc{\gr}{{\mathop{\operatorname{\rm gr}}}}
\nc{\Id}{{\mathop{\operatorname{\rm Id}}}}
\nc{\de}{{\mathop{\operatorname{\rm def}}}}
\nc{\length}{{\mathop{\operatorname{\rm length}}}}
\nc{\supp}{{\mathop{\operatorname{\rm supp}}}}

\nc{\Cliff}{{\mathsf{Cliff}}}
\nc{\Fl}{\on{Fl}}
\nc{\Fib}{{\mathsf{Fib}}}
\nc{\Coh}{{\mathsf{Coh}}}
\nc{\FCoh}{{\mathsf{FCoh}}}

\nc{\cplus}{{\mathbf{C}_+}}
\nc{\cminus}{{\mathbf{C}_-}}
\nc{\cthree}{{\mathbf{C}_*}}
\nc{\Qbar}{{\bar{Q}}}

\nc{\bh}{{\bar{h}}}

\nc{\seq}[1]{\stackrel{#1}{\sim}}

\nc{\wh}{\widehat}

\nc{\mc}{\mathcal}
\nc{\crit}{{\on{crit}}}
\nc{\reg}{{\on{reg}}}
\nc{\mer}{{\on{mer}}}
\renc{\int}{{\on{int}}}
\nc{\nilp}{{\on{nilp}}}
\nc{\nil}{\wt{reg}}
\nc{\ka}{\kappa}
\nc{\g}{\fg}
\nc{\mb}{\mathbf}
\nc{\ren}{ren}
\nc{\la}{\lambda}
\nc{\FZ}{{\mathfrak Z}}
\nc{\Z}{{\mathbb Z}}

\nc{\cG}{{\check G}}

\nc\Res{\on{Res}}
\renc\mod{{\text{-mod}\,{}}}

\nc\semiinf{\frac{\infty}{2}}

\nc\arrowtimes{\overset{\to}\otimes}
\nc\shriektimes{\overset{!}\otimes}

\nc{\ppart}{(\!(t)\!)}
\nc{\Op}{\on{Op}_{\cg}}
\nc{\nOp}{\on{Op}^{\nilp}_{\cg}}
\nc{\cg}{\check \fg}
\nc{\tg}{\wt{\check \fg}}

\begin{document}

\title[Fusion and convolution]{Fusion and convolution: applications to
affine Kac-Moody algebras at the critical level}

\dedicatory{Dedicated to Robert MacPherson on his 60th birthday}

\author{Edward Frenkel}

\address{Department of Mathematics, University of California,
  Berkeley, CA 94720, USA}

\email{frenkel@math.berkeley.edu}

\author{Dennis Gaitsgory}

\address{Department of Mathematics, Harvard University,
Cambridge, MA 02138, USA}

\email{gaitsgory@math.harvard.edu}

\date{November 2005; revised July 2006}

\maketitle

\section*{Introduction}

The goal of this article is two-fold. On the one hand, it constitutes
the second in the series of papers devoted to the study of the
category of representations of an affine Kac-Moody algebra at the
critical level in connection with the local geometric Langlands
correspondence. On the other hand, we study from the geometric point
of view the concept of fusion of modules over an affine Kac-Moody
algebra, or, more generally, chiral modules over a chiral algebra that
admits an algebraic group as a group of its symmetries.

\ssec{}

Let us explain the first perspective in some detail. In \cite{FG2} we
proposed a framework for the local geometric Langlands
correspondence. Namely, let $G$ be a semi-simple algebraic group (over
$\BC$ or any other algebraically closed field of characteristic $0$),
and let us consider $G\ppart$ as a group ind-scheme. We proposed that
to a ``local Langlands parameter'' $\sigma$, which is a $\cG$--bundle
with a connection on the punctured disc (where $\cG$ is the Langlands
dual group of $G$), there should correspond a category $\C_\sigma$
equipped with an action of $G\ppart$.

Unfortunately, in general we could not characterize $\C_\sigma$ in a
unique way by a universal property. However, we conjectured that this
category is closely connected to the category of representations of
the affine Kac-Moody algebra $\hg$ at the critical level.

\medskip

To describe this connection, let us recall that, according to
\cite{FF,F:wak}, the center $\fZ_{\fg}$ of the completed enveloping
algebra of $\hg_\crit$ is isomorphic to the algebra of functions on
the space $\Op(\D^\times)$ of $\cg$-opers on the punctured disc
$\D^\times$ (where $\cg$ is the Langlands dual Lie algebra of $\g$).
We recall that a $\cg$-oper on $\D^\times$ (or, more generally, on any
curve) is by definition a $\cG$-bundle $\CF_G$ equipped with a
connection $\nabla$, together with an additional datum, namely, a
reduction of $\CF_G$ to a Borel subgroup of $\cG$, which is in a
particular relative position with respect to $\nabla$.

For a fixed $\cg$-oper $\chi$, which we regard as a character of
$\fZ_\fg$, we can consider the sub-category $\hg_\crit\mod_\chi$
of the category $\hg_\crit\mod$ of all discrete modules at
the critical level, on which the center acts according to this character.
This category carries a canonical action of $G\ppart$ via its adjoint 
action on $\hg_\crit$.

We proposed in \cite{FG2} that $\hg_\crit\mod_\chi$ should be
equivalent to the sought-after category $\C_\sigma$, where $\sigma$
corresponds to the pair $(\CF_G,\nabla)$ underlying the oper
$\chi$. This guess entails a far-reaching corollary that the
categories $\hg_\crit\mod_{\chi^1}$ and $\hg_\crit\mod_{\chi^2}$ for
two different central characters $\chi^1$ and $\chi^2$ are equivalent
once we identify the underlying local systems $(\CF^1_G,\nabla^1)$ and
$(\CF^2_G,\nabla^2)$.

\medskip

In \cite{FG2} we discussed various consequences of this proposal, and
in particular, a concrete conjecture about the structure of the full
subcategory $\hg_\crit\mod_\nilp^{I,m}$ of the category of
$\hg_\crit$-modules. Its objects are $I$-monodromic $\hg_\crit$-modules
(where $I$ is the Iwahori subgroup of $G\ppart$), which are
supported, as $\fZ_{\fg}$-modules, over a certain sub-scheme 
$\nOp\subset\Op(\D^\times)$. 

Our conjecture was that the derived category of
$\hg_\crit\mod_\nilp^{I,m}$ is equivalent to the derived category of
quasi-coherent sheaves on the scheme $\tg/\cG\underset{\cg/\cG}\times
\nOp$, where $\tg\to \cg$ is Grothendieck's alteration, and $\nOp\to
\cg/\cG$ is a canonical residue map $\Res^\nilp$, introduced in
\cite{FG2}. In \cite{FG2} we proved a weaker statement that certain
quotients of these two categories are equivalent.

\medskip

In order to prove this conjecture in full we need to know a wide range of
results about the structure of the category $\hg_\crit\mod$.
One class of such results has to do with the Harish-Chandra
convolution action of D-modules on $G\ppart$ on the derived 
category $D(\hg_\crit\mod)$.

\ssec{}

In \cite{FG2} we explained, following \cite{BD}, that the appropriately
defined $I$-equivariant derived category $D^b(\on{D}(\Fl_G)_\crit\mod)^{I}$ 
of critically twisted D-modules on the affine flag scheme
$\Fl_G=G\ppart/I$ acts on the $I$-equivariant derived category 
$D^b(\hg_\crit\mod)^{I}$ of $\hg_\crit$-modules.

This structure is called the Harish-Chandra convolution action, and
its origin is the action of the affine Hecke algebra (whose
categorification is $D^b(\on{D}(\Fl_G)_\crit\mod)^{I}$) on the space
of $I$-invariant vectors in any representation of the group $G$ over a
local non-archimedian field.

\medskip

The category $D^b(\on{D}(\Fl_G)_\crit\mod)^{I}$ contains special objects,
the so-called central sheaves $\CZ_V, V \in \on{Rep}(\cG)$,
introduced in \cite{Ga}. They correspond to the central elements of
the affine Hecke algebra. 

The main result of the present paper, is that all objects of 
$\hg_\crit\mod^{I}$ are eigen-objects with respect to the 
functor of convolution with the central sheaves. This is one of the
crucial steps in our project describing $\hg_\crit\mod_\nilp^{I,m}$
in terms of quasi-coherent sheaves. Let us now explain more 
precisely what being an eigen-object means in our set-up. 

\medskip

In \secref{statement of Iwahori theorem}, we explain that the support
over the ind-scheme $\Spec(\fZ_\fg)$ of any object from 
$\hg_\crit\mod^{I}$ is contained in the ind-subscheme 
$\Spec(\wt\fZ^{\int,\nilp}_\fg)$ that corresponds to opers which,
as local systems on $\D^\times$, have regular singularities
and a unipotent monodromy. We show also that to each
$V\in \Rep(\cG)$ there corresponds a vector bundle
$\CV_{\wt\fZ^{\int,\nilp}_\fg}$ over $\Spec(\wt\fZ^{\int,\nilp}_\fg)$.
The geometric meaning of this vector bundle is the following:
for a $\BC$-point $\chi\in \Spec(\wt\fZ^{\int,\nilp}_\fg)\simeq 
\on{Op}_{\cg}(\D^\times)$, the fiber
$\CV_\chi$ of $\CV_{\wt\fZ^{\int,\nilp}_\fg}$ is
isomorphic to the fiber at the origin of the canonical (i.e., Deligne's)
extension from $\D^\times$ to $\D$ of the local system $(\CF_G,\nabla)$, 
underlying $\chi$.

The vector bundle $\CV_{\wt\fZ^{\int,\nilp}_\fg}$
is equipped with a nilpotent endomorphism,
corresponding to the monodromy of the underlying oper.

\medskip

The main result of the present paper, \thmref{main}, states that
for any $\CM\in \hg_\crit\mod^{I}$ and $V\in \Rep(\cG)$ we have
a canonical isomorphism of $\hg_\crit$-modules:
\begin{equation}    \label{main formula}
\CZ_V \star \CM \simeq \CV_{\wt\fZ^{\int,\nilp}_\fg} 
\underset{\wt\fZ^{\int,\nilp}_\fg}\otimes \CM,
\end{equation}
such that the action of the monodromy on the left hand side, coming
from the definition of $\CZ_V$ via the nearby cycles functor as in
\cite{Ga}, goes under this isomorphism to the nilpotent endomorphism
of $\CV_{\wt\fZ^{\int,\nilp}_\fg}$ mentioned above. Moreover, the
isomorphisms of \eqref{main formula} for different representations $V$
are compatible, in a natural sense, with the operation of tensor
product.

Thus, for example, if $\CM$ is an $I$-equivariant $\hg_\crit$--module
on which the center $\fZ_{\fg}$ acts via the character corresponding
to a particular $\chi\in \wt\fZ^{\int,\nilp}_\fg$, then the
isomorphism \eqref{main formula} becomes
$$
\CZ_V \star \CM \simeq \CV_\chi \underset{\CC}\otimes \CM,
$$
where $\CV_\chi$ is as above.

\ssec{}

The isomorphism stated in \eqref{main formula} has an analogue for the
category of $G[[t]]$-equivariant $\hg_\crit$-modules (rather than
$I$-equivariant ones). In this case the role of the category
$\on{D}(\Fl_G)_\crit\mod^{I}$ is played by the category
$\on{D}(\Gr_G)_\crit\mod^{G[[t]]}$ of $G[[t]]$-equivariant D-modules
on the affine Grassmannian $\Gr_G = G\ppart/G[[t]]$. This category is
a tensor category that is equivalent to the category $\Rep(\cG)$ of
finite-dimensional representations of $\cG$, see \cite{MV} (it may
also be thought of as a categorification of the spherical Hecke
algebra).

\medskip

Hence, for each $V \in \Rep(\cG)$ we have the corresponding 
object $\CF_V$ in the category $\on{D}(\Gr_G)_\crit\mod^{G[[t]]}$.
A spherical version of our main theorem, \thmref{main, spherical},
states that for every $\CM\in \hg_\crit\mod^{G[[t]]}$ there
is a canonical isomorphism
\begin{equation}    \label{main formula, spherical}
\CF_V \star \CM \simeq \CV_{\wt\fZ^{\int,\nilp}_\fg} 
\underset{\wt\fZ^{\int,\nilp}_\fg}\otimes \CM.
\end{equation}
Moreover, we obtain that the support of every such $\CM$
over $\Spec(\wt\fZ^{\int,\nilp}_\fg)$ is contained in the
ind-subscheme $\Spec(\wt\fZ^{\on{m.f}}_\fg)$,
corresponding to monodromy-free opers.

\medskip

Tautologically, the isomorphism \eqref{main formula, spherical}
is a particular case of that of \eqref{main formula}. Let us specialize
further to the case when $\CM\in \hg_\crit\mod^{G[[t]]}$ is the
vacuum module $\BV_\crit$. In the latter case, the corresponding
isomorphism 
\begin{equation}    \label{main formula, vacuum}
\CF_V \star \BV_\crit \simeq \CV_{\fz_\fg} 
\underset{\fz_\fg}\otimes \BV_\crit.
\end{equation}
is known, due to \cite{BD}.

The main idea of this paper is that one can derive results such as
\eqref{main formula} and \eqref{main formula, spherical} from the
special case when $\CM=\BV_\crit$. This is based on the operation
of fusion product, and this brings us to the discussion of the 
second perspective in which one can view this paper.

\ssec{}

Let $\CA$ be a chiral algebra on a curve $X$, and let $\CM_1,\CM_2$ and
$\CN$ be chiral $\CA$-modules. In this case one can consider the set
of {\it chiral pairings} $\{\CM_1,\CM_2\}\to \CN$. These notions were 
introduced by A. Beilinson and V. Drinfeld in \cite{CHA}. 

When $\CA$ is a chiral algebra attached to a conformal vertex algebra
${\mathsf V}$, and $\CM_1,\CM_2,\CN$ are obtained from ${\mathsf
V}$-modules, the notion of chiral pairing is similar to that of
intertwining operator between modules over ${\mathsf V}$, see
\cite{FHL}.

If for fixed $\CM_1$ and $\CM_2$ the functor that assigns to $\CN$ the
set of chiral pairings
$\{\CM_1,\CM_2\}\to \CN$ is representable, we shall call the
representing object {\it the fusion product} of $\CM_1$ and $\CM_2$.

\medskip

Assume now that $\CA$ is the chiral algebra $\CA_{\fg,\kappa}$,
attached to a semi-simple Lie algebra $\fg$ and a symmetric invariant
form $\kappa:\fg\otimes \fg\to \BC$. Then the category of chiral
$\CA_{\fg,\kappa}$-modules, supported at a given point $x\in X$ is
naturally equivalent to the category $\hg_\kappa\mod$ of
representations of the affine Kac-Moody algebra $\hg$ at the level
$\kappa$. (More generally, we can consider chiral algebras that admit
$G$ as a group of symmetries, see \secref{general alg}).

Let us suppose now that $\kappa$ is non-positive and integral. In this
case, to every $V\in \Rep(\cG)$ we can attach a chiral
$\CA_{\fg,\kappa}$-module, denoted $\Gamma(\Gr_{G,X},\CF_{V,X})$, or
$\CF_{V,X}\star \CA_{\fg,\kappa}$.  It is constructed using the
D-module $\CF_V\in \on{D}(\Gr_G)_\kappa\mod^{G[[t]]}$, which is
well-defined for every integral $\kappa$.

Let $\CM_2=\CM$ be any $G[[t]]$-equivariant (resp., $I$-equivariant)
$\hg_\kappa$-module. We set $\CM_1:=\Gamma(\Gr_{G,X},\CF_{V,X})$, 
or more generally, $\CM_1:=\Gamma(\Gr_{G,X},\CF_{V,X})\otimes \CE$,
where $\CE$ is a local system on the punctured curve $X-x$ with
a nilpotent monodromy around $x$.

\medskip

Our key technical tool is the assertion (see \thmref{representability,
spherical} and \thmref{representability, Iwahori}) that in this case
fusion products exist and are given by $\CF_V\star \CM$ in the
$G[[t]]$-equivariant case, and by $$\Bigl((\CZ_V\star \CM)\otimes
\Psi(\CE)\Bigr)_N$$ in the $I$-equivariant case, where $\Psi(\CE)$
denotes nearby cycles of $\CE$ at $x$, viewed as a vector space with a
nilpotent operator, and the subscript $N$ stands for taking
coinvariants of the monodromy.

\medskip

The proofs of the main Theorems, namely, \ref{main, spherical} and \ref{main},
are obtained by showing that the right-hand sides of the stated
isomorphisms also represent the above functors of chiral pairings, and
it is this last assertion that uses the result of \cite{BD} about the
isomorphism \eqref{main formula, vacuum}.

Let us remark that the proofs of the main theorems do not actually use
the full statements of Theorems \ref{representability, spherical} and
\ref{representability, Iwahori}, but only the existence of the
corresponding maps in one direction, and their compatibility
with tensor products of representations.



\ssec{}

Let us briefly describe the way this paper is organized. It is divided
into two parts.

\medskip

In Part I we discuss in detail the $G[[t]]$-equivariant situation.
In \secref{sect 1} we state \thmref{main, spherical}; in \secref{sect 2}
we discuss the relationship between the fusion product and Harish-Chandra
convolution and state \thmref{representability, spherical}. In
\secref{sect 3} we complete the proof of \thmref{main, spherical},
and finally in \secref{proof of rep} we prove a generalization of
\thmref{representability, spherical} in the framework of chiral algebras,
endowed with a Harish-Chandra action of the group $G$.

In Part II we show how to modify the material of Part I for the
$I$-equivariant situation. Thus, in \secref{sect 5} we state
\thmref{main}, and in \secref{sect 6} we derive it from (a part of)
\thmref{representability, Iwahori}.  In \secref{sect 7} we prove a
suitable generalization of \thmref{representability, Iwahori}.

\medskip

In the appendix, \secref{app A}, we give a proof of the fact
that the chiral bracket on chiral algebras such as $\CA_{\fg,\kappa}$
can be described using D-modules on the Beilinson-Drinfeld
Grassmannian.

\medskip

The notation in this paper follows closely that of \cite{FG2}.

\ssec{Acknowledgments}

D.G. would like to thank Sasha Beilinson for patient explanations
and stimulating discussions.

The research of E.F. was supported by the DARPA grant HR0011-04-1-0031
and by the NSF grant DMS-0303529.

\vspace*{10mm}

\centerline{\bf {\large Part I: The spherical case}}

\section{Convolution at the critical level}    \label{sect 1}
 
\subsection{Recollections}    \label{recol}

Let $\fg$ be a simple finite-dimensional Lie algebra. For an invariant
inner product $\kappa$ on $\fg$ (which is unique up to a scalar) define
the central extension $\hg_\kappa$ of the formal loop algebra $\fg
\otimes \BC\ppart$ which fits into the short exact sequence
$$
0 \to \BC {\mb 1} \to \hg_\kappa \to \fg \otimes \BC\ppart \to 0.
$$
This sequence is split as a vector space, and the commutation
relations read
\begin{equation}    \label{KM rel}
[x \otimes f(t),y \otimes g(t)] = [x,y] \otimes f(t) g(t) - \kappa(x,y)\cdot
\on{Res}(f \cdot dg)\cdot  {\mb 1},
\end{equation}
and ${\mb 1}$ is a central element. The Lie algebra $\hg_\kappa$ is
the {\em affine Kac-Moody algebra} associated to $\kappa$. We will
denote by $\hg_\kappa\mod$ the category of {\em discrete}
representations of $\hg_{\kappa}$ (i.e., such that any vector is
annihilated by $\fg \otimes t^n\BC[[t]]$ for sufficiently large $n$),
on which ${\mb 1}$ acts as the identity.

\medskip

Let $\kappa_{Kil}$ be the Killing form:
$$\kappa_{Kil}(x,y) = \on{Tr} (\on{ad} (x) \circ \on{ad}(y)).$$
A level $\kappa$ is called critical (resp., positive, negative, irrational)
if $\kappa=c\cdot \kappa_{Kil}$ and $c=-\frac{1}{2}$
(resp., $c+\frac{1}{2}\in \BQ^{>0}$, $c+\frac{1}{2}\in \BQ^{<0}$, $c\notin \BQ$). 
A level $\kappa$ is called integral is it is an integral multiple of the standard 
inner product, normalized so that the square length of the maximal root is
equal to $2$.

\medskip

Next, we recall some notation and terminology from the theory of
chiral algebras introduced in \cite{CHA}. Our chiral algebras will be
defined on a smooth algebraic curve $X$. We will fix a point $x\in X$,
and identify D-modules supported at $x$ with underlying
vector spaces. 

We will denote by $\D_x$ and $\D_x^\times$, respectively, the formal
disc and the formal punctured disc around $x$. If we choose a
coordinate $t$ near $x$, we obtain the identifications $\D_x\simeq
\D:=\Spec(\CC[[t]])$ and $\D_x^\times\simeq
\D^\times:=\Spec(\CC\ppart)$,

\medskip

Following \cite{CHA} (see also \cite{FG2}, Sect. 10), we associate to
$\fg$ the Lie-* algebra $L_\fg = \fg\otimes D_X$, and for each level
$\kappa$ its central extension $L_{\fg,\kappa}$ by means of
$\omega_X$. By definition, the chiral algebra $\CA_{\fg,\kappa}$ is
the quotient of the chiral universal envelope of $L_{\fg,\kappa}$
obtained by identifying the two copies of $\omega_X$.

\medskip

For the rest of this section we shall fix the level $\kappa$ to be 
critical. Let $\fz_\fg$ be the center of $\CA_{\fg,\crit}$, viewed as a
commutative D-algebra on a curve $X$. Denote by $\fz_{\fg,x}$ the
fiber of $\fz_\fg$ at the point $x\in X$, and let $\fZ_\fg$ be the
topological commutative algebra corresponding to $\fz_\fg$ and $x$, as
defined in \cite{CHA}, Sect. 3.6.18.

There exists a canonical map from $\fZ_\fg$ to the center of the
category $\hg_\crit\mod$. (The latter identifies tautologically with
the center of the corresponding completed universal enveloping algebra
$\wt{U}_\crit(\hg)$, see \cite{FG2}, Sect. 5.1 for more details.) One
can show that the map from $\fZ_\fg$ to $Z(\wt{U}_\crit(\hg))$ is in
fact an isomorphism, but we will not use this fact.

\medskip

Denote by $\check \fg$ the Langlands dual Lie algebra to $\fg$.  Let
$\on{Op}_{\cg,X}$ be the D-scheme of $\cg$-opers on $X$ \cite{BD}. The
following isomorphism is proved in \cite{FF,F:wak}:
\begin{equation} \label{isom with opers global}
\fz_{\fg} \simeq \on{Fun} (\on{Op}_{\cg,X}).
\end{equation}

Let consider also the spaces of $\cg$-opers $\on{Op}_{\cg}(\D_x)$ and
$\on{Op}_{\cg}(\D_x^\times)$ on the formal disc $\D_x$ and the
formal punctured disc $\D_x^\times$, respectively. 
The former is a scheme of infinite type, and the latter is an 
ind-scheme. We refer the reader to Part I of \cite{FG2} for a
detailed discussion of opers. 

The isomorphism \eqref{isom with opers global} implies that:
\begin{equation} \label{isom with opers reg}
\fz_{\fg,x} \simeq \on{Fun} (\on{Op}_{\cg}^\reg),
\end{equation}
where $\on{Op}_{\cg}^\reg:=\on{Op}_{\cg}(\D_x)$, and  
\begin{equation} \label{isom with opers}
\fZ_\fg \simeq \on{Fun} (\on{Op}_{\cg}(\D_x^\times)).
\end{equation}

\subsection{Vector bundles over opers} \label{bundle on opers}

Let us denote by $\cG$ the simple algebraic group of simply-connected
type, corresponding to $\cg$. Then there is a tautological
$\cG$-bundle $\CP_{\cG,\on{Op}_{\cg,X}}$ over $\on{Op}_{\cg,X}$
equipped with a connection along $X$. For any finite-dimensional
representation $V$ of $\cG$ we will denote by $\CV_{\on{Op}_{\cg,X}}$
the corresponding associated vector bundle over $\on{Op}_{\cg,X}$,
equipped with a connection along $X$.

Via the isomorphism \eqref{isom with opers global} we can view
$\CV_{\on{Op}_{\cg,X}}$ as a locally free $\fz_\fg$-module with
a connection, and we shall also denote it by $\CV_{\fz_\fg}$.
(We remark that as a $\fz_\fg$-module, $\CV_{\fz_\fg}$ is
actually free, but there is no canonical choice of generators.)
The fiber of $\CV_{\fz_\fg}$ over $x\in X$ is a locally free
$\fz_{\fg,x}$-module, which we will denote by $\CV_{\fz_{\fg,x}}$.

\medskip

Denote by $\wt\fz_{\fg,x}$ the topological algebra equal to the
completion of $\fZ_{\fg}$ with respect to the ideal defining
$\fz_{\fg,x}$, i.e.,
$$\wt\fz_{\fg,x}:=\underset{n}{\underset{\longleftarrow}{\lim}}\,
\fZ_{\fg}/\CJ_n,$$ where $\CJ_n$ is the closure of the $n$th power
of $\on{ker}(\fZ_{\fg}\to \fz_{\fg,x})$. 
Let $\Spec(\wt\fz_{\fg,x})$ be the resulting ind-subscheme of
$\Spec(\fZ_\fg)$. 

\medskip

(In what follows, for a topological algebra
$A\simeq \underset{\longleftarrow}{\lim}\, A_i$, we shall denote by 
$\Spec(A)$ the corresponding ind-scheme, i.e.,
$\underset{\longrightarrow}{"\lim"}\, \Spec(A_i)$. 
By a vector bundle over an ind-scheme we shall
mean a compatible system of vector bundles over its closed
subschemes.)

\begin{propconstr} \label{extension of bundle over center}  \hfill

\smallskip

\noindent{\em (1)} The vector bundle $\CV_{\fz_{\fg,x}}$ naturally
extends to a vector bundle $\CV_{\wt\fz_{\fg,x}}$ over
$\Spec(\wt\fz_{\fg,x})$.

\smallskip

\noindent{\em (2)} For $V,W\in \Rep(\cG)$ and $U=V\otimes W$, there is
a natural isomorphism
$\CV_{\wt\fz_{\fg,x}}\underset{\wt\fz_{\fg,x}}\otimes
\CW_{\wt\fz_{\fg,x}}\simeq \CU _{\wt\fz_{\fg,x}}$.

\smallskip

\noindent{\em (3)}
Each $\CV_{\wt\fz_{\fg,x}}$ is equipped with a (pro)-nilpotent endomorphism 
$N_{\CV_{\wt\fz_{\fg,x}}}$, and these endomorphisms are compatible with
the identifications of (2) above.

\end{propconstr}

The rest of this subsection is devoted to the
proof of this proposition.

\medskip

By definition, $$\wt\fz_{\fg,x}\simeq
\underset{\fz'_\fg}{\underset{\longleftarrow}{\lim}}\, \fz'_{\fg,x},$$
where $\fz'_\fg$ runs over the family of D-subalgebras of $\fz_\fg$,
such that the ideal $\on{ker}(\fz'_{\fg,x}\to \fz_{\fg,x})$ is
nilpotent (here $\fz'_{\fg,x}$ denotes the fiber of $\fz'_\fg$ at $x$).
Consider the following general set-up:

\medskip

Let $\CB$ be a commutative chiral algebra on $X$, and let $\CV_{X-x}$
be a free finite-rank $\CB$-module, equipped with a compatible
connection along $X$.  Assume that $\CV_{X-x}$ admits an extension to
a free finite rank $\CB$-module $\CV_X$, on which the connection along
$X$ has a pole of order $\leq 1$ at $x$ and its residue
$\Res(\nabla,\CV_X)$, thought of as an endomorphism of the fiber
$\CV_x$ of $\CV_X$ at $x$, is nilpotent.

Let $\CB'\hookrightarrow \CB$ be a chiral subalgebra, such that
$\CB'|_{X-x}\simeq \CB|_{X-x}$, and such that $\on{ker}(\CB'_x\to
\CB_x)$ is nilpotent.

\begin{lem}  \label{extension of chiral modules} 
Under the above circumstances, $\CV_{X-x}$, viewed as a $\CB'$-module
over $X-x$, admits a unique extension to a free finite rank
$\CB'$-module $\CV'_X$ such that $\CV_X\simeq
\CV'_X\underset{\CB'}\otimes \CB$, and such that the connection on
$\CV'_X$ has a pole of order $\leq 1$ at $x$. In this case
$\Res(\nabla,\CV'_X)$ is also nilpotent.
\end{lem}

Applying this lemma in our situation we obtain a locally free sheaf
$\CV_{\fz'_\fg}$ over each $\fz'_\fg$ as above. Moreover, the
formation of $\CV_{\fz'_\fg}$ is compatible with tensor products of
representations, by the uniqueness statement of the lemma.

By taking the fiber of $\CV_{\fz'_\fg}$ at $x$ we obtain a vector
bundle $\CV_{\fz'_\fg,x}$ over each $\Spec(\fz'_{\fg,x})$ as above,
i.e., a vector bundle over $\Spec(\wt\fz_{\fg,x})$. The endomorphisms
$N_{\CV_{\wt\fz_{\fg,x}}}|_{\Spec(\fz'_{\fg,x})}$ are equal to
$\Res(\nabla,\CV_{\fz'_\fg})$.

\medskip

Let us now prove \lemref{extension of chiral modules}. We have the
natural morphisms $\CB\to \CB_x[[t]]$ and $\CB'\to \CB'_x[[t]]$,
compatible with connections. It follows from the Beauville-Laszlo
theorem \cite{BL} (see also \cite{BD}) that the problem of
extension of $\CV_{X-x}$ to $X$ translates to a problem of extensions
of locally free modules with connections over $\CB_x[[t]]$ and
$\CB'_x[[t]]$. Thus, we have to prove the following:

\begin{lem}  \label{extensions of modules}
Let $B'\twoheadrightarrow B$ be a surjection of commutative algebras
with a nilpotent kernel. Let $T'$ be a free module over
$B'\ppart$, endowed with a connection along $t$, and let $T$ be the
corresponding $B\ppart$-module. Let $T_0\subset T$ be a
$B[[t]]$-lattice, such that $t\partial_t\cdot T_0\subset T_0$, and the
endomorphism induced on $T_0/t\cdot T_0$ is nilpotent.  Then there
exists a unique $B'[[t]]$-lattice $T'_0$ in $T'$ with the same
properties, such that $T'_0\underset{B'[[t]]}\otimes B[[t]]=T_0\subset
T$.
\end{lem}
 
\begin{proof}
 
By induction, we can assume that $I:=\on{ker}(B'\to B)$ is such that
$I^2=0$.  Then we have a short exact sequence
$$0\to I\cdot T'\to T'\to T'/I\cdot T'\to 0,$$ where both $I\cdot
T'\simeq I\underset{B}\otimes T$ and $T'/I\cdot T'\simeq T$ are
$B\ppart$-modules.

It is easy to see that $T'$ admits at least one $B'[[t]]$-lattice
$'T'_0$, satisfying
\begin{equation} \label{cond on lattice}
'T'_0\on{mod} \, I\cdot T'=T_0 \text{ and }
'T'_0\cap \, I\cdot T'=I[[t]]\underset{B[[t]]}\otimes T_0.
\end{equation}
 
The operator $t\partial_t$, acting on $'T_0'$ defines a
$B[[t]]$-linear map
$$\phi:T_0\to I[[t]]\underset{B[[t]]}\otimes T/T_0.$$

Any other lattice $T'_0$ in $T'$ that satisfies \eqref{cond on
lattice} differs from $'T'_0$ by a $B[[t]]$-linear operator
$$E:T_0\to I[[t]]\underset{B[[t]]}\otimes T/T_0.$$ The condition for
$T'_0$ to satisfy the assumption of the lemma reads as follows:
$$\phi=[E,t\partial_t].$$ This equation is uniquely solvable by
induction on the order of the pole, since $\phi-k\cdot \on{Id}$ is
invertible whenever $k\neq 0$.

\end{proof}
 
\ssec{A generalization}
 
Let $V^\mu$ be the irreducible $\fg$-module with highest weight $\mu$,
and set $\BV^\mu=\on{Ind}^{\hg_\crit}_{\fg[[t]]}(V^\mu)$.  Recall from
\cite{FG2}, Sect. 7.6 that to any dominant integral
weight $\lambda$ we have attached a subscheme
$\Spec(\fZ^{\lambda,\reg}_{\fg})\subset \Spec(\fZ_{\fg})$.  The
following will be proved in \secref{another proof mon-free}:
\begin{lem}  \label{sup of ind}
The support of $\BV^\mu$ over $\Spec(\fZ_{\fg})$ is contained in
$\Spec(\fZ^{\mu,\reg}_{\fg})$.
\end{lem}

Let $\fz^{\lambda,\reg}_\fg$ be the commutative chiral algebra on $X$,
isomorphic to $\fz_\fg$ over $X-x$ whose fiber at $x$ is
$\Spec(\fZ^{\lambda,\reg}_{\fg})$. According to \cite{FG2}, Sect. 2.9, 
for $V\in \Rep(\cG)$, the locally free sheaf $\CV_{\fz_\fg}$
with a connection, defined on $\fz^{\lambda,\reg}_\fg|_{X-x}$, extends
to a locally free sheaf with a {\it regular connection} over
$\fz^{\lambda,\reg}_\fg$. We will denote the resulting vector bundle
over $\Spec(\fZ^{\lambda,\reg}_{\fg})$ by
$\CV_{\fZ^{\lambda,\reg}_{\fg}}$.

Let $\wt\fZ^{\lambda,\reg}_{\fg}$ be the formal completion of
$\Spec(\fZ_{\fg})$ along $\Spec(\fZ^{\lambda,\reg}_{\fg})$. The
proof of \propconstrref{extension of bundle over center} implies that
$\CV_{\fZ^{\lambda,\reg}_{\fg}}$ naturally extends to a vector
bundle $\CV_{\wt\fZ^{\lambda,\reg}_{\fg}}$ over
$\Spec(\wt\fZ^{\lambda,\reg}_{\fg})$, equipped with a nilpotent
endomorphism $N_{\CV_{\wt\fZ^{\lambda,\reg}_{\fg}}}$, in a way
compatible with tensor products of objects of $\cG$.

\medskip

Let us denote by $\Spec(\fZ^{\int,\reg}_{\fg})$ (resp.,
$\Spec(\wt\fZ^{\int,\reg}_{\fg})$) the ind-subscheme of
$\Spec(\fZ_{\fg})$ equal to the disjoint union
$\underset{\lambda}\sqcup\, \Spec(\fZ^{\lambda,\reg}_{\fg})$ (resp.,
$\underset{\lambda}\sqcup\, \Spec(\wt\fZ^{\lambda,\reg}_{\fg})$),
and by $\CV_{\fZ^{\int,\reg}_{\fg}}$ (resp.,
$\CV_{\wt\fZ^{\int,\reg}_{\fg}}$) the vector bundle on it,
corresponding to $V\in \Rep(\cG)$. We shall denote by
$N_{\CV_{\wt\fZ^{\int,\reg}_{\fg}}}$ the endomorphism of the latter
vector bundle.

\medskip

Let $\Spec(\fZ_{\fg}^{\on{m.f.}}) \subset \Spec(\fZ_{\fg})$ be the
sub-functor that corresponds to opers that are monodromy-free as local
systems.  It is easy to see that we have the following inclusions of
functors:
$$\Spec(\fZ^{\int,\reg}_{\fg}) \subset \Spec(\fZ_{\fg}^{\on{m.f.}})
\subset \Spec(\wt\fZ^{\int,\reg}_{\fg}).$$ (In fact, one can show
that $\Spec(\fZ_{\fg}^{\on{m.f.}})$ is the minimal ind-subscheme of
$\Spec(\wt\fZ^{\int,\reg}_{\fg})$, containing $\Spec(\fZ^{\int,\reg}_{\fg})$,
stable under the action of the
Lie algebroid $\Omega^1(\fZ_{\fg})$, see \cite{FG1}, Sect. 6.12.)

Moreover, from the definitions it follows that
$\Spec(\fZ_{\fg}^{\on{m.f.}})$ is the closed ind-subscheme of
$\Spec(\wt\fZ^{\int,\reg}_{\fg})$, equal to the locus of vanishing
of $N_{\CV_{\wt\fZ^{\int,\reg}_{\fg}}}$ for all $V\in \cG$.

\medskip

Let $\CM$ be an object of $\fZ_{\fg}\mod$, i.e., a (discrete) vector
space, endowed with a continuous action of $\fZ_\fg$. We shall say
that $\CM$ is supported on $\Spec(\wt\fZ^{\lambda,\reg}_{\fg})$ if
every element of $\CM$ is annihilated by some power of the ideal of
$\fZ^{\lambda,\reg}_{\fg}$ in $\fZ_\fg$.

We shall say that a module $\CM$ is supported on
$\Spec(\wt\fZ^{\int,\reg}_{\fg})$ if it is a union (or,
automatically, a direct sum) of module $\CM^\lambda$, where each
$\CM^\lambda$ is supported on
$\Spec(\wt\fZ^{\lambda,\reg}_{\fg})$.\footnote{In fact, this notion
makes sense for any pair of affine ind-schemes, one being a closed
ind-subscheme of the other.}

Given $V\in \Rep(\cG)$ we have a well-defined functor on the category
of $\fZ_\fg$-modules, supported on
$\Spec(\wt\fZ^{\int,\reg}_{\fg})$:
$$\CM\mapsto
\CV_{\wt\fZ^{\int,\reg}_{\fg}}\underset{\wt\fZ^{\int,\reg}_{\fg}}
\otimes \CM.$$ Moreover, these functors come equipped with nilpotent
endomorphisms $N_{\CV_{\wt\fZ^{\int,\reg}_{\fg}}}$, compatible with
tensor products of objects of $\Rep(\cG)$.

The above constructions have a relevance for us due to the following:

\begin{lem}  \label{sup of sph}
If $\CM$ is a $G[[t]]$-integrable $\hg_\crit$-module, then, viewed as
a module over $\fZ_\fg$, it is supported on
$\Spec(\wt\fZ^{\int,\reg}_{\fg})$.
\end{lem}

The proof follows from \lemref{sup of ind}, since every object of
$\hg_\crit\mod^{G[[t]]}$ admits a filtration, such that each
subquotient is isomorphic to a quotient module of some $\BV^\lambda$.
(For a different argument see \secref{another proof mon-free}.)

\ssec{Statement of the main theorem (spherical case)} 

Let $G$ be the group of adjoint type with the Lie algebra $\fg$, so
that the Langlands dual group $\cG$ is simply-connected. We denote by
$\Gr_G = G\ppart/G[[t]]$ the corresponding affine
Grassmannian.\footnote{Note that $\Gr_G$ has connected components
labeled by elements of the fundamental group of $G$}

We will consider the category $\on{D}(\Gr_G)_\crit\mod$ of critically
twisted right D-modules on $\Gr_G$. Since the critical level
is integral, this category is equivalent to category of usual
right D-modules on $\Gr_G$, via the tensor product by the 
corresponding line bundle. Let $\on{D}(\Gr_G)_\crit\mod^{G[[t]]}$
be the corresponding category of $G[[t]]$-equivariant twisted D-modules.

\medskip

The geometric Satake equivalence (see \cite{MV}) defines a functor
$$V\in \Rep(\cG)\mapsto \CF_V\in \on{D}(\Gr_G)_\crit\mod^{G[[t]]}.$$

Moreover, $\on{D}(\Gr_G)_\crit\mod^{G[[t]]}$ is endowed with a structure
of tensor category via the convolution product, such that the above
functor becomes an equivalence of tensor categories.

\medskip

To $\CM\in \hg_\crit\mod^{G[[t]]}$ and $\CF\in \on{D}(\Gr_G)_\crit\mod$
we can associate their convolution $$\CF\star \CM\in D(\hg_\crit),$$
(see \cite{FG2}, Sect. 22.5 for the corresponding definitions).

If $\CF$ is an object of $\on{D}(\Gr_G)_\crit\mod^{G[[t]]}$, then
$\CF\star \CM$ will be naturally an object of $D(\hg_\crit)^{G[[t]]}$.

\begin{thm}   \label{main, spherical}
For $\CM\in \hg_\crit\mod^{G[[t]]}$ and $V\in \Rep(\cG)$, the
convolution $\CF_V\star \CM$ is acyclic off cohomological degree $0$,
and we have a functorial isomorphism
\begin{equation} \label{main eq}
\fs_V:\CF_V\star \CM \simeq
\CV_{\wt\fZ^{\int,\reg}_{\fg}}\underset{\wt\fZ^{\int,\reg}_{\fg}}
\otimes \CM,
\end{equation}
compatible with tensor products of $\cG$-representations, i.e., for
$V,W\in \cG$ the diagrams
\begin{equation} \label{comp1}
\CD \CF_V\star (\CF_W\star \CM) @>{\fs_V}>>
\CV_{\wt\fZ^{\int,\reg}_{\fg}}\underset{\wt\fZ^{\int,\reg}_{\fg}}\otimes
(\CF_W\star \CM) \\ @V{\sim}VV @V{\on{id}_V\otimes \fs_W}VV \\
(\CF_V\star \CF_W)\star \CM @>{\fs_{V\otimes W}}>>
\CV_{\wt\fZ^{\int,\reg}_{\fg}}\underset{\wt\fZ^{\int,\reg}_{\fg}}\otimes
\CW_{\wt\fZ^{\int,\reg}_{\fg}}\underset{\wt\fZ^{\int,\reg}_{\fg}}\otimes
\CM, \endCD
\end{equation}
and
\begin{equation} \label{comp2}
\CD \CF_V\star (\CF_W\star \CM) @>{\on{id}_V\star \fs_W}>> \CF_V\star
(\CW_{\wt\fZ^{\int,\reg}_{\fg}}\underset{\wt\fZ^{\int,\reg}_{\fg}}
\otimes \CM) \\ @V{\sim}VV @V{\fs_V}VV \\ (\CF_V\star \CF_W)\star \CM
@>{\fs_{V\otimes W}}>>
\CV_{\wt\fZ^{\int,\reg}_{\fg}}\underset{\wt\fZ^{\int,\reg}_{\fg}}
\otimes
\CW_{\wt\fZ^{\int,\reg}_{\fg}}\underset{\wt\fZ^{\int,\reg}_{\fg}}
\otimes \CM, \endCD
\end{equation}
are commutative. The endomorphism on RHS of \eqref{main eq}, given by
$N_{\CV_{\wt\fZ^{\int,\reg}_{\fg}}}$, is identically equal to $0$.
\end{thm}

The last statement of \thmref{main, spherical} implies the following
statement, conjectured by A.~Beilinson (see \cite{FG1}, Conjecture
6.13):

\begin{cor}   \label{monodromy-free}
The support in $\Spec(\fZ_\fg)$ of any $G[[t]]$-integrable
$\hg_\crit$-module is contained in $\Spec(\fZ_\fg^{\on{m.f.}})$.
\end{cor}
 
\ssec{The case of differential operators on the group} \label{thm for
diff op}

Let $\fD_{G,\crit}$ denote the chiral algebra of differential
operators on $G$ at the critical level introduced in \cite{AG}.
It comes equipped homomorphisms of chiral algebras
$$\fl:\CA_{\fg,\crit}\to \fD_{G,\crit}\rightarrow \CA_{\fg,\crit}:\fr$$
whose images mutually Lie-* commute.
Let $\fD_{G,\crit,x}$ be its fiber at $x$, which we view as a
$G[[t]]$-equivariant bimodule at the critical level.
 
Using the map $\pi:G\ppart\to \Gr_G$, starting from $\CF\in
\on{D}(\Gr_G)_\kappa\mod$ (for any level $\kappa$), we produce a
chiral $\fD_{G,\crit}$-module supported at $x$ by considering
$$\Gamma(G\ppart,\pi^*(\CF)).$$ We can also view it as
$\hg_\crit$-bimodule, which is $G[[t]]$-equivariant with respect to
the $\fr$ action.
 
We have:
\begin{equation} \label{sections on D as conv}
\Gamma(G\ppart,\pi^*(\CF))\simeq \CF\star \fD_{G,\crit,x},
\end{equation}
as $\fD_{G,\crit}$-modules with respect to the action of $G$ on itself
by left translations. Hence, the corresponding isomorphism holds also
on the level of $\hg_\crit$-bimodules.  \thmref{main, spherical} then
immediately implies the following:
 
\begin{thm} \label{main, spherical, diff}
We have a canonical isomorphism of $\hg_\crit$-bimodules,
\begin{equation} \label{main eq diff}
\Gamma(G\ppart,\pi^*(\CF_V))\simeq 
\CV_{\wt\fZ^{\int,\reg}_{\fg}}\underset{\wt\fZ^{\int,\reg}_{\fg}}
\otimes \fD_{G,\crit,x},
\end{equation}
where the tensor product is taken with respect to the
$\fZ_\fg$-module structure on $\fD_{G,\crit,x}$, given by
$\fl$. These isomorphisms are compatible with tensor products of
objects of $\Rep(\cG)$. The endomorphism on the RHS of \eqref{main eq
diff}, given by $N_{\CV_{\wt\fZ^{\int,\reg}_{\fg}}}$, is zero.
\end{thm}
 
Conversely, \thmref{main, spherical, diff} implies the first two
statements of \thmref{main, spherical}:
 
\begin{proof}
 
Recall (see \cite{FG2}, Sect. 21.13) that for a $G[[t]]$-integrable
$\hg_\kappa$-module $\CM$ (at any level $\kappa$) and $\CF\in
\on{D}(\Gr_G)_\kappa\mod$ there exists a canonical isomorphism of
individual cohomologies
\begin{equation} \label{conv as semiinf}
h^i(\CF\star \CM) \simeq h^i\Bigl(
\Gamma(G\ppart,\pi^*(\CF))\underset{\fg\ppart,\fg}{\overset{\semiinf}
\otimes}\CM\Bigr),
\end{equation}
where the semi-infinite Tor is taken with respect to action of
$\hg_{\kappa'}$ on $\Gamma(G\ppart,\CF)$ given by $\fr$ (here
$\kappa'$ is the opposite level). We refer the reader to \cite{FG2} for
the precise definition of this functor.
 
\medskip
  
Applying this to $\CF=\CF_V$, and using the isomorphism of
\thmref{main, spherical, diff}, we obtain that the cohomologies of
$\CF_V\star \CM$ are isomorphic to those of
\begin{equation} \label{co 1}
\Bigl(\CV_{\wt\fZ^{\int,\reg}_{\fg}}
\underset{\wt\fZ^{\int,\reg}_{\fg}}\otimes \fD_{G,\crit,x}\Bigr)
\underset{\fg\ppart,\fg}{\overset{\semiinf}\otimes}\CM.
\end{equation} 

However, since $\CV_{\wt\fZ^{\int,\reg}_\fg}$ is locally free over
$\wt\fZ^{\int,\reg}_\fg$, we obtain that the complex \eqref{co 1} is
quasi-isomorphic to
\begin{equation} \label{co 2}
\CV_{\wt\fZ^{\int,\reg}_{\fg}}
\underset{\wt\fZ^{\int,\reg}_{\fg}}\otimes \Bigl(\fD_{G,\crit,x}
\underset{\fg\ppart,\fg}{\overset{\semiinf}\otimes}\CM\Bigr).
\end{equation}

However, by \cite{FG2}, Corollary 21.14, at any level $\kappa$, we
have a quasi-isomorphism of $\hg_\kappa$-modules
\begin{equation} \label{semi-jective}
\fD_{G,\kappa,x}\underset{\fg\ppart,\fg}{\overset{\semiinf}\otimes}
\CM\simeq \CM.
\end{equation}
 
Combining \eqref{co 2} and \eqref{semi-jective} we obtain that
$\CF\star \CM$ has the same cohomologies as
$\CV_{\wt\fZ^{\int,\reg}_{\fg}}
\underset{\wt\fZ^{\int,\reg}_{\fg}}\otimes \CM $, as contended.

\end{proof}

\section{Interpretation of convolution as fusion}      \label{sect 2}

Our strategy of proof of \thmref{main, spherical} is to show that both
sides of \eqref{main eq} represent the same functor on the category
$\hg_\crit\mod$. This functor has to do with the notion of {\it
chiral pairing} or {\it fusion} of modules over a chiral algebra,
which we shall presently define. Having introduced this
functor, in the rest of the section we consider various properties of 
fusion and its connection to the convolution functors.

\ssec{Chiral pairings}
 
Let $\CA$ be a chiral algebra on a smooth curve $X$, and $\CM_1,\CM_2$
and $\CN$ be chiral $\CA$-modules. As usual, we shall denote by
$\Delta$ the embedding of $X$ into $X^n$ as the main diagonal, and by
$j$ the embedding of the complement of the diagonal divisor.
 
A chiral pairing $\{\CM_1,\CM_2\}\to \CN$ is by definition a map of
D-modules on $X\times X$:
$$\phi:j_*j^*(\CM_1\boxtimes \CM_2)\to \Delta_!(\CN),$$
which is compatible with the $\CA$ action in the following sense:
 
We need that the sum of the three morphisms
$$j_*j^*(\CA\boxtimes \CM_1\boxtimes \CM_2)\to \Delta_!(\CN)$$ is
zero, where the first morphism is
$$j_*j^*(\CA\boxtimes \CM_1\boxtimes \CM_2)
\overset{\CA\text{-action on } \CM_1}\longrightarrow
\Delta_{x_1=x_2}{}_!(j_*j^*(\CM_1\boxtimes\CM_2))
\overset{\phi}\to \Delta_!(\CN),$$ the second morphism is the negative of
$$j_*j^*(\CA\boxtimes \CM_1\boxtimes \CM_2)
\overset{\CA\text{-action on } \CM_2}\longrightarrow
\Delta_{x_1=x_3}{}_!(j_*j^*(\CM_1\boxtimes\CM_2))
\overset{\phi}\to \Delta_!(\CN),$$ and the third morphism is
$$j_*j^*(\CA\boxtimes \CM_1\boxtimes \CM_2)\overset{\phi}\to
\Delta_{x_1=x_3}{}_!(j_*j^*(\CA\boxtimes \CN))
\overset{\CA\text{-action on } \CN}\longrightarrow \Delta_!(\CN).$$

Chiral pairings evidently form a functor $\CA\mod^o\times
\CA\mod^o\times \CA\mod\to \on{Vect}$.

\medskip

\noindent{\it Remark.}  The notion of chiral pairing was introduced in
\cite{CHA}. As was mentioned in the introduction, in the case when
$\CA$ is obtained from a conformal vertex algebra ${\mathsf V}$, and
the modules $\CM_1,\CM_2$ and $\CN$ correspond to some ${\mathsf
V}$-modules $M_1,M_2$ and $N$, the set of chiral pairings
$\{\CM_1,\CM_2\}\to \CN$ is in a natural bijection with the set of
intertwining operators $M_1\otimes M_2\to N$ in the sense of
\cite{FHL}, such that in the corresponding fields all powers of the
formal variable $z$ are integral.

\medskip

Let us regard $\CA$ as a chiral module over itself. Then for any 
other chiral $\CA$-module $\CM$, we have a canonical chiral pairing
$\{\CA,\CM\} \to \CM$ given by the action of $\CA$ on $\CM$. The
following result is obtained directly from the definitions:
 
\begin{lem} \label{A is unit}
The canonical chiral pairing $\{\CA,\CM\}\to \CM$ establishes a
bijection between the set of maps of chiral modules $\CM\to \CN$ and
the set of chiral pairings $\{\CA,\CM\}\to \CN$.
\end{lem}
 
\ssec{A canonical chiral pairing}
 
Let $\kappa$ be any level, and let $\CA_{\fg,\kappa}$ be the
corresponding chiral algebra on $X$.
 
Consider the relative version $\Gr_{G,X}$ of the affine Grassmannian
over $X$.  Let $\CF_X$ be a (right, $\kappa$-twisted) D-module on
$\Gr_{G,X}$.  In what follows we will assume that $\CF_X$ is
torsion-free with respect to $X$.  By a slight abuse of notation, we
will denote by $\Gamma(\Gr_{G,X},\CF_X)$ the quasi-coherent direct
image of $\CF_X$ on $X$; this is a chiral $\CA_{\fg,\kappa}$-module.
By replacing the subscript "$X$" by either "$X-x$" or "$x$" we will
denote the fibers of the above objects over $X-x$ and $x$
respectively.

\begin{propconstr}   \label{const of basic pairing}
Given an object $\CM\in \hg_\kappa\mod^{G[[t]]}$ there exists a canonical
chiral pairing
\begin{equation} \label{basic pairing for modules}
\{\Gamma(\Gr_{G,X},\CF_X),\CM\}\to h^0(\CF_x\star \CM).
\end{equation}
\end{propconstr}

The rest of this subsection is devoted to the proof of this proposition. 

\medskip

Let $\on{Jets}^{\mer}(G)_X$ be the D-ind scheme over $X$ of meromorphic
jets into $G$. This is a relative version of the ind-scheme
$G\ppart$. By definition, the category
$\on{D}(\on{Jets}^{\mer}(G)_X)_\kappa\mod$ of (right $\kappa$-twisted)
D-modules on $\on{Jets}^{\mer}(G)_X$ is equivalent to that of chiral
$\fD_{G,\kappa}$-modules.
 
\medskip
 
We have a natural projection $\pi:\on{Jets}^{\mer}(G)_X\to \Gr_{G,X}$
and by considering the quasi-coherent pull-back, we obtain a functor
$$\pi^*:\on{D}(\Gr_{G,X})_\kappa\mod\to
\on{D}(\on{Jets}^{\mer}(G)_X)_\kappa\mod.$$ For $\CF_X\in
\on{D}(\Gr_{G,X})_\kappa\mod$, we will denote by
$\Gamma(\on{Jets}^{\mer}(G)_X,\pi^*(\CF_X))$ the resulting
$\fD_{G,\kappa}$-module.
 
One reconstructs $\Gamma(\Gr_{G,X},\CF_X)$ as a subset of
$\Gamma(\on{Jets}^{\mer}(G)_X,\pi^*(\CF_X))$ as follows: this is the
D-submodule consisting of sections that *-commute with
$\fr(\CA_{\fg,\kappa'})$, cf. \cite{AG} or \cite{FG1}, Theorem 2.5.
 
\medskip
 
For $\CF_X\in \on{D}(\Gr_{G,X})_\kappa\mod$ consider the action map
\begin{equation} \label{act of diff op}
j_*j^*\Bigl(\Gamma(\on{Jets}^{\mer}(G)_X,\pi^*(\CF_X))\boxtimes
\fD_{G,\kappa}\Bigr)\to
\Delta_!\Bigl(\Gamma(\on{Jets}^{\mer}(G)_X,\pi^*(\CF_X))\Bigr).
\end{equation} 
 
This is a $\fD_{G,\kappa}$-chiral pairing, and hence, in particular, a
chiral pairing with respect to the action of $\CA_{\fg,\kappa}$ via
$\fl$.  In particular, we obtain a chiral $\CA_{\fg,\kappa}$-pairing
\begin{equation}  \label{basic pairing univ}
j_*j^*(\Gamma(\Gr_{G,X},\CF_X)\boxtimes \fD_{G,\kappa})\to 
\Delta_!(\Gamma(\on{Jets}^\mer(G)_X,\pi^*(\CF_X)).
\end{equation}
 
This map commutes with the right $\CA_{\fg,\kappa'}$-action in the
sense that the map
\begin{align*}
&\wt{ j}_*\wt{j}^*(\Gamma(\Gr_{G,X},\CF_X)\boxtimes 
\fD_{G,\kappa}\boxtimes \CA_{\fg,\kappa'})
\to \Delta_{x_1=x_2}{}_!\Bigl(\Gamma(\on{Jets}^{\mer}(G)_X,\pi^*(\CF_X))
\boxtimes \CA_{\fg,\kappa'}\Bigr) \to \\
&\to \Delta_!\Bigl(\Gamma(\on{Jets}^{\mer}(G)_X,\pi^*(\CF_X))\Bigr)
\end{align*}
coincides with
\begin{align*}
&\wt{ j}_*\wt{j}^*(\Gamma(\Gr_{G,X},\CF_X) \boxtimes
\fD_{G,\kappa}\boxtimes \CA_{\fg,\kappa'}) \to
\Delta_{x_2=x_3}{}_!\Bigl(\Gamma(\Gr_{G,X},\CF_X)\boxtimes
\fD_{G,\kappa}\Bigr) \to \\ &\to
\Delta_!\Bigl(\Gamma(\on{Jets}^{\mer}(G)_X,\pi^*(\CF_X))\Bigr),
\end{align*}
where $\wt{j}$ denotes the embedding of the complement of the union of
divisors corresponding to $x_1=x_2$ and $x_2=x_3$.
  
\medskip
  
Let us restrict both sides of \eqref{basic pairing univ} to $X\times
x\subset X\times X$.  We obtain a chiral pairing
\begin{equation} \label{univ pairing at point}
j_x{}_*j_x^*\Bigl(\Gamma(\Gr_{G,X-x},\CF_{X-x})\Bigr)\otimes
\fD_{G,\kappa,x}\to i_x{}_!  \Bigl(\Gamma(G\ppart,\pi^*(\CF_x))\Bigr),
\end{equation}
where we view $\CF_x$ as a $\kappa$-twisted D-module over the affine
Grassmannian, and $i_x$ (resp., $j_x$) denotes the embedding of $x$
into $X$ (resp., the embedding of this complement).

\medskip

By the above, the map in \eqref{univ pairing at point} commutes with
the action of $\hg_{\kappa'}$ via $\fr$. (On the LHS this action
affects only the $\fD_{G,\kappa,x}$ multiple.)

Let now $\CM$ be a $G[[t]]$-integrable $\hg_\kappa$-module. By
considering the complex, computing the semi-infinite Tor,
$$\fC^\semiinf(\fg\ppart;\fg,?\otimes \CM)$$ against the two sides of
\eqref{univ pairing at point}, considered as $\hg_{\kappa'}$-modules,
we obtain a chiral pairing of complexes
of $\CA_{\fg,\kappa}$-modules:
$$\{\Gamma(\Gr_{G,X},\CF_X),
\fC^\semiinf\Bigl(\fg\ppart;\fg,\fD_{G,\kappa,x}\otimes \CM\Bigr)\}\to
\fC^\semiinf\Bigl(\fg\ppart;\fg,\Gamma(G\ppart,\pi^*(\CF_x))\otimes
\CM\Bigr).$$

By passing to the $0$th cohomology, and taking into account
\eqref{conv as semiinf} we obtain the chiral pairing of \eqref{basic
pairing for modules}.

\ssec{Chiral pairing with a module of global sections is represented 
by the convolution}

Let us now specialize to the case when $\kappa$ is integral
non-positive (see \secref{recol} for the definition of what this
means). For $V\in \Rep(\cG)$, let $\CF_{V,X}$ be the corresponding
object of $\on{D}(\Gr_{G,X})_\kappa\mod$. Then from \eqref{basic pairing
for modules} we obtain a canonical chiral pairing
\begin{equation} \label{pairing with sph}
\{\Gamma(\Gr_{G,X},\CF_{V,X}),\CM\}\to h^0(\CF_V\star \CM).
\end{equation}
 
In \secref{proof of rep} we will prove the following theorem:

\begin{thm}  \label{representability, spherical} \hfill

\smallskip

\noindent{\em(1)} For any $\CM\in \hg_\kappa\mod^{G[[t]]}$ and $V\in
\Rep(\cG)$, the convolution $\CF_V\star \CM$ is acyclic away from
cohomological degree $0$.

\smallskip

\noindent{\em(2)}
The functor on $\hg_\kappa\mod$ that sends $\CN$ to the set 
of chiral parings $\{\Gamma(\Gr_{G,X},\CF_{V,X}),\CM\}\to \CN$ is
representable by $\CF_V\star \CM$.

\end{thm}

Our proof of \thmref{main, spherical} will be independent of
\thmref{representability, spherical}. However, it is useful to 
keep \thmref{representability, spherical} in mind as it gives us 
an important insight into the connection between fusion and convolution.

\ssec{}  \label{sect pairing with sph}

We will need the following result about the associativity
property of the map given by \eqref{pairing with sph}.

For any two objects $V,W\in \Rep(\cG)$ there exists a canonical chiral
pairing
$$\{\Gamma(\Gr_{G,X},\CF_{V,X}),\Gamma(\Gr_{G,X},\CF_{W,X})\}\to
\Gamma(\Gr_{G,X},\CF_{V\otimes W,X}).$$ (Its construction will be
recalled in the sequel.) Let $\CM$ be an object of $\CM\in
\hg_\kappa\mod^{G[[t]]}$, such that $\CF_{V'}\star \CM$ is acyclic
away from cohomological degree $0$ for any $V'\in \Rep(\cG)$ (we
choose not to rely here on \thmref{representability, spherical}(1),
which says that the latter assumption is satisfied automatically).

Then there are three maps
$$j_*j^*\Bigl(\Gamma(\Gr_{G,X},\CF_{V,X})\boxtimes
\Gamma(\Gr_{G,X},\CF_{W,X})\boxtimes \CM\Bigr)\to
\Delta_!(\CF_{V\otimes W}\star \CM):$$
 
The first one is the composition
\begin{align*}
& j_*j^*\Bigl(\Gamma(\Gr_{G,X},\CF_{V,X})\boxtimes
\Gamma(\Gr_{G,X},\CF_{W,X})\boxtimes \CM\Bigr)\to \\ &\to
\Delta_{x_2=x_3}{}_!\Bigl(j_*j^*(\Gamma(\Gr_{G,X},\CF_{V,X})\boxtimes
(\CF_W\star \CM))\Bigr) \to \Delta_!(\CF_V\star (\CF_W\star
\CM))\simeq \Delta_!(\CF_{V\otimes W}\star \CM).
\end{align*}
The second map is the negative of a similar map with the roles of the
first and the second factor swapped. The third map is the composition
\begin{align*}
&j_*j^*\Bigl(\Gamma(\Gr_{G,X},\CF_{V,X})\boxtimes
\Gamma(\Gr_{G,X},\CF_{W,X})\boxtimes \CM\Bigr)\to \\ &\to
\Delta_{x_1=x_2}{}_!\Bigl(j_*j^*(\Gamma(\Gr_{G,X},\CF_{V\otimes W,X})
\boxtimes \CM)\Bigr)\to \Delta_!(\CF_{V\otimes W}\star \CM).
\end{align*}

\begin{prop}  \label{comp with ten products}
The sum of the above three maps is $0$.
\end{prop}
 
\ssec{Proof of \propref{comp with ten products}}  \label{Jets n}

By the construction of the chiral pairing \eqref{pairing with sph}, it
is sufficient to consider the universal case, i.e., when $\CM\simeq
\fD_{G,\kappa}$ as a chiral $\CA_{\fg,\kappa}$-module under $\fl$.

Let us rewrite the map \eqref{act of diff op} in more geometric terms.
For any $n$, let $\Gr_{G,X^n}$ be the Beilinson-Drinfeld affine
Grassmannian over $X^n$. I.e., this is the ind-scheme classifying the
data of
$$(x_1,...,x_n,\CP_G,\beta),$$ where $x_1,...,x_n$ is an $n$-tuple of
points on $X$, $\CP_G$ is a principal $G$-bundle on $X$, and $\beta$
is a trivialization of $\CP_G$ on $X-\{x_1,...,x_n\}$.

For a decomposition $n=n_1+n_2$ there exists a natural map
$${\bf 1}_{n_1,n_2}:X^{n_1}\times \Gr_{G,X^{n_2}}\to \Gr_{G,X^n},$$
which sends $((x_{1},...x_{n_1}), (x_{n_1+1},...,x_{n},\CP_G,\beta))
\mapsto (x_1,...,x_{n_1},x_{n_1+1},...,x_n,\CP_G,\beta)$.

\medskip

Consider now the ind-scheme $\on{Jets}^\mer(G)_{X^n}$, fibered over
$\Gr_{G,X^n}$, where in addition to the data
$(x_1,...,x_n,\CP_G,\beta)$ we have that of a trivialization of
$\CP_G$ on a formal neighborhood of $\underset{i}\cup\, x_i\subset X$.
Let us denote by $\pi$ the corresponding projection. This ind-scheme
also carries a flat connection along $X^n$.

The scheme $\on{Jets}^\mer(G)_{X^n}$ has the usual factorization
pattern:
\begin{equation} \label{factorization}
j^*\Bigl(\on{Jets}^\mer(G)_{X^n}\Bigr)\simeq
j^*\Bigl(\on{Jets}^\mer(G)_X^{\times n}\Bigr) \text{ and }
\Delta^*\Bigl(\on{Jets}^\mer(G)_{X^n}\Bigr)\simeq
\on{Jets}^{\mer}(G)_X.
 \end{equation}

\medskip

For two objects $V,W\in \Rep(\cG)$ we consider the D-module
$j^*(\CF_{V,X}\boxtimes \CF_{W,X})$ on $\Gr_{G,X^2}|_{X^2-X}$. As is
well-known, its Goresky-MacPherson extension
$j_{!*}\left(j^*(\CF_{V,X}\boxtimes \CF_{W,X})\right)$ onto
$\Gr_{G,X^2}$ has the property that
 $$\Delta^!(j_{!*}\left(j^*(\CF_{V,X}\boxtimes \CF_{W,X})\right))[1]\simeq
 \CF_{V\otimes W,X}.$$
 
Consider the (right $\kappa$-twisted) D-module
$$({\bf 1}_{1,2})_!\Bigl(\omega_X\boxtimes j_{!*}\left(j^*(\CF_{V,X}\boxtimes
\CF_{W,X})\right)\Bigr)$$ on $\Gr_{G,X^3}$. We have
three maps

$$j_*j^*\Biggl(({\bf 1}_{1,2})_!\Bigl(j_{!*}\left(\omega_X\boxtimes
j^*(\CF_{V,X}\boxtimes \CF_{W,X})\right)\Bigr)\Biggr)\to
\Delta_!(\CF_{V\otimes W,X}),$$ corresponding to the three diagonals
in $X^3$, whose sum is equal to $0$.

Applying the pull-back by means of $\pi$ to the two sides of the above
formula to $\on{Jets}^\mer(G)_{X^3}$, followed by the direct image
onto $X^3$, we obtain three maps
\begin{align*} \label{triple map downstairs}
&j_*j^*\Bigl(\fD_{G,\kappa}\boxtimes
  \Gamma(\on{Jets}^\mer(G)_X,\pi^*(\CF_{V,X}))\boxtimes
  \Gamma(\on{Jets}^\mer(G)_X,\pi^*(\CF_{W,X}))\Bigr)\to \\ &\to
  \Delta_!(\Gamma(\on{Jets}^\mer(G)_X,\pi^*(\CF_{V\otimes W})).
\end{align*}

Finally, note that we have a natural inclusion
\begin{align*}
&j_*j^*\Bigl(\fD_{G,\kappa}\boxtimes
\Gamma(\Gr_{G,X},\CF_{V,X})\boxtimes
\Gamma(\Gr_{G,X},\CF_{W,X})\Bigr)\hookrightarrow \\
&j_*j^*\Bigl(\fD_{G,\kappa}\boxtimes
\Gamma(\on{Jets}^\mer(G)_X,\pi^*(\CF_{V,X}))\boxtimes
\Gamma(\on{Jets}^\mer(G)_X,\pi^*(\CF_{W,X}))\Bigr).
\end{align*}

By composing this inclusion with the three maps above
we obtain three maps
\begin{equation} \label{same three maps}
j_*j^*\Bigl(\fD_{G,\kappa}\boxtimes
\Gamma(\Gr_{G,X},\CF_{V,X})\boxtimes
\Gamma(\Gr_{G,X},\CF_{W,X})\Bigr)\to
\Delta_!\Bigl(\Gamma(\on{Jets}^\mer(G)_X,
\pi^*(\CF_{V\otimes W,X}))\Bigr).
\end{equation}

\begin{lem}
The three maps of \eqref{same three maps} coincide with those that
appear in the statement of \propref{comp with ten products}. 
\end{lem}

Clearly, this implies the statement of the proposition. To prove the
above lemma we will consider a slightly more general framework.

\ssec{}  \label{sect action as fusion}

Let $\CF_1$ be a ($\kappa$-twisted, right) D-module on $G\ppart/K$,
where $K\subset G[[t]]$ is a subgroup of finite codimension.  Let $\CF'_1$ 
be the pull-back of $\CF_1$ to $G\ppart$.

Consider the ind-scheme $\Gr_{G,X;K}$, classifying the data
$(x_1,\CP_G,\beta,\alpha)$, where $(x_1,x,\CP_G,\beta)$ is a point of
$\Gr_{G,X\times x}$, and $\alpha$ is a reduction of $\CP_G|_{\D_x}$ to
$K$. We have the isomorphisms
$$\Gr_{G,X-x;K}\simeq \Gr_{G,X-x}\times G\ppart/K \text{ and }
\Gr_{G,x;K}\simeq G\ppart/K.$$ We also have a natural projection
$\pi_K:\on{Jets}^\mer(G)_{X\times x}\to \Gr_{G,X;K}$, where
$\on{Jets}^\mer(G)_{X\times x}$ denotes the preimage of $X\times
x\subset X^2$ in $\on{Jets}^\mer(G)_{X^2}$.

\medskip

Consider the direct image of $\omega_X\boxtimes \CF_1$ under the
tautological embedding ${\bf 1}_{1,1}:X\times G\ppart/K\to
\Gr_{G,X;K}$.  It gives rise to a map of twisted D-modules on
$\Gr_{G,X;K}$:
\begin{equation}  \label{basic map down}
j_x{}_*j_x^*(\delta_{1,\Gr_{G,X}}\boxtimes \CF_1)\to i_x{}_!(\CF_1).
\end{equation}

Lifting the two sides of \eqref{basic map down} by means of $\pi_K$
and taking the (quasi-coherent direct image onto $X\simeq X\times x$,
we obtain a map
\begin{equation} \label{basic map up}
j_x{}_*j_x^*(\fD_{G,\kappa}) \otimes \Gamma(G\ppart,\CF'_1)\to
i_x{}_!\Bigl(\Gamma(G\ppart,\CF'_1)\Bigr).
\end{equation}

The following assertion is built into the interpretation of twisted
D-modules on $G(\ppart$ as chiral $\fD_{G,\kappa}$-modules:

\begin{prop} \label{action as fusion}
The map of \eqref{basic map up} coincides with the chiral action map
for $\Gamma(G\ppart,\CF'_1)$ considered as a chiral
$\fD_{G,\kappa}$-module supported at $x$.
\end{prop}

In the appendix (\secref{app A}) we shall give a proof of this result for the sake of
completeness.

\medskip

Assume now that the D-module $\CF_1$ on $G\ppart/K$ is
$G[[t]]$-equivariant with respect to the left action of $G[[t]]$ on
$G\ppart/K$.  Let now $\CF_X$ be a ($\kappa$-twisted, right) D-module
on $\Gr_{G,X}$.  In this case we can form the convolution
$\CF_X\underset{G[[t]]}\star \CF_1$, which will be an object of the
derived category of ($\kappa$-twisted, right) D-modules on
$\Gr_{G,X;K}$.  We have:
$$\CF_X\underset{G[[t]]}\star \CF_1|_{\Gr_{G,X-x;K}}\simeq
\CF_{X-x}\boxtimes \CF_1 \text{ and } \CF_x\underset{G[[t]]}\star
\CF_1|_{\Gr_{G,x;K}}\simeq \CF_{x}\underset{G[[t]]}\star \CF_1.$$

In particular, we obtain a map of D-modules on $\Gr_{G,X;K}$:
\begin{equation}  \label{F map down}
j_x{}_*j_x^*\Bigl(\CF_{X-x}\boxtimes \CF_1\Bigr)\to
i_x{}_!\Bigl(h^0(\CF_{x}\underset{G[[t]]}\star \CF_1)\Bigr).
\end{equation}
Pulling the two sides of \eqref{F map down} under the map
$\pi_K:\on{Jets}^\mer(G)_{X\times x}\to \Gr_{G,X;K}$, and taking the
(quasi-coherent) direct image onto $X$, we obtain a map of D-modules
on $X$:
\begin{align*}  
&j_x{}_*j_x^*\Bigl(\Gamma(\on{Jets}^\mer(G)_X,\pi^*(\CF_X))\Bigr)
\otimes \Gamma(G\ppart,\CF'_1) \to
i_x{}_!\Bigl(\Gamma(G\ppart,h^0(\CF_x\star \CF'_1))\Bigr)\simeq \\
&\simeq i_x{}_!\Bigl(h^0\Bigl(\CF_x\star
\Gamma(G\ppart,\CF'_1)\Bigr)\Bigr).
\end{align*}

\begin{lem} \label{two pair}
The above map coincides with \eqref{pairing with sph} for
$\CM=\Gamma(G\ppart,\CF'_1)$.
\end{lem}

It is easy to see that this lemma, combined with a version of
\propref{action as fusion}, where the point $x$ varies on $X$, imply
the required property of the three maps of \eqref{same three maps}.
 
\begin{proof} (of \lemref{two pair}) 
 
Consider first the case when $K=G[[t]]$ and
$\CF_1=\delta_{1,\Gr_G}$. Then we are dealing with a map
$$j_x{}_*j_x^*\Bigl(\Gamma(\on{Jets}^\mer(G)_X,\pi^*(\CF_X))\Bigr)\otimes
\fD_{G,\kappa,x}\to i_x{}_!\Bigl(\Gamma(G\ppart,\pi^*(\CF_x))\Bigr),$$
and we claim that it coincides with the map, obtained by restriction
to $X\times x$ from the map of \eqref{act of diff op}. This follows
from a version of \propref{action as fusion} for $K=G[[t]]$, when
instead of the fixed point $x$ we consider a morphism of D-modules on
$X\times X$, and the map \eqref{act of diff op} itself is described as
a Cousin map corresponding to $({\bf 1}_{1,1})_!(\CF_X)$, viewed as
D-module on $\Gr_{G,X^2}$.

\medskip

To pass to the general case, consider the category of $\kappa$-twisted
D-modules on $\on{Jets}^\mer(G)_{X,x}$ as acted on by the group
ind-scheme $G\ppart$ on the right, which is of Harish-Chandra type
with respect to the extension $\hg_{\kappa'}$ of $\fg\ppart$. Then the
convolution action of $\CF_1$ on
$$j_*{}_*j_x^*(\pi^*(\CF_{X-x})\boxtimes \pi^*(\delta_{1,\Gr_G}))
\text{ (resp., } i_x{}_!(\CF_x))$$ equals
$$j_*{}_*j_x^*(\pi^*(\CF_{X-x})\boxtimes \CF'_1) \text{ and }
i_x{}_!(\CF_x\star \CF'_1),$$
respectively.

Consider also the category of D-modules on $X$, endowed with an action
of $\hg_{\kappa'}$ as a category with a Harish-Chandra action of
$G\ppart$.  Then the convolution action of $\CF_1$ on a module of the
type $j_x{}_*j_x^*(\CM)\otimes \fD_{G,\kappa,x}$ (resp.,
$i_x{}_!(\CN)$) (where $\CM$ is a D-module, and $\CN$ is a
representation of $\hg_{\kappa'}$) equals
$$j_x{}_*j_x^*(\CM)\otimes \Gamma(G\ppart,\CF'_1) \text{ and }
i_x{}_!\Bigl(\fC^\semiinf(\fg\ppart;\fg,\CN\otimes
\Gamma(G\ppart,\CF'_1))\Bigr),$$ respectively.

Note now that the (quasi-coherent) direct image is a functor between
the above two categories that respects the $G\ppart$-actions. Hence,
the assertion of the lemma follows from the case when $K=G[[t]]$ and
$\CF_1=\delta_{1,\Gr_G}$ by the functoriality of the convolution.

\end{proof}

\section{Proof of \thmref{main, spherical}}  \label{sect 3}

Having established some preliminary results in the previous section,
we are now ready to prove \thmref{main, spherical}. The plan of the
proof is as follows: we will first show that the chiral maps
$\{\CM_1,\CM \} \to \CN$ in the special case when $\CM_1 =
\CV_{\fz_\fg}\underset{\fz_\fg}\otimes \CA_{\fg,\crit}$ are the same
as ordinary homomorphisms of $\hg_\crit$-modules $\CM' \to \CN$ where
$\CM'\simeq
\CV_{\wt\fZ_\fg^{\int,\reg}}\underset{\wt\fZ_\fg^{\int,\reg}}\otimes
\CM$.  Roughly speaking, this means that we can ``swap'', under the chiral
product, the tensor product with $\CV_{\fz_\fg}$ from
$\CA_{\fg,\crit}$ to $\CM$. This is the content of \thmref{repr of
pair from ten product}.

Once we have that, we obtain, for each $V \in \Rep(\cG)$, a map in one
direction in \eqref{main eq}, provided that $\CF_V \star \CM$ is
acyclic away from cohomological dimension $0$. But we do know this for
$\CM = \fD_{G,\crit,x}$, according to \thmref{main, spherical, diff}. 
Moreover, according to \secref{thm for diff op}, it
is sufficient to prove the statement of \thmref{main, spherical} just
in this special case. Hence all we need is to check that the maps we
have constructed in this case are indeed isomorphisms and that they
satisfy the properties listed in \thmref{main, spherical}. We first
check that the maps in question do satisfy the required properties; then
we show that this already implies that the above maps are necessarily 
isomorphisms, thus completing the proof.

\ssec{Fusion at the critical level}   \label{crit fusion}

In this section we fix $\kappa=\kappa_\crit$, and study the map of
\eqref{pairing with sph} in this case. We will use the following
crucial result of \cite{BD}, Sects. 5.5-5.6: there is a canonical
isomorphism of chiral $\CA_{\fg,\crit}$-modules,
\begin{equation}  \label{BD isom}
\Gamma(\Gr_{G,X},\CF_{V,X})\simeq
\CV_{\fz_\fg}\underset{\fz_\fg}\otimes \CA_{\fg,\crit},
\end{equation}
where the tensor product makes sense since $\fz_\fg$ is the center of
$\CA_{\fg,\crit}$.

Given a chiral $\CA_{\fg,\crit}$-module $\CM$, let us study the
functor on the category of chiral $\CA_{\fg,\crit}$-modules that
assigns to a chiral $\CA_{\fg,\crit}$-module $\CN$ the set of chiral
pairings
\begin{equation} \label{pairings from ten product}
\{\CV_{\fz_\fg}\underset{\fz_\fg}\otimes \CA_{\fg,\crit},\CM\}\to \CN.
\end{equation}

\begin{thm}  \label{repr of pair from ten product}
Let $\CM$ be concentrated at $x\in X$, and assume that, when viewed as
a module over $\fZ_\fg$, it is supported on
$\wt\fZ^{\int,\reg}_{\fg}$. Then chiral pairings \eqref{pairings
from ten product} for $\CN\in \CA_{\fg,\crit}\mod_x$ are in bijection
with maps of $\hg_\crit$-modules
$$\Bigl(\CV_{\wt\fZ^{\int,\reg}_{\fg}}
\underset{\wt\fZ^{\int,\reg}_{\fg}}\otimes
\CM\Bigr)_{N_{\CV_{\wt\fZ^{\int,\reg}_{\fg}}}}\to \CN.$$
\end{thm} 

In the statement of the theorem the subscript
$N_{\CV_{\wt\fZ^{\int,\reg}_{\fg}}}$ stands for the coinvariants
with respect to the action of this nilpotent operator. Before giving a
proof of this theorem, which will occupy the next few subsections, let
us note that \thmref{repr of pair from ten product} implies the
existence of the map in one direction in \thmref{main, spherical}:

Indeed, the chiral pairing \eqref{pairing with sph} gives rise to a
map
\begin{equation} \label{desired map almost}
\Bigl(\CV_{\wt\fZ^{\int,\reg}_{\fg}}
\underset{\wt\fZ^{\int,\reg}_{\fg}}\otimes
\CM\Bigr)_{N_{\CV_{\wt\fZ^{\int,\reg}_{\fg}}}} \to h^0(\CF_V\star
\CM).
\end{equation}

Note that if we assume \thmref{representability, spherical}, then we
obtain that the map of \eqref{desired map almost} is an isomorphism
(and moreover, that $h^0(\CF_V\star \CM) = \CF_V\star \CM$). However,
this would still not be enough to prove \thmref{main, spherical}: we
would also need to show that ${N_{\CV_{\wt\fZ^{\int,\reg}_{\fg}}}}$
acts trivially on the left-hand side, which requires a separate
argument. The proof of \thmref{main, spherical} that we give below
will be independent of \thmref{representability, spherical}.

\ssec{Fusion over commutative chiral algebras}

To prove \thmref{repr of pair from ten product} we need to study chiral
pairings between modules over commutative chiral algebras.

\medskip

Recall the set-up of the proof of \propconstrref{extension of bundle
over center}, i.e., let $\CB$ be a commutative chiral algebra on $X$,
and let $\CV_{X-x}$ be a *-commutative module over $\CB|_{X-x}$,
extended to the entire $X$ such that the connection has a pole of
order $\leq 1$ and the residue $\Res(\nabla,\CV_X)\in \End(\CV_x)$ is
nilpotent.

\begin{prop}  \label{commutative fusion}
Let $\CM$ be a $\CB_x$-module, viewed as *-commutative $\CB$-module,
supported at $x$. Let $\CN$ be another chiral $\CB$-module, supported
at $x$.  Then chiral pairings
\begin{equation} \label{com pair one}
\{\CV_{X-x},\CM\}\to \CN
\end{equation}
are in bijection with maps of $\hat\CB_x$-modules 
\begin{equation} \label{com pair two}
\Bigl(\CV_x\underset{\CB_x}\otimes \CM\Bigr)_{\Res(\nabla,\CV_X)}\to
\CN,
\end{equation}
where $\hat\CB_x$ denotes the topological commutative algebra,
corresponding to $\CB$ at $x$.
\end{prop}

\begin{proof}

First, from the definition of chiral pairings it is easy to see that
any pairing $\{\CV_{X-x},\CM\}\to \CN$ factors through $\CN'\subset
\CN$, where $\CN'$ is the maximal submodule of $\CN$ on which
$\hat\CB_x$ acts via $\hat\CB_x\twoheadrightarrow \CB_x$. In other
words, we can assume that $\CN$ is also *-commutative.

Secondly, for a fixed $\CN$, both \eqref{com pair one} and \eqref{com
pair two}, regarded as contravariant functors with respect to $\CM$
are both right exact and commute with direct limits.  Hence, we can
replace $\CM$ by $\CB_x$.

\medskip

Let us denote by $B$ the fiber $\CB_x$ and by $T_0$ the module over
$B[[t]]$, corresponding to $\CV_X$. Let $T$ be the localization of
$T_0$ with respect to $t$.

Then chiral pairings $\{\CV_{X-x},\CB_x\}\to \CN$ are in bijection
with $B[[t]]$-linear maps
\begin{equation} \label{all poles}
\phi:T \to \CN\underset{B}\otimes B\ppart dt/B[[t]]dt,
\end{equation}
that respect that action of $\partial_t$.

Given such $\phi$, the finite rank assumption on $T_0$ implies
that for some integer $k$ the composition
$$t^k\cdot T_0\hookrightarrow T\overset{\phi}\to 
\CN\underset{B}\otimes B\ppart/B[[t]]$$
is zero. However, the nilpotency assumption on $t\partial_t$ acting
on $T_0/t\cdot T_0$ implies that for any $k>0$ the map
$$\partial_t:\Bigl(t^k\cdot T_0\Bigr)\to \Bigl(t^{k-1}\cdot
T_0\Bigr)$$ is surjective.

Therefore, $\phi$ in fact annihilates $T_0$. Restricting $\phi$ to
$t^{-1}\cdot T_0$ we obtain a map
\begin{equation} \label{pole one}
t^{-1}\cdot T_0/T_0\to (t^{-1}\BC[[t]]/\BC[[t]])\otimes \CN,
\end{equation}
which is zero on the image of $\partial_t(T_0)$. Hence, we obtain
a map 
$$\on{coker}\Bigl(t\partial_t:T_0/t\cdot T_0\to T_0/t\cdot
T_0\Bigr)\to \CN,$$ as desired.

\medskip

Conversely, starting from a map as in\eqref{pole one}, it is easy to
see that it uniquely extends to a map as in \eqref{all poles}.

\end{proof}

Let us denote by $can_{\CV}$ the resulting canonical map
$$j_*j^*(\CV_{X-x}\boxtimes \CM)\to 
\Delta_!(\CV_x\underset{\CB_x}\otimes \CM)_{\Res(\nabla,\CV_X)}.$$
For the proof of \thmref{main, spherical}, we will need the following 
additional property of this map.

Let $(\CV_{X-x},\CV_X)$ and $(\CW_{X-x},\CW_X)$ be two pairs of
$\CB$-modules, satisfying the assumptions of \propref{commutative
fusion}. Note that the tensor product
$(\CV_{X-x}\underset{\CB_{X-x}}\otimes \CW_{X-x},
\CV_X\underset{\CB}\otimes \CW_X)$ also has the same property.

We have three maps
$$j_*j^*(\CV_{X-x}\boxtimes \CW_{X-x}\boxtimes \CM)\to
\Delta_!\Bigl(\CV_x\underset{\CB_x}\otimes \CW_x
\underset{\CB_x}\otimes\CM
\Bigr)_{\Res(\nabla,\CV_X),\Res(\nabla,\CW_X)},$$ where the subscript
refers to the fact that we are taking coinvariants of the above two
endomorphisms acting on $\CV_x$ and $\CW_x$.

The first map is the composition
\begin{align*}
&j_*j^*(\CV_{X-x}\boxtimes \CW_{X-x}\boxtimes \CM)
\overset{\on{id}_{\CV_{X-x}}\boxtimes can_{\CW}}\longrightarrow
\Delta_{x_2=x_3}{}_!\Bigl(j_*j^*(\CV_{X-x}\boxtimes
(\CW_x\underset{\CB_x}\otimes\CM)_{\Res(\nabla,\CW_X)})\Bigr)
\overset{can_{\CV}}\longrightarrow \\ &\to
\Delta_!\Bigl(\CV_x\underset{\CB_x}\otimes
\CW_x\underset{\CB_x}\otimes \CM\Bigr)
_{\Res(\nabla,\CV_X),\Res(\nabla,\CW_X)}.
\end{align*}
The second map is obtained by interchanging the roles of $\CV$ and
$\CW$.  The third map if the composition
\begin{align*}
&j_*j^*(\CV_{X-x}\boxtimes \CW_{X-x}\boxtimes \CM_x) \to
\Delta_{x_1=x_2}{}_!\Bigl(j_*j^*((\CV_{X-x}\underset{\CB_{X-x}}\otimes
\CW_{X-x}) \boxtimes \CM_x)\Bigr)\overset{can_{\CV\otimes
\CW}}\longrightarrow \\ &\to
\Delta_!\Bigl(\CV_x\underset{\CB_x}\otimes
\CW_x\underset{\CB_x}\otimes \CM
\Bigr)_{\Res(\nabla,\CV_X)+\Res(\nabla,\CW_X)}\to
\Delta_!\Bigl(\CV_x\underset{\CB_x}\otimes
\CW_x\underset{\CB_x}\otimes \CM
\Bigr)_{\Res(\nabla,\CV_X),\Res(\nabla,\CW_X)}.
\end{align*}

\begin{lem}  \label{add assoc}
The sum of the three maps above is equal to zero.
\end{lem}

The proof follows by unfolding the construction of the bijection in
the proof of \propref{commutative fusion}.

\ssec{End of proof of \thmref{repr of pair from ten product}}

Consider the following general set-up. Let $\CA$ be a chiral algebra,
and let $\fz$ be its center. Let $\CV$ be a *-commutative chiral
module over $\fz$, free of finite rank as a $\fz$-module. Let
$\CM,\CN$ be two chiral $\CA$-modules.

Given a chiral pairing $\{\CV\underset{\fz}\otimes \CA,\CM\}\to
\CN$ we can restrict it and obtain a chiral pairing
$$j_*j^*(\CV\boxtimes \CM)\to \Delta_!(\CN)$$ of chiral
$\fz$-modules. This map has the following commutation property with
respect to $\CA$: the map
$$\wt{j}_*\wt{j}{}^*(\CV\boxtimes \CA\boxtimes \CM)\to
\Delta_{x_2=x_3}{}_!  \Bigl(j_*j^*(\CV\boxtimes \CM)\Bigr)\to
\Delta_!(\CN)$$ equals
$$\wt{j}_*\wt{j}{}^*(\CV\boxtimes \CA\boxtimes \CM)\to
\Delta_{x_1=x_3}{}_!  \Bigl(j_*j^*(\CA\boxtimes \CN)\Bigr)\to
\Delta_!(\CN),$$ where $\wt{j}$ denotes the embedding of the
complement of the union of the divisors $x_1=x_3$ and $x_2=x_3$.

We have the following general assertion:
\begin{lem}  \label{pairings over center}
The set of $\CA$-chiral pairings $\{\CV\underset{\fz}\otimes
\CA,\CM\}\to \CN$ is in a bijection with the set of $\fz$-chiral
pairings $\{\CV,\CM\}\to \CN$, satisfying the additional condition
above.
\end{lem}

Together with \propref{commutative fusion}, this lemma implies
\thmref{repr of pair from ten product}.

\ssec{Proof of \thmref{main, spherical}}

According to \secref{thm for diff op}, in order to prove \thmref{main,
spherical}, it suffices to consider the universal case, namely, the one
when $\CM=\fD_{G,\crit,x}$.  We will work in a more general framework,
assuming only that $\CM$ is such that $\CF_{V'}\star \CM$ is acyclic
away from cohomological degree $0$ for any $V'\in \Rep(\cG)$, which is
satisfied in the above case by \thmref{main, spherical, diff}. (Of
course, \thmref{main, spherical} will imply that this assumption is
satisfied for any $\CM$.)
 
Composing the map \eqref{desired map almost} with the tautological
projection
$$\CV_{\wt\fZ^{\int,\reg}_{\fg}}
\underset{\wt\fZ^{\int,\reg}_{\fg}}\otimes \CM \twoheadrightarrow
\Bigl(\CV_{\wt\fZ^{\int,\reg}_{\fg}}
\underset{\wt\fZ^{\int,\reg}_{\fg}}\otimes
\CM\Bigr)_{N_{\CV_{\wt\fZ^{\int,\reg}_{\fg}}}},$$ we obtain a map
\begin{equation} \label{desired map}
\CV_{\wt\fZ^{\int,\reg}_{\fg}}
\underset{\wt\fZ^{\int,\reg}_{\fg}}\otimes \CM \to \CF_V\star \CM.
 \end{equation}
 
We will denote this map by $\fs_V^{-1}$, as it will be the inverse of
the desired map of \thmref{main, spherical}. We will deduce
\thmref{main, spherical} from the following:
 
\begin{prop}  \label{map is assoc}  \hfill
The map of \eqref{desired map} is compatible with tensor products of
representations, in the sense that for $V,W\in \Rep(\cG)$ the 3
compositions
\begin{equation} \label{compo 1}
\CV_{\wt\fZ^{\int,\reg}_{\fg}}\underset{\wt\fZ^{\int,\reg}_{\fg}}
\otimes
\CW_{\wt\fZ^{\int,\reg}_{\fg}}\underset{\wt\fZ^{\int,\reg}_{\fg}}
\otimes \CM\overset{\fs^{-1}_V}\longrightarrow \CF_V\star
\Bigl(\CW_{\wt\fZ^{\int,\reg}_{\fg}}\underset{\wt\fZ^{\int,\reg}_{\fg}}
\otimes \CM\Bigr) \overset{\on{id}_{\CF_V}\star
\fs^{-1}_W}\longrightarrow \CF_V\star (\CF_W\star \CM),
\end{equation}
\begin{equation} \label{compo 2}
\CV_{\wt\fZ^{\int,\reg}_{\fg}}\underset{\wt\fZ^{\int,\reg}_{\fg}}
\otimes
\CW_{\wt\fZ^{\int,\reg}_{\fg}}\underset{\wt\fZ^{\int,\reg}_{\fg}}
\otimes \CM \overset{\on{id}_{\CV_{\wt\fZ^{\int,\reg}_{\fg}}}\otimes
\fs^{-1}_W}\longrightarrow
\CV_{\wt\fZ^{\int,\reg}_{\fg}}\underset{\wt\fZ^{\int,\reg}_{\fg}}
\otimes \Bigl(\CF_W\star \CM\Bigr) \overset{\fs^{-1}_V}\longrightarrow
\CF_V\star (\CF_W\star \CM)
\end{equation}
and
\begin{equation} \label{compo 3}
\CV_{\wt\fZ^{\int,\reg}_{\fg}}\underset{\wt\fZ^{\int,\reg}_{\fg}}
\otimes
\CW_{\wt\fZ^{\int,\reg}_{\fg}}\underset{\wt\fZ^{\int,\reg}_{\fg}}
\otimes \CM \overset{\fs^{-1}_{V\otimes W}}\to \CF_{V\otimes W}\star
\CM\simeq \CF_V\star (\CF_W\star \CM)
\end{equation}
coincide.
\end{prop}

Let us finish the proof of \thmref{main, spherical} modulo this
proposition. Without restriction of generality, we can assume that
the representation $V$ is finite-dimensional.

Suppose the map of \eqref{desired map} is not injective,
and let $\CN$ be its kernel. Let $V^*$ be representation dual to
$V$. Consider the diagram
$$ \CD \CM & @>{\on{id}}>> & \CM \\ @AAA & & @AAA \\
\CV^*_{\wt\fZ^{\int,\reg}_{\fg}}\underset{\wt\fZ^{\int,\reg}_{\fg}}
\otimes
\CV_{\wt\fZ^{\int,\reg}_{\fg}}\underset{\wt\fZ^{\int,\reg}_{\fg}}
\otimes \CM @>{\on{id}_{\CV^*}\otimes \fs^{-1}_V}>>
\CV^*_{\wt\fZ^{\int,\reg}_{\fg}}\underset{\wt\fZ^{\int,\reg}_{\fg}}
\otimes
(\CF_V\star \CM) @>{\fs^{-1}_{V^*}}>> \CF_{V^*}\star (\CF_V\star \CM)
\\ @AAA \\
\CV^*_{\wt\fZ^{\int,\reg}_{\fg}}\underset{\wt\fZ^{\int,\reg}_{\fg}}
\otimes
\CN.  \endCD
$$
It is commutative since the maps \eqref{compo 2}  and \eqref{compo 3} of
\propref{map is assoc} coincide. Hence, we obtain that the composed map
$$\CV^*_{\wt\fZ^{\int,\reg}_{\fg}}\underset{\wt\fZ^{\int,\reg}_{\fg}}
\otimes \CN\to \CM$$ is zero. But this is a contradiction since this
map is adjoint to the tautological embedding $\CN\to
\CV_{\wt\fZ^{\int,\reg}_{\fg}}\underset{\wt\fZ^{\int,\reg}_{\fg}}
\otimes \CM$.

In particular, we obtain that
$\CV_{\wt\fZ^{\int,\reg}_{\fg}}\underset{\wt\fZ^{\int,\reg}_{\fg}}
\otimes \CM \to
\Bigl(\CV_{\wt\fZ^{\int,\reg}_{\fg}}\underset{\wt\fZ^{\int,\reg}_{\fg}}
\otimes \CM\Bigr)_ {N_{\CV_{\wt\fZ^{\int,\reg}_{\fg}}}}$ is
injective, implying that the action of
$N_{\CV_{\wt\fZ^{\int,\reg}_{\fg}}}$ on
$\CV_{\wt\fZ^{\int,\reg}_{\fg}}\underset{\wt\fZ^{\int,\reg}_{\fg}}
\otimes \CM$ is trivial.

\medskip

Suppose now that the map of \eqref{desired map} is not surjective, and
let $\CN'$ be its kernel.  Let us recall that the functors $\CF_V\star
?$ and $\CF_{V^*}\star ?$ are mutually (both left and right) adjoint
on the category $D(\hg_\crit\mod)^{G[[t]]}$. Consider the diagram
$$ \CD & & & & \CF_{V^*}\star \CN' \\ & & & & @AAA \\
\CV^*_{\wt\fZ^{\int,\reg}_{\fg}}\underset{\wt\fZ^{\int,\reg}_{\fg}}
\otimes
\CV_{\wt\fZ^{\int,\reg}_{\fg}}\underset{\wt\fZ^{\int,\reg}_{\fg}}
\otimes \CM @>{\fs^{-1}_{V^*}}>> \CF_{V^*}\star
(\CV_{\wt\fZ^{\int,\reg}_{\fg}}\underset{\wt\fZ^{\int,\reg}_{\fg}}
\otimes \CM) @>{\on{id}_{\CF_{V^*}}\star \fs^{-1}_V}>> \CF_{V^*}\star
(\CF_V\star \CM) \\ @AAA & & @AAA \\ \CM & @>{\on{id}}>> & \CM.
\endCD
$$
This diagram is commutative because the maps \eqref{compo 1} and
\eqref{compo 3} of \propref{map is assoc} coincide.  Hence, we obtain
that the composed map $\CM\to \CF_{V^*}\star \CN'$ is zero, which is a
contradiction.

\ssec{Proof of \propref{map is assoc}}

Consider the diagram
$$ \CD j_*j^*\Bigl(\Gamma(\Gr_{G,X},\CF_{V,X})\boxtimes
(\CW_{\wt\fZ^{\int,\reg}_{\fg}}\underset{\wt\fZ^{\int,\reg}_{\fg}}
\otimes \CM)\Bigr) @>{\fs^{-1}_W}>>
j_*j^*\Bigl(\Gamma(\Gr_{G,X},\CF_{V,X})\boxtimes (\CF_W\star
\CM)\Bigr) \\ @V{can_V}VV @V{can_V}VV \\
\Delta_!\Bigl(\CV_{\wt\fZ^{\int,\reg}_{\fg}}
\underset{\wt\fZ^{\int,\reg}_{\fg}}\otimes
(\CW_{\wt\fZ^{\int,\reg}_{\fg}}\underset{\wt\fZ^{\int,\reg}_{\fg}}
\otimes \CM)\Bigr)_{N_{\CV_{\wt\fZ^{\int,\reg}_{\fg}}}}
@>{\fs^{-1}_W}>> \Delta_!
\Bigl(\CV_{\wt\fZ^{\int,\reg}_{\fg}}\underset{\wt\fZ^{\int,\reg}_{\fg}}
\otimes (\CF_W\star \CM)\Bigr)_{N_{\CV_{\wt\fZ^{\int,\reg}_{\fg}}}}
\\ @V{\fs^{-1}_V}VV @V{\fs^{-1}_V}VV \\
\Delta_!\Bigl(\CF_V\star(\CW_{\wt\fZ^{\int,\reg}_{\fg}}
\underset{\wt\fZ^{\int,\reg}_{\fg}}\otimes \CM)\Bigr)
@>{\fs^{-1}_W}>> \Delta_!(\CF_V\star (\CF_W\star \CM)).  \endCD
$$

To prove that \eqref{compo 1} and \eqref{compo 2} coincide we have to
show the commutativity of the lower square in this diagram.

Note that the upper square is commutative by functoriality of the map
$can_V$; moreover the vertical maps in this square are surjective, by
\propref{commutative fusion}.  Hence, it suffices to see that the
outer square is commutative, but this follows from the functoriality
of the map \eqref{pairing with sph}.

\medskip

Consider now the composed map
$$ \CD j_*j^*\Bigl(\Gamma(\Gr_{G,X},\CF_{V,X})\boxtimes
\Gamma(\Gr_{G,X},\CF_{W,X})\boxtimes \CM\Bigr) \\ @VVV \\
\Delta_{x_2=x_3}{}_!\Bigl(j_*j^*\Bigl(\Gamma(\Gr_{G,X},\CF_{V,X})
\boxtimes
(\CW_{\wt\fZ^{\int,\reg}_{\fg}}\underset{\wt\fZ^{\int,\reg}_{\fg}}
\otimes \CM)_{N_{\CW_{\wt\fZ^{\int,\reg}_{\fg}}}} \Bigr)\Bigr) \\
@VVV \\ \Delta_!\Bigl(\CV_{\wt\fZ^{\int,\reg}_{\fg}}
\underset{\wt\fZ^{\int,\reg}_{\fg}}\otimes
\CW_{\wt\fZ^{\int,\reg}_{\fg}}\underset{\wt\fZ^{\int,\reg}_{\fg}}
\otimes \CM\Bigr)
_{N_{\CV_{\wt\fZ^{\int,\reg}_{\fg}}},N_{\CW_{\wt\fZ^{\int,\reg}_{\fg}}}}
\\ @V{\text{\eqref{compo 1}}}VV \\ \Delta_!(\CF_V\star \CF_W\star
\CM), \endCD
$$
and a similar composition when the roles of $V$ and $W$ are
interchanged.  The resulting maps are equal to the first and the
second maps, respectively, of \propref{comp with ten products}.

Hence, their sum equals the composed map from the commutative diagram
$$ \CD j_*j^*\Bigl(\Gamma(\Gr_{G,X},\CF_{V,X})\boxtimes
\Gamma(\Gr_{G,X},\CF_{W,X})\boxtimes \CM\Bigr) \\ @VVV \\
\Delta_{x_1=x_2}{}_!\Bigl(j_*j^*\Bigl(\Gamma(\Gr_{G,X},\CF_{V\otimes
W,X})\boxtimes \CM \Bigr)\Bigr) @>>> \Delta_!(\CF_{V\otimes W}\star
\CM) \\ @VVV @V{\sim}VV \\ \Delta_!\Bigl(
\CV_{\wt\fZ^{\int,\reg}_{\fg}}\underset{\wt\fZ^{\int,\reg}_{\fg}}\otimes
\CW_{\wt\fZ^{\int,\reg}_{\fg}}\underset{\wt\fZ^{\int,\reg}_{\fg}}\otimes
\CM\Bigr)_
{N_{\CV_{\wt\fZ^{\int,\reg}_{\fg}}}+N_{\CW_{\wt\fZ^{\int,\reg}_{\fg}}}}
@>{\text{\eqref{compo 3}}}>> \Delta_!(\CF_V\star \CF_W\star \CM),
\endCD
$$
by \propref{comp with ten products}. Moreover, by \lemref{add
assoc}, the sum of the three composed maps
\begin{align*}
&j_*j^*\Bigl(\Gamma(\Gr_{G,X},\CF_{V,X})\boxtimes
\Gamma(\Gr_{G,X},\CF_{W,X})\boxtimes \CM\Bigr) \to \\ &\to
\Delta_!\Bigl(\CV_{\wt\fZ^{\int,\reg}_{\fg}}
\underset{\wt\fZ^{\int,\reg}_{\fg}}\otimes
\CW_{\wt\fZ^{\int,\reg}_{\fg}}\underset{\wt\fZ^{\int,\reg}_{\fg}}
\otimes \CM\Bigr)
_{N_{\CV_{\wt\fZ^{\int,\reg}_{\fg}}},N_{\CW_{\wt\fZ^{\int,\reg}_{\fg}}}}
\end{align*}
is also zero.

Let $\CM_{V,W}\subset j_*j^*\Bigl(\Gamma(\Gr_{G,X},\CF_{V,X})\boxtimes
\Gamma(\Gr_{G,X},\CF_{W,X})\boxtimes \CM\Bigr)$ be the kernel of the
map to
$$\Delta_{x_1=x_3}{}_!\Bigl(j_*j^*\Bigl(\Gamma(\Gr_{G,X},\CF_{W,X})
\boxtimes
(\CV_{\wt\fZ^{\int,\reg}_{\fg}}\underset{\wt\fZ^{\int,\reg}_{\fg}}
\otimes \CM)_ {N_{\CV_{\wt\fZ^{\int,\reg}_{\fg}}}}\Bigr)\Bigr).$$

Then the map 
$$\CM_{V,W}\to \Delta_{x_2=x_3}{}_!
\Bigl(j_*j^*\Bigl(\Gamma(\Gr_{G,X},\CF_{V,X})\boxtimes
(\CW_{\wt\fZ^{\int,\reg}_{\fg}}\underset{\wt\fZ^{\int,\reg}_{\fg}}
\otimes \CM)_ {N_{\CW_{\wt\fZ^{\int,\reg}_{\fg}}}}\Bigr)\Bigr)$$ is
still surjective.  Hence, the two surviving maps
$$\CM_{V,W}\to \Delta_!\Bigl(\CV_{\wt\fZ^{\int,\reg}_{\fg}}
\underset{\wt\fZ^{\int,\reg}_{\fg}}\otimes
\CW_{\wt\fZ^{\int,\reg}_{\fg}}\underset{\wt\fZ^{\int,\reg}_{\fg}}
\otimes \CM\Bigr) _{N_{\CV_{\wt\fZ^{\int,\reg}_{\fg}}},
N_{\CW_{\wt\fZ^{\int,\reg}_{\fg}}}}$$ coincide and are
surjective. By the above, the two surviving maps
$$\CM_{V,W}\to \Delta_!(\CF_V\star \CF_W\star \CM)$$ coincide as well,
and the equality of the maps \eqref{compo 1} and \eqref{compo 3}
follows.

\section{Convolution and fusion for general chiral algebras}
\label{proof of rep}

The goal of this section is to prove \thmref{representability,
spherical}. In fact, we will prove a more general result, valid for
any chiral algebra, endowed with a Harish-Chandra action of
$\on{Jets}(G)_X$. We will make a more extensive use of the formalism
developed in the appendix to \cite{FG2}, but we should note that the
main results of this paper, \thmref{main, spherical} and
\thmref{main}, are independent of this section.

\ssec{Twisting of chiral modules}    \label{general alg}

Recall the Lie-* algebras $L_\fg$ and $L_{\fg,\kappa}$ introduced in
\secref{recol}. Let $\on{Jets}(G)_X$ denote the group-like object in
the category of D-schemes on $X$ corresponding to jets into $G$. Its
relative cotangent sheaf is the D-module on $X$ canonically isomorphic
to the dual $L_\fg^\vee$ of $L_\fg$.

\medskip

Let $\CA$ be a chiral algebra on $X$, endowed with an action of
$\on{Jets}(G)_X$, such that the adjoint action of its Lie algebra is
inner at the level $\kappa$. In other words, we assume being given a
homomorphism of Lie-* algebras $L_{\fg,\kappa}\to \CA$, such that the
map
$$\CA\to \CA\otimes L_\fg^\vee,$$ corresponding to the Lie-* bracket
equals the derivative of the action of $\on{Jets}(G)$ on $\CA$.

\begin{propconstr}  \label{group act on cat}
The category of chiral $\CA$-modules, supported at $x\in X$ carries an
action of the group-scheme $G\ppart$ of Harish-Chandra type with
respect to the central extension $\hg_\kappa$ (cf. \cite{FG2},
Sect. 22).
\end{propconstr}

We will give two proofs. One, discussed below, is local in terms of
the De Rham cohomology of $\CA$ on the formal punctured disc
$\D_x^\times$ around $x$.  Another proof will be of chiral nature.

\begin{proof}

Let us first recall the following general construction,
cf. \cite{CHA}, 3.6.9.  For a chiral algebra $\CA$ let us consider the
topological vector space $H_{DR}^0(\D_x^\times,\CA)$. We remind that
for a D-module $\CV$ on $X$, we define
$$H_{DR}^0(\D_x^\times,\CV):=
\underset{\CL}{\underset{\longrightarrow}{\lim}}\,
H_{DR}^0(\D_x^\times,\CL),$$ where $\CL$ runs over the filtered set of
finitely-generated D-submodules of $\CV$, and the direct limit is
taken in the category of topological vector spaces.
\footnote{In other words, we first take the limit in the category of
  vector spaces, and then complete it in the natural topology.}
  Consider in addition the topological vector space
  $H_{DR}^0(\D_x^\times\times \D_x^\times-\Delta_{\D_x^\times},
  \CA\boxtimes \CA)$, which is acted on by the transposition
  $\sigma$. We have a canonical map
$$H_{DR}^0(\D_x^\times\times
  \D_x^\times-\Delta_{\D_x^\times},\CA\boxtimes \CA)\to
  H_{DR}^0(\D_x^\times,\CA)\arrowtimes H_{DR}^0(\D_x^\times,\CA),$$
  defined in fact for any pair of D-modules on $X$, cf. \cite{CHA}
  3.6.9.

In addition, the
structure of chiral algebra defines a map
$$\{\cdot,\cdot\}:H_{DR}^0(\D_x^\times\times
\D_x^\times-\Delta_{\D_x^\times},\CA\boxtimes \CA)\to
H_{DR}^0(\D_x^\times,\CA).$$

The category $\CA\mod_x$ of chiral $\CA$-modules supported at $x$
identifies with the category of vector spaces $\CM$, endowed with an
action map
$$\on{act}_\CM:H_{DR}^0(\D_x^\times,\CA)\arrowtimes \CM\to \CM,$$ such
that the difference of
\begin{align*}
&H_{DR}^0(\D_x^\times\times
  \D_x^\times-\Delta_{\D_x^\times},\CA\boxtimes \CA)\arrowtimes \CM\to
  H_{DR}^0(\D_x^\times,\CA)\arrowtimes
  H_{DR}^0(\D_x^\times,\CA)\arrowtimes \CM
  \overset{\on{act}_\CM}\longrightarrow \\ &\to
  H_{DR}^0(\D_x^\times,\CA)\arrowtimes \CM
  \overset{\on{act}_\CM}\longrightarrow \CM
\end{align*}
and the map obtained by first acting by $\sigma$ on the first factor, equals 
$$H_{DR}^0(\D_x^\times\times
\D_x^\times-\Delta_{\D_x^\times},\CA\boxtimes \CA)\arrowtimes \CM
\overset{\{\cdot,\cdot\}\otimes \on{id}}\longrightarrow
H_{DR}^0(\D_x^\times,\CA)\arrowtimes \CM
\overset{\on{act}_\CM}\longrightarrow \CM.$$

\medskip

Let $\bg$ be an $S$-point of $G\ppart$ for some base-scheme $S$. It
gives rise to a map $$\CO_{\on{Jets}(G)_X}\to \CO_S\ppart,$$
compatible with the connection. Composing it with the map $\CA\to
\CO_{\on{Jets}(G)_X}\underset{\CO_X}\otimes \CA$, given by the action
of $\on{Jets}(G)_X$ on $\CA$, we obtain a map
$$H_{DR}^0(\D_x^\times,\CA) \to \CO_S\shriektimes
H_{DR}^0(\D_x^\times,\CA).$$
Similarly, we have a map
$$H_{DR}^0(\D_x^\times\times
\D_x^\times-\Delta_{\D_x^\times},\CA\boxtimes \CA)\to
\CO_S\shriektimes H_{DR}^0(\D_x^\times\times
\D_x^\times-\Delta_{\D_x^\times},\CA\boxtimes \CA).$$

\medskip

Given an object $\CM\in \CA\mod_x$, we define an action map
$$H_{DR}^0(\D_x^\times,\CA) \arrowtimes (\CM\otimes \CO_S)\to
\CM\otimes \CO_S$$ as a composition
\begin{align*}
&H_{DR}^0(\D_x^\times,\CA) \arrowtimes (\CM\otimes \CO_S)\to
\Bigl(\CO_S\shriektimes H_{DR}^0(\D_x^\times,\CA)\Bigr)\arrowtimes
(\CM\otimes \CO_S)\to \\ &\to \CO_S\shriektimes
\Bigl(H_{DR}^0(\D_x^\times,\CA)\arrowtimes \CM)\shriektimes
\CO_S\overset{\on{act}_\CM}\to \CO_S\otimes \CM\otimes \CO_S\to
\CM\otimes \CO_S,
\end{align*}
where the last arrow is given by the multiplication on $\CO_S$.

By construction, it follows that the relation, that singles out
representations among all vector spaces endowed with an action of
$H_{DR}^0(\D_x^\times,\CA)$ (cf. above), holds. Thus, we obtain a
$G\ppart$-action on $\CA\mod_x$, which is of Harish-Chandra type by
construction.

\end{proof}

\ssec{Twisting and fusion}

Let now $\CM'$ be a torsion-free chiral $\CA$-module on $X$. Assume
that $\CM'$ is weakly $\on{Jets}(G)_X$-equivariant. I.e., we have an
action of $\on{Jets}(G)_X$ on $\CM'$, compatible with its action on
$\CA$ in the natural sense.

Let $\CM$, $\CN$ be two chiral $\CA$-modules, both supported at the
point $x\in X$. Let $\bg$ be an $S$-point of $G\ppart$ and let
$\CM^{\bg}$ and $\CN^{\bg}$ the corresponding $S$-families of objects
of $\CA\mod_x$, defined by \propconstrref{group act on cat}.

\begin{propconstr} \label{twisting pairing}
To every chiral pairing 
\begin{equation} \label{initial pairing}
\phi:j_*j^*(\CM'\boxtimes \CM)\to \Delta_!(\CN)
\end{equation}
there functorially corresponds a chiral pairing 
$$\phi:j_*j^*(\CM'\boxtimes \CM^\bg)\to \Delta_!(\CN^\bg).$$
\end{propconstr}

\begin{proof}

The proof is largely parallel to that of \propconstrref{group act on
cat} above.

First, let $\CM'$ be any torsion-free module over a chiral algebra
$\CA$.  Consider the topological vector spaces
$H^0_{DR}(\D_x^\times,\CM')$ and $H_{DR}^0(\D_x^\times\times
\D_x^\times-\Delta_{\D_x^\times},\CA\boxtimes \CM')$; the action of
$\CA$ on $\CM'$ gives rise to a map
$$\on{act}_{\CM'}:
H_{DR}^0(\D_x^\times\times \D_x^\times-\Delta_{\D_x^\times},
\CA\boxtimes \CM')\to H^0_{DR}(\D_x^\times,\CM').$$

For two objects $\CM,\CN\in \CA\mod_x$, chiral pairings
$\{\CM',\CM\}\to \CN$ are in bijection with maps
$$\phi:H^0_{DR}(\D_x^\times,\CM')\arrowtimes \CM\to \CN,$$
such that two difference of the two compositions
\begin{align*}
&H_{DR}^0(\D_x^\times\times
\D_x^\times-\Delta_{\D_x^\times},\CA\boxtimes \CM')\arrowtimes \CM \to
H^0_{DR}(\D_x^\times,\CA)\arrowtimes
H^0_{DR}(\D_x^\times,\CM')\arrowtimes \CM \overset{\on{id}\CA\otimes
\phi} \longrightarrow \\ &\to H^0_{DR}(\D_x^\times,\CA)\arrowtimes \CN
\overset{\on{act}_\CN}\longrightarrow \CN
\end{align*}
and
\begin{align*}
&H_{DR}^0(\D_x^\times\times
\D_x^\times-\Delta_{\D_x^\times},\CA\boxtimes \CM')\arrowtimes \CM \to
H^0_{DR}(\D_x^\times,\CM') \arrowtimes
H^0_{DR}(\D_x^\times,\CA)\arrowtimes \CM
\overset{\on{id}_{\CM'}\otimes \on{act}_\CM}\longrightarrow \\ &\to
H^0_{DR}(\D_x^\times,\CM') \arrowtimes \CM \overset{\phi}\to \CN
\end{align*}
equals
$$H_{DR}^0(\D_x^\times\times
\D_x^\times-\Delta_{\D_x^\times},\CA\boxtimes \CM')\arrowtimes \CM
\overset{\on{act}_{\CM'}}\longrightarrow
H^0_{DR}(\D_x^\times,\CM')\arrowtimes \CM\overset{\phi}\to \CN.$$

\medskip

Suppose now that $\CM'$ is weakly $\on{Jets}(G)_X$-equivariant. Given
an $S$-point $\bg$ of $G\ppart$ as above, as in the case of $\CA$, we
obtain a map
$$H^0_{DR}(\D_x^\times,\CM')\to \CO_S\shriektimes
H^0_{DR}(\D_x^\times,\CM').$$

\medskip

For a chiral pairing $\phi$ as above, we 
define a chiral pairing
$$H^0_{DR}(\D_x^\times,\CM')\arrowtimes (\CM\otimes \CO_S)\to
(\CN\otimes \CO_S)$$ as a composition
\begin{align*}
&H^0_{DR}(\D_x^\times,\CM')\arrowtimes (\CM\otimes \CO_S)\to
\Bigl(\CO_S\shriektimes H^0_{DR}(\D_x^\times,\CM')\Bigr) \arrowtimes
(\CM\otimes \CO_S) \to \\ &\to\CO_S\shriektimes
\Bigl(H^0_{DR}(\D_x^\times,\CM')\arrowtimes \CM\Bigr) \shriektimes
\CO_S\overset{\on{id}\otimes \phi\otimes \on{id}}\longrightarrow
\CO_S\otimes \CN\otimes \CO_S\to \CN\otimes \CO_S.
\end{align*}

It is easy to see that the relation involving the $\CA$-actions on
$\CM$ and $\CN$, described above, holds.

\end{proof}

Let us assume now that in the circumstances of the above proposition,
module $\CM'$ is strongly $\on{Jets}(G)_X$-equivariant. By definition,
this means that the derivative of the group-action, which is a map
$$\CM'\to \CM'\otimes L_\fg^\vee,$$
coincides with the one coming from the Lie-* action of
$L_{\fg,\kappa}$ via $L_{\fg,\kappa}\to \CA$.

Let $\bg_1$ and $\bg_2$ be two $S$-points of $G\ppart$, whose ratio is
a map from $S$ to the first infinitesimal neighborhood of the identity
in $G\ppart$. In particular, for every choice of the splitting
$\fg\ppart\to \hg_\kappa$ we have the canonical isomorphisms
$\CM^{\bg_1}\simeq \CM^{\bg_2}$ and $\CN^{\bg_1}\simeq \CN^{\bg_2}$.

{}From the proof of \propconstrref{twisting pairing} we obtain:

\begin{cor}  \label{strong twist}
Under the above circumstances, the diagram
$$
\CD
j_*j^*(\CM'\boxtimes \CM^{\bg_1})@>>> \Delta_!(\CN^{\bg_1}) \\
@V{\sim}VV   @V{\sim}VV  \\
j_*j^*(\CM'\boxtimes \CM^{\bg_2})@>>> \Delta_!(\CN^{\bg_2})
\endCD
$$
is commutative.
\end{cor}

Let now $\CF'$ be a $\kappa$-twisted D-module on $G\ppart$, which
strongly $K$-equivariant on the right, where $K$ is an "open-compact"
group-subscheme of $G\ppart$. By \cite{FG2}, Sect. 22.4, given an
object $\CM\in \CA\mod_x$, which is strongly $K$-equivariant, there
exists a well-defined complex of objects in $\CA\mod_x$:
$$\fC^\semiinf\Bigl(\fg\ppart;K_{red},\CF'\otimes \CM\Bigr),$$
where $K_{red}$ denotes the reductive quotient of $K$.

The image of this complex in the derived category is by definition the
convolution $\CF\star \CM$, where $\CF$ is the twisted D-module on
$G/K$, corresponding to $\CF'$.

{}From \propconstrref{twisting pairing} above, we obtain that given two 
strongly $K$-equivariant objects
$\CM,\CN\in \CA\mod_x$ and a
chiral pairing $\{\CM',\CM\}\to \CN$ we obtain a chiral pairing of graded
objects
\begin{equation} \label{conv pairing}
\{\CM',\fC^\semiinf\Bigl(\fg\ppart;K_{red},\CF'\otimes \CM\Bigr)\}\to 
\fC^\semiinf\Bigl(\fg\ppart;K_{red},\CF'\otimes \CN\Bigr).
\end{equation}
Moreover, by \corref{strong twist}, the above pairing is a map of
complexes, i.e., it respects the differentials on both sides.

\medskip

As an application we shall now establish the following result.  Let
$\on{Av}_{G[[t]]}$ denote the functor $D^+(\CA\mod_x)\to
D^+(\CA\mod_x)^{G[[t]]}$, right adjoint to the forgetful functor. In
particular, it is left-exact and the functor
$$\CN\to h^0\Bigl(\on{Av}_{G[[t]]}(\CN)\Bigr)$$ is the right adjoint
to the forgetful functor $\CA\mod_x^{G[[t]]}\to \CA\mod_x$.

\begin{prop}  \label{equivariance preserved}
Let $\{\CM',\CM\}\to \CN$ be a chiral pairing with $\CM'$ being
strongly $\on{Jets}(G)_X$-equivariant, and $\CM\in \CA\mod_x$ being
$G[[t]]$-equivariant.  Then this pairing canonically factors through
$\{\CM',\CM\}\to h^0\Bigl(\on{Av}_{G[[t]]}(\CN)\Bigr)$.
\end{prop}

\begin{proof}

Let us recall that the functor $\on{Av}_{G[[t]]}$ is represented by
the complex
$$\CN\mapsto \fC^\bullet\Bigl(\fg[[t]],\CO_{G[[t]]}\otimes
\CN\Bigr).$$ Hence, as in \eqref{conv pairing}, given a chiral pairing
$\{\CM',\CM\}\to \CN$, we obtain a chiral pairing of complexes
$$\{\CM',\on{Av}_{G[[t]]}(\CM)\}\to \on{Av}_{G[[t]]}(\CN),$$
compatible with the differential. Composing with the canonical map
$\CM\to \on{Av}_{G[[t]]}(\CM)$, we obtain the desired chiral pairing
$$\{\CM',\CM\}\to \on{Av}_{G[[t]]}(\CN).$$

\end{proof}

\ssec{The global case}   \label{global case}

We shall now generalize the discussion of the previous subsections to
the case of chiral $\CA$-modules, which are not necessarily supported
at a single point $x\in X$.

Recall that in addition to $\on{Jets}(G)_X$, we have the group D-ind
scheme $\on{Jets}^{mer}(G)_X$.

\begin{propconstr} \label{global}  \hfill

\smallskip

\noindent{\em (1)} Let $\CM$ be a chiral $\CA$-module on $X$, and let
$\bg$ be an $S$-point of $\on{Jets}^{mer}(G)$, where $S$ is an affine
D-scheme on $X$. Then the D-module
$\CM^\bg:=\CM\underset{\CO_X}\otimes \CO_S$ acquires a natural
structure of chiral $\CA$-module.

\smallskip

\noindent{\em (2)} If the ratio of two points $\bg_1$ and $\bg_2$ lies
in the first infinitesimal neighborhood of the unit section of
$\on{Jets}^{mer}(G)$, then for every choice of the splitting as a
D-module $L_\fg\to L_{\fg,\kappa}$ there exists a functorial
isomorphism $\CM^{\bg_1}\simeq \CM^{\bg_2}$.

\smallskip

\noindent{\em (3)} If $\CM'$ is a chiral $\CA$-module, which is weakly
$\on{Jets}(G)$-equivariant, then to every chiral pairing
$\{\CM',\CM\}\to \CN$ there functorially corresponds a chiral pairing
$\{\CM',\CM^\bg\}\to \CN^\bg$.

\smallskip

\noindent{\em (4)} In the circumstances of points (2) and (3) above
assume in addition that $\CM'$ is strongly
$\on{Jets}(G)$-equivariant. Then the diagram of \corref{strong twist}
commutes.

\end{propconstr}

\begin{proof}

Let us first recall the following general construction. Let 
$$j_*j^*(\CM_1^i\boxtimes \CM_2^i)\to \Delta_!(\CN^i),$$
be maps of D-modules, where $i$ runs over some finite set of indices $i$.
Then we have a map
\begin{equation} \label{mult chiral}
j_*j^*\Bigl((\underset{i}\otimes \CM_1^i)\boxtimes
(\underset{i}\otimes \CM_2^i)\Bigr) \to \Delta_!(\underset{i}\otimes
\CN^i).
\end{equation}

Let us recall also that the data of an $S$-point of
$\on{Jets}^{mer}(G)$ is equivalent to that of a map
$$j_*j^*(\CO_{\on{Jets}(G)_X}\boxtimes \CO_X)\to \Delta_!(\CO_S),$$
such that the diagram
$$ \CD j_*j^*\Bigl((\CO_{\on{Jets}(G)_X}\otimes
\CO_{\on{Jets}(G)_X})\boxtimes \CO_X\Bigr) @>>>
\Delta_!\Bigl((\CO_S\otimes \CO_S)\boxtimes \CO_X\Bigr) \\ @VVV @VVV
\\ j_*j^*(\CO_{\on{Jets}(G)_X}\boxtimes \CO_X) @>>> \Delta_!(\CO_S)
\endCD
$$
commutes, where the upper horizontal arrow comes from the map \eqref{mult chiral},
and vertical arrows are given by the algebra multiplication. 

\medskip

As in the proofs of \propconstrref{group act on cat} and
\propconstrref{twisting pairing}, the proof follows from the next
general construction. Let
$$j_*j^*(\CM_1\boxtimes \CM_2)\to \Delta_!(\CN) \text{ and }  
\CM_1\to \CM_1\otimes \CO_{\on{Jets}(G)_X}$$ be maps of D-modules
and $\bg$ be as above. Then from \eqref{mult chiral} we obtain a map
$$j_*j^*(\CM_1\boxtimes \CM_2)\to 
j_*j^*\Bigl((\CM_1\otimes \CO_{\on{Jets}(G)_X})\boxtimes
(\CM_2\boxtimes \CO_S)\Bigr) 
\to \Delta_!(\CN\otimes \CO_S),$$
which, in turn, gives rise to a map
\begin{equation} \label{acted map}
j_*j^*\Bigl(\CM_1\boxtimes (\CM_2\otimes \CO_S)\Bigr)\to
\Delta_!(\CN\otimes \CO_S).
\end{equation}

By putting first $\CM_1:=\CA$, $\CM_2:=\CM$ and $\CN:=\CM$, we arrive
to the chiral action map of point (1). By putting $\CM_1:=\CM'$,
$\CM_2:=\CM$ and $\CN:=\CN$, we arrive to the chiral pairing of point
(3).

\medskip

To prove point (2), we can assume being given a map 
$$j_*j^*(L^\vee_\fg\boxtimes \CO_X)\to \Delta_!(\CO_S),$$
and we have to construct a map $\varphi:\CM\to \CM\otimes \CO_S$,
which fits into the commutative diagram:
$$ \CD j_*j^*(\CA\boxtimes \CM) @>>> j_*j^*\Bigl((\CA\otimes
L^\vee_\fg)\boxtimes \CM\Bigr) \\ @V{\on{id}\boxtimes \varphi}VV @VVV
\\ j_*j^*\Bigl(\CA\boxtimes (\CM\otimes \CO_S)\Bigr) @>>>
\Delta_!(\CM\otimes \CO_S), \endCD
$$
where the right vertical arrow comes from \eqref{mult chiral}, and
the bottom horizontal arrow comes from the initial chiral action of
$\CA$ on $\CM$.

The desired map $\varphi$ is constructed as follows. The chiral
bracket of $L_\fg$ with $\CM$ (which is well-defined since we chose a
splitting of $L_{\fg,\kappa}$) and \eqref{mult chiral} give rise to a
map
$$j_*j^*\Bigl((L_\fg\otimes L^\vee_\fg)\boxtimes \CM\Bigr)\to
\Delta_!(\CM\otimes \CO_S).$$ The required map if the Lie-* bracket
induced by the above map applied to the canonical element ${\bf 1}\in
H^0_{DR}(L_\fg\otimes L^\vee_\fg)$.

\medskip

The fact that the axioms are satisfied, and point (4) of the
proposition follow from the construction.

\end{proof}

Let $\CM$ be a chiral $\CA$-module as above, and let $\CF'$ be a chiral
module over $\fD_{G,\kappa}$. Then \propconstrref{global} implies that
proceeding as in \cite{FG2}, Sect. 22.4, we can form a
twisted product of $\CF'$ and $\CM$, denoted $\CF'\tboxtimes \CM$.

Let us write down this construction explicitly. As a D-module, this
will be the usual tensor product $\CF'\otimes \CM$, but it will carry
a new action of $\CA$, and a commuting action of
$\CA_{\fg,2\kappa_\crit}$.  (Note that $2\kappa_\crit$ equals the
negative of the Killing form on $\fg$.)

Namely, the action of $\CA$ is the composition:
$$j_*j^*(\CA\boxtimes (\CF'\otimes \CM))\to
j_*j^*((\CO_{\on{Jets}(G)_X}\otimes \CA)\boxtimes (\CF'\otimes
\CM))\to \Delta_!(\CF'\otimes \CM),$$ where the last arrow comes from
the chiral action of $\CO_{\on{Jets}(G)_X}$ on $\CF$ and $\CA$ on
$\CM$ via \eqref{mult chiral}.

The chiral action of $L_{\fg,2\kappa_\crit}$ is the diagonal one with
respect to the $\fr$-action of $\CA_{\fg,\kappa'}$ on $\CF'$ and the
action of $L_{\fg,\kappa}$ on $\CM$ that comes from $L_{\fg,\kappa}\to
\CA$.

Hence, by tensoring with tensoring with the Clifford chiral algebra,
we obtain a well-defined complex $\fC^\semiinf(L_\fg,\CF'\otimes \CM)$
of chiral $\CA$-modules.  More generally, if $\{\CM',\CM\}\to \CN$ is
a chiral pairing of $\CA$-modules, with $\CM'$ being
$\on{Jets}(G)_X$-equivariant, as in point (3) of the above
proposition, we obtain a chiral pairing of complexes
$$j_*j^*\Bigl(\CM'\boxtimes \fC^\semiinf(L_\fg,\CF'\otimes \CM)\Bigr)\to
\Delta_!\Bigl(\fC^\semiinf(L_\fg,\CF'\otimes \CN)\Bigr).$$

In particular, if $\CF'$ is strongly equivariant with respect to the
right action of $\on{Jets}(G)_X$ and $\CM$ is also
$\on{Jets}(G)_X$-equivariant, by considering the corresponding
subcomplex of chains relative to $\fg\in \Gamma(X,\L_\fg)$, we obtain
a map
\begin{equation}   \label{global pairing}
j_*j^*\Bigl(\CM'\boxtimes \fC^\semiinf(L_\fg;\fg,\CF'\otimes \CM)\Bigr)\to
\Delta_!\Bigl(\fC^\semiinf(L_\fg;\fg,\CF'\otimes \CN)\Bigr),
\end{equation}
and on the level of individual cohomologies the chiral pairings
\begin{equation} \label{global, cohomologies}
j_*j^*(\CM'\boxtimes h^i(\CF\star \CM))\to \Delta_!(h^i(\CF\star \CN)),
\end{equation}
where $\CF$ denotes the corresponding twisted D-module on $\Gr_{G,X}$.

\ssec{Some compatibilities}

Let us take as an example the case when $\CM=\CA$, and the canonical
chiral pairing $\{\CM',\CA\}\to \CM'$. We obtain the chiral pairings
\begin{equation}  \label{general basic pairing}
\{\CM',h^i(\CF\star \CA)\}\to h^i(\CF\star \CM').
\end{equation}

Consider the case $\CA\simeq \CA_{\fg,\kappa}$. Then, by construction,
$$\CF\star \CA_{\fg,\kappa}\simeq \Gamma(\Gr_{G,X},\CF).$$

\begin{lem}  \label{two def}
For $\CA\simeq \CA_{\fg,\kappa}$, the chiral pairings of
\eqref{general basic pairing} coincide with those of \eqref{basic
pairing for modules}.
\end{lem}

\begin{proof}

By the construction of the pairings in \eqref{general basic pairing},
it is sufficient to consider the case when $\CM'\simeq
\fD_{G,\kappa}$. The latter reduces to the case of the chiral algebra
$\fD_{G,\kappa}$ rather than $\CA_{\fg,\kappa}$.

Now the assertion of the lemma follows from the next general
observation: if $\CM'=\CA$, then the map of \propconstrref{global}(3)
coincides with the chiral action of $\CA$ on $\CM^\bg$. In particular,
the maps of \eqref{global, cohomologies} are also given by the chiral
action of $\CA$ on $h^i(\CF\star \CM)$.

\end{proof}

Let us now establish some further compatibilities, satisfied by the maps 
of \propconstrref{global}.

\medskip

Let $\CM_1,\CM_2,\CM_3$ be chiral modules over $\CA$, and let
$\{\CM_1,\CM_2\}\to \CM_{1,2}$, $\{\CM_2,\CM_3\}\to \CM_{2,3}$,
$\{\CM_1,\CM_3\}\to \CM_{1,3}$ be chiral pairings. In addition,
let $\{\CM_{1,2},\CM_3\}\to \CN$, $\{\CM_{2,3},\CM_1\}\to \CN$,
$\{\CM_{1,3},\CM_2\}\to \CN$ be chiral pairings such that
the sum of the three maps
$$j_*j^*(\CM_1\boxtimes \CM_2\boxtimes \CM_3)\to \Delta_!(\CN)$$
equals zero. 

Assume also that all of the above modules are
$\on{Jets}(G)_X$-equivariant, and let $\CF$ be a twisted D-module on
$\Gr_{G,X}$.

\begin{lem}   \label{comp with triple}
Under the above circumstances the sum of the three induced maps
$$j_*j^*\Bigl(h^i(\CF\star \CM_1)\boxtimes \CM_2\boxtimes \CM_3\Bigr)\to 
\Delta_!(h^i(\CF\star \CN)\Bigr)$$
is zero.
\end{lem}

\begin{proof}

By the construction of convolution, it suffices to note the
following. Let $\CM_1,\CM_2,\CM_3$, $\CM_{1,2}, \CM_{2,3}, \CM_{1,3},
\CN$ be as above, but let us only assume that $\CM_2,\CM_3$ and
$\CM_{2,3}$ are $\on{Jets}(G)_X$-equivariant. Let $\bg$ be an
$S$-point of $\on{Jets}^{mer}(G)_X$ for some D-scheme $S$.

Then we have the chiral pairings,
\begin{align*}
&\{\CM_1^\bg,\CM_2\}\to \CM_{1,2}^\bg,\,\, \{\CM_1,\CM_3\}\to
\CM_{1,3}^\bg,\,\, \\ &\{\CM^\bg_{1,2},\CM_3\}\to \CN^\bg, \,\,
\{\CM^\bg_{1,3},\CM_2\}\to \CN^\bg, \,\, \{\CM_1^\bg,\CM_{2,3}\}\to
\CN^\bg,
\end{align*}
and the sum of the resulting three maps
$$j_*j^*(\CM_1^\bg\boxtimes \CM_2\boxtimes \CM_3)\to \Delta_!(\CN^\bg)$$
is zero.

\end{proof}

Let now $\CF'_1$, $\CF'_2$ be two chiral $\fD_{G,\kappa}$-modules with
$\CF'_1$ being strongly $\on{Jets}(G)_X$-equivariant on the right and
$\CF'_2$ being strongly $\on{Jets}(G)_X$-equivariant on the left.  Let
$\CF'_1\star \CF'_2$ be their convolution, which is by definition
represented by the complex
$$\fC^\semiinf(L_\fg;\fg,\CF'_1\otimes \CF'_2).$$

Given a chiral pairing $\{\CM',\CM\}\to \CN$ with $\CM'$ being
strongly $\on{Jets}(G)$-equivariant, we on the one hand, obtain 
a chiral pairing of complexes
$$\{\CM', \CF'_2\star \CM\}\to \CF'_2\star \CN,$$
from which we further obtain a chiral pairing of bi-complexes
$$\{\CM', \CF'_1\star (\CF'_2\star \CM)\}\to
\CF'_1\star (\CF'_2\star \CN).$$
On the other hand, from the original pairing we obtain another
pairing of bi-complexes
$$\{\CM', (\CF'_1\star \CF'_2)\star \CM\}\to
(\CF'_1\star \CF'_2)\star \CN.$$

However, by \cite{FG2}, Sect. 22.9.1, the complexes associated to
$$(\CF'_1\star \CF'_2)\star \CM'' \text{ and }
\CF'_1\star (\CF'_2\star \CM'')$$
for $\CM''=\CM$ or $\CM''=\CN$ are isomorphic. 
The next assertion follows from the construction:

\begin{lem}   \label{assoc of fusion}
Under the above circumstances, the diagram of complexes
$$ \CD j_*j^*\Bigl(\CM'\boxtimes \left(\CF'_1\star (\CF'_2\star
\CM)\right)\Bigr) @>>> \Delta_!\Bigl(\CF'_1\star (\CF'_2\star
\CN)\Bigr) \\ @V{\sim}VV @V{\sim}VV \\ j_*j^*\Bigl(\CM'\boxtimes
\left((\CF'_1\star \CF'_2)\star \CM)\right)\Bigr) @>>>
\Delta_!\Bigl((\CF'_1\star \CF'_2)\star \CN)\Bigr) \endCD
$$
commutes.
\end{lem}

\ssec{Proof of \thmref{representability, spherical}}

We shall now prove the following generalization of
\thmref{representability, spherical}.  Let $\CA$ be a chiral algebra
as in \secref{general alg}. Let us assume that the level $\kappa$ is
integral, i.e., the spherical D-modules $\CF_{V,X}$ on $\Gr_{G,X}$ for
$V\in \Rep(\cG)$ make sense.

Assume that for any $V$ as above the convolution $\CF_{V,X}\star \CA$
is acyclic away from degree $0$.

\begin{thm} \label{repr, general}
Let $\CM$ be a strongly $\on{Jets}(G)_X$-equivariant chiral module.

\smallskip

\noindent{\em (1)} The convolution $\CF_{V,X}\star \CM$ is acyclic
away from cohomological degree $0$.

\smallskip

\noindent{\em (2)} Chiral pairings $\{\CF_{V,X}\star \CA,\CM\}\to
\CN$, where $\CN$ is any other chiral $\CA$-module, are in bijection
with maps of chiral $\CA$-modules $\CF_{V,X}\star \CM\to \CN$.

\end{thm}

In this subsection we will prove the first point of the theorem. 

\medskip

The functor $\CM\mapsto \CF_{V,X}\star \CM$ on the derived category of
strongly $\on{Jets}(G)_X$-equivariant chiral modules over $\CA$ is
both left and right adjoint to $\CN\mapsto \CF_{V^*,X}\star
\CN$. Hence, it is enough to show that it is right-exact. Suppose not,
and let $k$ be the maximal integer, for which $h^k(\CF_{V,X}\star
\CM)$ is non-zero for some $\CM\in \CA\mod^{\on{Jets}(G)_X}$.  Then
$k$ is also the maximal integer, for which $h^{-k}(\CF_{V^*,X}\star
\CN)\neq 0$ for $\CN\in \CA\mod^{\on{Jets}(G)_X}$.

Let us choose an object $\CM$ as above that saturates this bound, and
let us denote by $\CN$ the $k$-th cohomology of $\CF_{V,X}\star
\CM$. By adjunction we have a non-zero map $\CM\to
h^{-k}(\CF_{V^*,X}\star \CN)$.

We can represent $\tau^{\geq 0}(\CF_{V,X}\star \CM)$ by a complex
$\CM_1^\bullet$, supported in degrees $\geq 0$, such that the map of
\eqref{global pairing} gives rise to a chiral pairing
$\{\CF_{V,X}\star \CA,\CM\}\to \CM_1^\bullet$. Moreover, we can
represent $\CN$ by a complex $\CN^\bullet$, also supported in degrees
$\geq 0$, such that the map $\CF_{V,X}\star \CM\to \CN[-k]$ is
represented by a map of complexes $\CM_1^\bullet\to \CN^\bullet[-k]$.

Consider the diagram of complexes:
\begin{equation} \label{one comp}
\CD & & \Delta_!\Bigl(h^0(\CF_{V^*,X}\star \CN^\bullet[-k])\Bigr) \\ &
& @AAA \\ j_*j^*\Bigl(h^0(\CF_{V^*,X}\star (\CF_{V,X}\star
\CA))\boxtimes \CM\Bigr) @>>> \Delta_!\Bigl(h^0(\CF_{V^*,X}\star
\CM^\bullet_1)\Bigr) \\ @AAA @AAA \\ j_*j^*(\CA\boxtimes \CM) @>>>
\Delta_!(\CM), \endCD
\end{equation}
where the left vertical arrow comes from the map $\delta_{1,\Gr_G}\to
\CF_{V^*,X}\star \CF_{V,X}$, the lower right vertical arrow comes from
$$\CM\to \CF_{V^*,X}\star \CF_{V,X}\star \CM\to \CF_{V^*,X}\star \CM_1,$$
and the upper horizontal arrow comes by functoriality from 
\propconstrref{global}. This diagram is commutative by \lemref{assoc
  of fusion}.

The composed arrow from the lower left corner to the upper right
corner of the above diagram vanishes, since the composition
$$j_*j^*\Bigl((\CF_{V,X}\star \CA)\boxtimes \CM\Bigr)\to
\Delta_!(\CM^\bullet_1)\to \CN^\bullet[-k]$$ is evidently equal to
$0$.

This is a contradiction, since the lower horizontal arrow is
surjective, and the composition
$$\Delta_!(\CM)\to \Delta_!\Bigl(h^0(\CF_{V^*,X}\star \CM_1)\Bigr)\to
\Delta_!\Bigl(h^0(\CF_{V^*,X}\star \CN^\bullet[-k])\Bigr)$$
equals the map $\CM\to h^{-k}(\CF_{V^*,X}\star \CN)$ above.

\ssec{Proof of \thmref{repr, general}(2)}

The canonical map of \eqref{global, cohomologies} assigns to every map
of chiral modules $\CF_{V,X}\star \CM\to \CN$ a chiral pairing
$\phi:\{\CF_{V,X}\star \CA,\CM\}\to \CF_{V,X}\star \CM$.

Let us construct a map in the opposite direction. By
\propref{equivariance preserved}, it is sufficient to consider the
case when $\CN$ is also $\on{Jets}(G)_X$-equivariant.  For $\phi$ as
above, consider the chiral pairing
$$j_*j^*\Bigl((\CF_{V^*,X}\star \CF_{V,X}\star \CA)\boxtimes
\CM\Bigr)\to \Delta_!(\CF_{V^*,X}\star \CN),$$ obtained from
\propconstrref{global}. Using the canonical map $\CA\to
\CF_{V^*,X}\star \CF_{V,X}\star \CA$ we thus obtain a chiral pairing
$\{\CA,\CM\}\to \CF_{V^*,X}\star \CN$. By \lemref{A is unit}, the
latter gives rise to a map of chiral modules $\CM\to \CF_{V^*,X}\star
\CN$.  By adjunction, we obtain a map $\psi:\CF_{V,X}\star \CM\to
\CN$, as required.

The fact that, if the initial chiral pairing $\phi$ came from a map
$\CF_{V,X}\star \CM\to \CN$, then the resulting map $\psi$ equals the
initial one, follows from \eqref{one comp}.

Let us start with a map $\phi$, and show that the pairing
$\phi':\{\CF_{V,X}\star \CA,\CM\}\to \CM$ obtained from the
corresponding $\psi$, equals the initial $\phi$.

\medskip

Consider the three maps
$$j_*j^*\Bigl(\CA\boxtimes (\CF_{V,X}\star \CA)\boxtimes \CM\Bigr)\to
\Delta_!(\CN),$$ obtained from the map $\phi$. By the definition of
chiral pairings, their sum equals $0$. Using \propconstrref{global},
we obtain three maps
$$j_*j^*\Bigl(\CA\boxtimes (\CF_{V^*,X}\star \CF_{V,X}\star
\CA)\boxtimes \CM\Bigr)\to \Delta_!(\CF_{V^*,X}\star \CN),$$ whose sum
is still equal to $0$, by \lemref{comp with triple}.

Furthermore, using \propconstrref{global} again, we obtain three maps
$$j_*j^*\Bigl((\CF_{V,X}\star \CA)\boxtimes (\CF_{V^*,X}\star
\CF_{V,X}\star \CA) \boxtimes \CM\Bigr)\to \Delta_!(\CF_{V,X}\star
\CF_{V^*,X}\star \CN).$$ Composing these maps with $\CA\to
\CF_{V^*,X}\star \CF_{V,X}\star \CA$ and $\CF_{V,X}\star
\CF_{V^*,X}\star \CN\to \CN$ we obtain three maps
\begin{equation} \label{second comp}
j_*j^*\Bigl((\CF_{V,X}\star \CA)\boxtimes \CA\boxtimes \CM\Bigr)\to
\Delta_!(\CN),
\end{equation}
that sum up to zero.

\medskip

Let us calculate the resulting maps explicitly. It is easy to see that
the first of these maps, namely, the one obtained by first fusing the
$x_1$ and $x_2$ coordinates equals
$$j_*j^*\Bigl((\CF_{V,X}\star \CA)\boxtimes \CA\boxtimes \CM\Bigr) \to
\Delta_{x_1=x_2}{}_!\Bigl(j_*j^*((\CF_{V,X}\star \CA)\boxtimes
\CM)\Bigr)\overset{\phi}\to \Delta_!(\CN),$$ where the first arrow is
given by the chiral action of $\CA$ on $\CF_{V,X}\star \CA$.

The second map, namely, the one obtained by first fusing the $x_2$ and
$x_3$ coordinates equals
$$j_*j^*\Bigl((\CF_{V,X}\star \CA)\boxtimes \CA\boxtimes \CM\Bigr) \to
\Delta_{x_2=x_3}{}_!  \Bigl(j_*j^*((\CF_{V,X}\star \CA)\boxtimes
\CM)\Bigr)\overset{\phi'}\to \Delta_!(\CN),$$ where the first arrow is
the chiral action of $\CA$ on $\CM$.

The third arrow, by construction, factors through a D-module supported
on the diagonal $x_1=x_3$.

\medskip

Let us now compare the three maps of \eqref{second comp} with the
three maps between the same objects that correspond to the initial
chiral pairing $\phi$, and subtract one from another. 

We have two equal maps
$$j_*j^*\Bigl((\CF_{V,X}\star \CA)\boxtimes \CA\boxtimes \CM\Bigr) \to
\Delta_!(\CN),$$ such that one factors as
$$j_*j^*\Bigl((\CF_{V,X}\star \CA)\boxtimes \CA\boxtimes \CM\Bigr)
\twoheadrightarrow \Delta_{x_2=x_3}{}_!  \Bigl(j_*j^*((\CF_{V,X}\star
\CA)\boxtimes \CM)\Bigr) \overset{\phi'-\phi}\longrightarrow
\Delta_!(\CN),$$ and the other through a D-module, supported on the
diagonal $x_1=x_3$.  this implies that both maps are in fact $0$, and
in particular, $\phi=\phi'$.

\bigskip

\bigskip

\centerline{\bf {\large Part II: The Iwahori case}}

\bigskip

In this part of the paper we consider convolutions of the central
sheaves on the affine flag variety with objects of the category of
$I$-equivariant $\hg_\crit$-modules.

\section{Convolution with central sheaves}     \label{sect 5}

\ssec{Opers with nilpotent singularities}

{}From now on we shall fix a point $x\in X$. Let $\lambda$ be an
integral weight such that $\lambda+\rho$ is dominant. Recall that to
such $\lambda$ in \cite{FG2}, Sect. 7.6, we have attached a
subscheme
$$\Spec(\fZ_\fg^{\lambda,\nilp})\subset \Spec(\fZ_\fg).$$
In terms of the isomorphism \eqref{isom with opers} between
$\Spec(\fZ_\fg)$ and the ind-scheme of $\check\fg$-opers on
$\D_x^\times$, the subscheme $\Spec(\fZ_\fg^{\lambda,\nilp})$
corresponds to opers with a regular singularity and residue
$\varpi(-\lambda-\rho)$, see \cite{FG2}, Sect. 2.9.

According to \cite{FG2}, Sect. 7.6, if a weight $\mu$ is of the form
$w(\lambda+\rho)-\rho$ for some $w\in W$, then the support of the
$\hg_\crit$-module $\BM^\mu:=\on{Ind}^{\hg_\crit}_{\fg[[t]]}(M^\mu)$,
where $M^\mu$ denotes the Verma module of highest weight $\mu$ over
$\fg$, is contained (and in fact equal to)
$\Spec(\fZ_\fg^{\lambda,\nilp})$.

\medskip

Let $\fz_\fg^{\lambda,\nilp}$ be the modification of the D-algebra
$\fz_\fg$ at $x$, corresponding to $\fZ_\fg^{\lambda,\nilp}$. From
\cite{FG2}, Sect. 2.9, we obtain that for every $V\in
\Rep(\cG)$, the module $\CV_{X-x}$ extends naturally to a free
$\fZ_\fg^{\lambda,\nilp}$-module, such that the connection has a
pole of order $\leq 1$ at $x$ and a nilpotent monodromy.

\medskip

Let $\wt\fZ_\fg^{\lambda,\nilp}$ be the completion of
$\fZ_\fg$ with respect to the ideal of
$\fZ_\fg^{\lambda,\nilp}$. Let $\Spec(\fZ_\fg^{\int,\nilp})$
(resp., $\Spec(\wt\fZ_\fg^{\int,\nilp})$) be the sub-ind scheme of
$\Spec(\fZ_\fg)$ defined to the disjoint union
$\underset{\lambda}\sqcup\, \Spec(\fZ_\fg^{\lambda,\nilp})$
(resp., $\underset{\lambda}\sqcup\,
\Spec(\wt\fZ_\fg^{\lambda,\nilp})$).

Applying \lemref{extension of chiral modules}, we obtain that for each
$V\in \Rep(\cG)$, the module $\CV_{X-x}$ gives rise to a vector bundle
$\CV_{\wt\fZ_\fg^{\int,\nilp}}$ over
$\Spec(\wt\fZ_\fg^{\int,\nilp})$, equipped with a nilpotent
endomorphism, which we will denote by
$N_{\CV_{\wt\fZ_\fg^{\int,\nilp}}}$. Moreover, for $U\simeq
V\otimes W$ we have an isomorphism
$$\CU_{\wt\fZ_\fg^{\int,\nilp}} \simeq
\CV_{\wt\fZ_\fg^{\int,\nilp}}
\underset{\wt\fZ_\fg^{\int,\nilp}} \otimes
\CW_{\wt\fZ_\fg^{\int,\nilp}},$$ so that
$N_{\CU_{\wt\fZ_\fg^{\int,\nilp}}} =
N_{\CV_{\wt\fZ_\fg^{\int,\nilp}}}+
N_{\CW_{\wt\fZ_\fg^{\int,\nilp}}}$.

\ssec{Statement of the theorem}  \label{statement of Iwahori theorem}

Let $I\subset G\ppart$ be the Iwahori subgroup, i.e., the preimage of
$B\subset G$ under the evaluation map $G[[t]]\to G$. Following the
notation of \cite{FG2}, we shall denote by $\hg_\kappa\mod^I$ the
category of $I$-integrable representations of $\hg$ at the level
$\kappa$.

\medskip

The above estimate on the support of the modules $\BM^\mu$ implies:
\begin{lem}
The support of every object $\CM\in \hg_\crit\mod^I$ over
$\Spec(\fZ_\fg)$ is contained in
$\Spec(\wt\fZ_\fg^{\int,\nilp})$.
\end{lem}

Thus, for every $\CM\in \hg_\crit\mod^I$ and $V\in \Rep(\cG)$ we can
functorially attach another object of $\hg_\crit\mod^I$:
$$\CV_{\wt\fZ_\fg^{\int,\nilp}}
\underset{\wt\fZ_\fg^{\int,\nilp}}\otimes \CM,$$ which carries a
nilpotent endomorphism $N_{\CV_{\wt\fZ_\fg^{\int,\nilp}}}$.

\medskip

Let $\Fl_G=G\ppart/I$ be the affine flag scheme of $G$. For a level
$\kappa$ we shall denote by $\on{D}(\Fl_G)_\kappa\mod$ the category of
$\kappa$-twisted right D-modules on $\Fl_G$. By
$\on{D}(\Fl_G)_\kappa\mod^I$ we shall denote the corresponding
category of $I$-equivariant objects in
$\on{D}(\Fl_G)_\kappa\mod$. When $\kappa$ is integral (e.g., critical)
we will sometimes identify $\on{D}(\Fl_G)_\kappa\mod$ with the usual
category of D-modules by means of the tensor product with the
corresponding line bundle.

\medskip

Given an object $\CF\in \on{D}(\Fl_G)_\kappa\mod$ and $\CM\in
\fg_\kappa\mod^I$ we can form their convolution, denoted
$\CF\underset{I}\star \CM$, which is an object of
$D(\fg_\kappa\mod)$. When no confusion is likely to occur, we will
omit the subscript $I$ from $\underset{I}\star$.

\medskip

Let us recall from \cite{Ga} that to every object $V\in \Rep(\cG)$
there corresponds an object $\CZ_V\in \on{D}(\Fl_G)\mod^I$ called a
central sheaf. Each central sheaf is endowed with a functorial
endomorphism $N_V$. The construction of D-modules $\CZ_V$ will be
reviewed in some detail in the sequel. Slightly abusing the notation,
we shall denote by the same symbol $\CZ_V$ the corresponding
$I$-equivariant object of $\on{D}(\Fl_G)_\crit\mod^I$.

Our main result is the following:

\begin{thm} \label{main}
For every $\CM\in \hg_\crit\mod^I$ and $V\in \Rep(\cG)$ the
convolution $\CZ_V\star \CM$ is acyclic away from cohomological degree
$0$, and we have a canonical isomorphism
$$\fs_V:\CZ_V\underset{I}\star \CM \simeq
\CV_{\wt\fZ^{\int,\nilp}_{\fg,x}}\underset{\wt\fZ^{\int,\nilp}_{\fg,x}}
\otimes \CM,$$ such that the endomorphism induced by $N_V$ on the LHS
goes over to the endomorphism, induced by
$N_{\CV_{\wt\fZ_\fg^{\int,\nilp}}}$ on the RHS. This system of
isomorphisms is compatible with tensor products of
$\cG$-representations in the same sense as in \thmref{main, spherical}
\end{thm}

\ssec{Spherical case, revisited}  \label{another proof mon-free}

Let us note that \thmref{main}, whose prove is parallel to, but independent
of, the proof of \thmref{main, spherical}, implies the latter theorem.
Indeed, let $\fp$ denote the natural projection
$\Fl_G\to \Gr_G$; then by \cite{Ga}, there is a canonical isomorphism
$$\fp_!(\CZ_V)\simeq \CF_V.$$
For an object $\CM\in \hg_\crit\mod^{G[[t]]}$ we have:
$$\CZ_V\underset{I}\star \CM\simeq \CF_V\underset{G[[t]]}\star \CM,$$
and the assertion follows from the fact that the restriction of
$\CV_{\wt\fZ^{\int,\nilp}_{\fg,x}}$ to
$\Spec(\wt\fZ^{\int,\reg}_{\fg})\subset
\Spec(\wt\fZ^{\int,\nilp}_{\fg,x})$ identifies canonically with
$\CV_{\wt\fz^{\int,\reg}_{\fg,x}}$.

In particular, since the nilpotent endomorphism that $N_V$ induces on
$\CF_V$ is zero, we obtain the assertion of \corref{monodromy-free}.
This, in turn, gives an alternative proof of \lemref{sup of sph}, as 
was promised earlier. Let us now prove \lemref{sup of ind}:

\begin{proof}

One the one hand, by \cite{FG2}, Corollary 7.6.2, the support of $\BV^\mu$ over $\Spec(\fZ_\fg)$
is contained in the subscheme $\Spec(\fZ_\fg^{\mu,\nilp})$. 
On the other hand, by \corref{monodromy-free}, which was proved independently,
the support of $\BV^\mu$ is contained in 
the ind-subscheme $\Spec(\fZ_{\fg}^{\on{m.f.}})$. The assertion of the lemma follows
now from the fact that
$$\Spec(\fZ_\fg^{\mu,\nilp})\cap \Spec(\fZ_{\fg}^{\on{m.f.}})=\Spec(\fZ_\fg^{\mu,\reg}),$$
(see \cite{FG2}, Sect. 2.9).

\end{proof}

\ssec{The case of differential operators}

Parallel to the spherical situation, let us consider a particular case
of the above theorem, corresponding to differential
operators. Consider the object of $\delta_{I,G\ppart}\in
\fD_{G,\crit}\mod_x$, corresponding to distributions on $I$.

In other words, if $\pi_{\Fl}$ denotes the projection $G\ppart\to
\Fl_G$, then
$$\delta_{I,G\ppart}\simeq
\Gamma(G\ppart,\pi_{\Fl}^*(\delta_{1,\Fl_G})),$$ where
$\delta_{1,\Fl_G}$ is the $\delta$-function at $1\in \Fl_G$. As a
$\hg_\crit$-module, $\delta_{I,G\ppart}$ can be described as
$\on{Ind}^{\hg_\crit}_{\fg[[t]]}(\CO_I)$. We have:

\begin{thm} \label{main for diff op}
We have a canonical isomorphism of $\hg_\crit$-bimodules
$$\Gamma(G\ppart,\pi_{\Fl_G}^*(\CZ_V))\simeq
\CV_{\wt\fZ^{\int,\nilp}_{\fg,x}}\underset{\wt\fZ^{\int,\nilp}_{\fg,x}}
\otimes \delta_{I,G\ppart},$$ which intertwines the endomorphism $N_V$
on the LHS with $N_{\CV_{\wt\fZ_\fg^{\int,\nilp}}}$ on the RHS.
This system of isomorphisms is compatible with tensor products of
$\cG$-representations.
\end{thm}

The same argument as in the spherical case shows that \thmref{main for
diff op} actually implies \thmref{main}.

\section{Fusion in the Iwahori case}    \label{sect 6}

Our strategy of proof of \thmref{main} will be parallel to that of
\thmref{main, spherical}.  Namely, we will construct a functor on the
category of $\hg_\crit$-modules, and show that it is represented by
both sides of the isomorphism stated in the theorem.

An additional ingredient in the present situation is that we will have
to twist the chiral $\CA_{\fg,\kappa}$-modules
$\Gamma(\Gr_{G,X-x},\CF_{V,X-x})$ by certain local systems on the
punctured curve.

\ssec{Nilpotent local systems}

Let us fix a family of local systems $\CE_n$, $n\in \BN$, defined on a
punctured Zariski neighborhood of $x\in X$, each $\CE_n$ being an
$n$-fold extension of the trivial local system (i.e., $\CE_1$ is {\it
the} trivial local system), endowed with a nilpotent endomorphism
$N_{\CE}$ of order $n-1$, and such that we are given a compatible
system of surjections
\begin{equation} \label{change of n}
\CE_n\twoheadrightarrow \CE_m, 
\end{equation}
defined for every $m\leq n$ and compatible with the action of
$N_{\CE}$.  Then the action of $N_\CE$ gives rise to a system of exact
sequences:
$$0\to \CE_m\to \CE_n\overset{N^m_{\CE}}\to \CE_n\to \CE_{n-m}\to 0.$$

For example, if we choose a local coordinate $t$ near $x$, such a
system can be obtained as follows. Let $E_n$ be the standard
$n$-dimensional Jordan block, i.e., the vector space $\BC[s]/s^{n-1}$ with
the nilpotent operator $N_E$ equal to the multiplication by $s$.  Set
$\CE_n:=E_n\otimes \CO_{X-x}$, with $t\nabla_t$ acting as
$N_E+t\partial_t$.

Let us denote by $\CE_{n,x}$ the fiber at $x$ of the extension of
$\CE_n$, given by \lemref{extension of chiral modules}. This is an
$n$-dimensional vector space, endowed with a nilpotent operator of order
$n-1$, which coincides with the one induced by $N_\CE$.  We shall fix
a system of identifications $\CE_{n,x} \simeq E_n$, such that $N_\CE$
goes over to $N_E$, and which is compatible with the morphisms
\eqref{change of n}.

\ssec{A canonical chiral pairing}

Recall now that if $\CA$ is any chiral algebra and $\CM$ is a chiral
module over it, and $\CE$ is any D-module on $X$, then the D-module
$\CM\otimes \CE$ is naturally a chiral module over $\CA$.

\begin{propconstr} \label{can pairing with monod}
Let $\kappa$ be any integral level. Then for any $\CM\in
\hg_\kappa\mod^I$ and $n\in \BN$ there exists a canonical chiral
pairing of $\fD_{G,\crit}$-modules:
$$j_x{}_*j_x^*\biggl(\Gamma\Bigl(G\ppart,\pi^*(\CF_V) \Bigr)\otimes
\CE_n\biggr) \otimes \delta_{I,G\ppart}\to
i_x{}_!\biggl(\Gamma\Bigl(G\ppart,\pi_{\Fl}^*(\CZ_V)
\Bigr)\biggr)_{N_V^n},$$ such that for $n\geq m$ the diagram
$$ \CD j_x{}_*j_x^*\biggl(\Gamma \Bigl(G\ppart,\pi^*(\CF_V)\Bigr)
\otimes \CE_n\biggr)\otimes \delta_{I,G\ppart} @>>>
i_x{}_!(\Gamma\Bigl(G\ppart,\pi_{\Fl}^*(\CZ_V)\Bigr)\biggr)_{N_V^n} \\
@VVV @VVV \\
j_x{}_*j_x^*\biggl(\Gamma\Bigl(G\ppart,\pi^*(\CF_V)\Bigr)\otimes
\CE_m\biggr)\otimes \delta_{I,G\ppart} @>>>
i_x{}_!\biggl(\Gamma\Bigl(G\ppart,\pi_{\Fl}^*(\CZ_V) \Bigr)
\biggr)_{N_V^m} \endCD
$$
is commutative.
\end{propconstr}

Before giving the proof we need to review the construction of the
D-modules $\CZ_V$. Consider the ind-scheme $\Fl_{G,X}$ over $X$, whose
fiber over $x_1\in X$ is the set of triples $(\CP_G,\beta, \alpha)$,
where $\beta$ is a trivialization of $\CP_G$ on $X-\{x,x_1\}$ and
$\alpha$ is a reduction of the fiber $\CP_{G,x}$ at $x$ to $B$. Note
that $\Fl_{G,X}$ is the same as $\Gr_{G,X;K}$ for $K=I$.

The preimage of $X-x$ in $\Fl_{G,X}$, denoted $\Fl_{G,X-x}$, is
isomorphic to $\Gr_{G,X-x}\times \Fl_G$, whereas the preimage of $x\in
X$ is isomorphic to $\Fl_G$.  More generally, we will consider the
ind-scheme $\Fl_{G,X^n}$ over $X^n$, whose fiber over
$(x_1,...,x_n)\in X^n$ is the set of triples $(\CP_G,\beta, \alpha)$,
where $\beta$ now is a trivialization defined away from
$x_1\cup...\cup x_n\cup x$.

Recall the ind-scheme $\on{Jets}^\mer(G)_{X^n}$ from \secref{Jets n},
and let $\on{Jets}^\mer(G)_{X^{n-1}\times x}$ be the closed
ind-subscheme defined by the condition that the last point in the
$n$-tuple $x_1,...,x_n$ is fixed to be $x$.

Note that we have a natural projection
$\pi_{\Fl}:\on{Jets}^\mer(G)_{X^n\times x}\to \Fl_{G,X^n}$, which
corresponds to remembering the reduction to $B$ at $x$ out of the
trivialization of $\CP_G$ on the formal neighborhood of
$x_1\cup...\cup x_n\cup x$.

\medskip

Given an object $V\in \Rep(\cG)$ and $n\in \BN$, we consider
$(\CF_{V,X-x}\otimes \CE_n)\boxtimes \delta_{1,\Fl_G}$ as a D-module
on $\Fl_{G,X-x}$.  For $n\in \BN$ consider the fiber over $x$ of its
intermediate extension onto the entire $\Fl_{G,X}$, i.e.,
$$\CZ_{V,n}:=i_x^! j_x{}_{!*}\Bigl((\CF_{V,X-x}\otimes \CE_n)\boxtimes
\delta_{1,\Fl_G}\Bigr).$$

This is a D-module on $\Fl_G$, endowed with an action of $N_\CE$, by
the transport of structure. The following summarizes the main
construction of \cite{Ga}:

\medskip

\begin{itemize}

\item
For $n$ large enough and $n'\geq n$ the maps $\CZ_{V,n'}\to \CZ_{V,n}$
are isomorphisms.

\item
The D-module $\CZ_V$ is isomorphic to $\CZ_{V,n}$ for all $n$ large
enough.  Under this identification, the endomorphism $N_V$ equals the
one induced by $N_{\CE}$.

\item
For any other $n''$, $\CZ_{V,n''}\simeq (\CZ_V)_{N^{n''}_V}$.

\end{itemize}

In fact, the above is a paraphrase of the fact that $\CZ_V$ is
obtained by applying the functor of unipotent nearby cycles
$\Psi^{un}$ to the D-module $(\CF_{V,X-x}\boxtimes \delta_{1,\Fl_G})$
on $\Fl_{G,X-x}$.

Note that if $\CE'$ is any local system on $X-x$ with a nilpotent
monodromy around $x$, then the formalism of nearby cycles, developed
in \cite{Be} implies that
\begin{equation} \label{extension against local system}
i_x^! j_x{}_{!*}\Bigl((\CF_{V,X-x}\otimes \CE')\boxtimes
\delta_{1,\Fl_G}\Bigr)[1]\simeq \Bigl(\CZ_V\otimes
\CE'_x\Bigr)_{N_V+N_{\CE'}},
\end{equation}
where $\CE'_x$ is the fiber of the extension of $\CE'$ across $x$,
given by \lemref{extension of chiral modules}, and $N_{\CE'}$ is its
canonical nilpotent endomorphism.

\ssec{Proof of \propconstrref{can pairing with monod}}

Consider the canonical map of twisted D-modules on $\Fl_{G,X}$
\begin{equation} \label{map downstairs, Iwahori}
j_x{}_*\Bigl((\CF_{V,X-x}\otimes \CE_n)\boxtimes
\delta_{1,\Fl_G}\Bigr)\to i_x{}_!(\CZ_{V,n}).
\end{equation}

By taking the pull-back of \eqref{map downstairs, Iwahori} under
$\pi_{\Fl}:\on{Jets}^\mer(G)_{X\times x}\to \Fl_{G,X}$, we obtain a
map of D-modules on $X$:
\begin{equation} \label{map upstairs, Iwahori}
j_x{}_*j_x^*\biggl(\Gamma\Bigl(\on{Jets}^\mer(G)_X,\pi^*(\CF_{V,X})\Bigr)
\otimes \CE_n\biggr) \otimes \delta_{I,G\ppart} \to 
i_x{}_!\biggl(\Gamma\Bigl(G\ppart,\pi_{\Fl}^*(\CZ_V)\Bigr)\biggr)_{N_V^n},
\end{equation}
as was stated in \propconstrref{can pairing with monod}. Thus, in
order to finish the proof, we need to show that the above map respects
the action of $\fD_{G,\kappa}$, i.e., that it is indeed a chiral
pairing.

\medskip

Consider the ind-scheme $\Fl_{G,X^2}$ over $X^2$; let
$\Fl_{G,(X-x)^2-\Delta_{X-x}}$ denote its open subscheme equal to the
preimage of the corresponding open subscheme in $X^2$; let $j$ denote
the corresponding open embedding. Let $\Delta$ be the embedding of
$\Fl_{G,x\times x}\simeq \Fl_G$.

We have an isomorphism:
$$\Fl_{G,(X-x)^2-\Delta_{X-x}}\simeq \Bigl((\Gr_{G,X-x}\times
\Gr_{G,X-x})\underset{(X-x)^2}\times ((X-x)^2-\Delta_{X-x})\Bigr)
\times \Fl_G.$$

As in the case of $\Gr_{G,X^n}$, we have a map ${\bf 1}_{1,1}:X\times
\Fl_{G,X}\to \Fl_{G,X^2}$, and consider the twisted D-module on
$\Fl_{G,X^2}$ equal to
$$({\bf 1}_{1,1})_!\biggl(\omega_X\boxtimes
\Bigl(j_x{}_{!*}(\CF_{V,X-x}\otimes \CE_n)\Bigr)\biggr).$$ It gives
rise to three maps
\begin{equation} \label{three maps down flags}
j_*j^*\Bigl(\delta_{1,\Gr_{G,X}}\boxtimes (\CF_{V,X-x}\otimes
\CE_n)\boxtimes \delta_{1,\Fl_G}\Bigr)\to \Delta_!(\CZ_V)_{N_V^n},
\end{equation}
which sum up to zero.

Pulling back to the two sides of \eqref{three maps down flags} to
$\on{Jets}^\mer(G)_{X^2\times x}$ and taking the (quasi-coherent)
direct image onto $X^2$ we obtain three maps
$$j_*j^*\biggl(\fD_{G,\kappa}\boxtimes
\Bigl(\Gamma(\on{Jets}^\mer(G)_X,\pi^*(\CF_{V,X})) \otimes
\CE_n\Bigr)\boxtimes \delta_{I,G\ppart}\biggr)\to
\Delta_!\biggl(\Gamma\Bigl(G\ppart,\pi_{\Fl}^*(\CZ_V)\Bigr)\biggr)_{N_V^n}.$$

{}From \propref{action as fusion}, we obtain that these three maps are
equal to those that appear in the definition of chiral pairings.

\ssec{Chiral pairings with other representations}

As in the spherical situation, from the chiral pairing given by
\propconstrref{can pairing with monod}, we obtain a chiral pairing
over $\CA_{\fg,\kappa}$:
\begin{equation}  \label{another can pairing, Iwahori}
j_x{}_*j_x^*\biggl(\Bigl(\Gamma(\Gr_{G,X},\CF_{V,X} \Bigr)\otimes
\CE_n\biggr)\otimes \delta_{I,G\ppart}\to
i_x{}_!\biggl(\Gamma\Bigl(G\ppart,\pi_{\Fl}^*(\CZ_V)
\Bigr)\biggr)_{N_V^n},
\end{equation}
which commutes with the action of $\hg_{\kappa'}$ on
$\delta_{I,G\ppart}$ and
$\Gamma\Bigl(G\ppart,\pi_{\Fl}^*(\CZ_V)\Bigr)$, given by $\fr$.

Given an $I$-integrable $\fg_\kappa$-module $\CM$ we obtain a chiral
pairing of complexes of $\CA_{\fg,\kappa}$-modules:
$$\{\Gamma(\Gr_{G,X},\CF_{V,X})\otimes \CE_n,
\fC^\semiinf(\fg\ppart;\fh,\delta_{I,G\ppart}\otimes \CM)\} \to
\fC^\semiinf\Bigl(\fg\ppart;\fh,\Gamma(G\ppart,\pi_{\Fl}^*(\CZ_V))\otimes
\CM\Bigr)_{N_V^n}.$$ In particular, by taking $n>>0$ so that
$N_V^n=0$, we obtain a chiral pairing
\begin{equation} \label{can pairing, Iwahori}
\{\Gamma(\Gr_{G,X-x},\CF_{V,X-x})\otimes \CE_n,\CM\}\to
h^0(\CZ_V\underset{I}\star \CM),
\end{equation}
such that the action of $N_{\CE}$ on LHS corresponds to the action of
$N_V$ on the RHS.

More generally, for an arbitrary local system $\CE'$ with a nilpotent
monodromy around $x$, we obtain a chiral pairing
\begin{equation} \label{can pairing against loc sys}
\{\Gamma(\Gr_{G,X-x},\CF_{V,X-x})\otimes \CE',\CM\}\to 
\Bigl(h^0(\CZ_V\underset{I}\star \CM)\otimes \CE'_x\Bigr)_{N_V+N_{\CE'}}.
\end{equation}

Parallel to \thmref{representability, spherical}, we will prove:

\begin{thm} \label{representability, Iwahori}  
Assume that $\kappa$ is non-positive, Then:

\smallskip

\noindent{\em (1)} For every $\CM\in \hg_\kappa\mod^I$, the
convolution $\CZ_V\underset{I}\star \CM$ is acyclic away from
cohomological degree $0$.

\smallskip

\noindent{\em (2)} The functor on $\hg_\kappa\mod$ that sends $\CN$ to
the set of chiral pairings
$$\{\Gamma(\Gr_{G,X},\CF_{V,X})\otimes \CE',\CM\}\to\CN,$$ is
representable by $\Bigl((\CZ_V\underset{I}\star \CM)\otimes
\CE'_x\Bigr)_{N_V+N_{\CE'}}$.

\end{thm}

We shall now proceed with the proof of \thmref{main}, which does not
rely on \thmref{representability, Iwahori}.  We will need to establish
the following generalization of \propref{comp with ten products}:

For $V,W\in \Rep(\cG)$, recall from \cite{Ga} that there exists a
natural isomorphism
\begin{equation} \label{tensor center}
\CZ_V\underset{I}\star \CZ_W\simeq \CZ_{V\otimes W},
\end{equation}
such that the endomorphism $N_V+N_W$ on the LHS goes over to
the endomorphism $N_{V\otimes W}$.

\medskip

Let $n,m\in \BN$ be such that $\CZ_V\to (\CZ_V)_{N_V^n}$ and $\CZ_W\to
(\CZ_W)_{N_V^n}$ are isomorphisms. Let $\CM$ be an object of
$\hg_\crit\mod^I$, such that all convolutions
$\CZ_{V'}\star \CM$, $V'\in \Rep(\cG)$ are acyclic away from
cohomological degree $0$.

In this case, as in \propref{comp with ten products}, we have three
maps
\begin{equation} \label{three maps, Iw}
j_*j^*\Bigl((\Gamma(\Gr_{G,X},\CF_{V,X})\otimes \CE_n)\boxtimes
(\Gamma(\Gr_{G,X},\CF_{W,X})\otimes \CE_m)\boxtimes \CM\Bigr)\to 
i_x{}_!(\CZ_{V\otimes W}\underset{I}\star \CM),
\end{equation}
defined as follows.

The first map is the composition
\begin{align*}
&j_*j^*\biggl(\Bigl(\Gamma(\Gr_{G,X},\CF_{V,X})\otimes
\CE_n\Bigr)\boxtimes \Bigl(\Gamma(\Gr_{G,X},\CF_{W,X})\otimes
\CE_m\Bigr)\boxtimes \CM\biggr)\to \\ &
\Delta_{x_2=x_3}{}_!\biggl(j_*j^*\Bigl(\Gamma(\Gr_{G,X},\CF_{V,X})\otimes
\CE_n) \boxtimes (\CZ_W\star \CM)\Bigr)\biggr)\to
\Delta_!\Bigl(\CZ_V\underset{I}\star \CZ_W\underset{I}\star \CM\Bigr).
\end{align*}

The second map is the negative of the one defined by interchanging the
roles of $V$ and $W$.  To define the third map note that from
\eqref{extension against local system} we obtain a map of D-modules on
$\Fl_{G,X}$:
\begin{equation} \label{n & m}
j_x{}_*j_x^*\biggl(\Bigl(\Gamma(\Gr_{G,X},\CF_{V\otimes W})\otimes
\CE_n\otimes \CE_m\Bigr) \boxtimes \delta_{1,\Fl_G}\biggr)\to
\Delta_!(\CZ_V\underset{I}\star \CZ_W),
\end{equation}
which comes from the projection
\begin{align*} 
&\Bigl(\CZ_{V\otimes W}\otimes E_n\otimes E_m\Bigr) _{N_{V\otimes
W}+N_\CE|_{E_n}+N_\CE|_{E_m}}\simeq \Bigl((\CZ_V\otimes E_n)\star
(\CZ_W\otimes E_m)\Bigr) _{N_{V}+N_{W}+N_\CE|_{E_n}+N_\CE|_{E_m}} \\
&\twoheadrightarrow (\CZ_V\otimes E_n)_{N_V+N_\CE|_{E_n}}\star
(\CZ_W\otimes E_m)_{N_W+N_\CE|_{E_m}}\simeq \CZ_V\underset{I}\star
\CZ_W.
\end{align*}

Corresponding to it there is a chiral pairing, defined for every $\CM$
as above:
\begin{equation} \label{n & m, modules}
j_x{}_*j_x^*\biggl(\Bigl(\Gamma(\Gr_{G,X},\CF_{V\otimes W})\otimes
\CE_n\otimes \CE_m\Bigr) \boxtimes \CM\biggr)\biggr)\to
\Delta_!(\CZ_V\underset{I}\star \CZ_W\underset{I}\star \CM).
\end{equation}

The third map in \eqref{three maps, Iw} equals the composition
\begin{align*}
&j_*j^*\biggl(\Bigl(\Gamma(\Gr_{G,X},\CF_{V,X})\otimes
\CE_n\Bigr)\boxtimes \Bigl(\Gamma(\Gr_{G,X},\CF_{W,X})\otimes
\CE_m\Bigr)\boxtimes \CM\biggr)\to \\ &\to \Delta_{x_1=x_2}{}_!
\biggl(j_x{}_*j_x^*\biggl(\Bigl(\Gamma(\Gr_{G,X},\CF_{V\otimes
W})\otimes \CE_n\otimes \CE_m\Bigr) \boxtimes
\CM\biggr)\biggr)
\overset{\text{\eqref{n & m, modules}}}\longrightarrow 
\Delta_!(\CZ_V\underset{I}\star \CZ_W\underset{I}\star \CM).
\end{align*}

\begin{prop} \label{assoc, Iwahori}
The sum of the three maps above is $0$.
\end{prop}

\begin{proof}

Recall that $$\Fl_{G,(X-x)^2-\Delta_{X-x}}\simeq
\Bigl((\Gr_{G,X-x}\times \Gr_{G,X-x})\underset{(X-x)^2}\times
((X-x)^2-\Delta_{X-x})\Bigr) \times \Fl_G,$$ and consider the twisted
D-module
\begin{equation} \label{on square}
j_{!*}\Bigl((\CF_{V,X-x}\otimes \CE_n)\boxtimes 
(\CF_{W,X-x}\otimes \CE_m)\boxtimes \delta_{1,\Fl_G}\Bigr)
\end{equation}
on $\Fl_{G,X^2}$.

\medskip

Let $\Fl_{G,X\times x}$ (resp., $\Fl_{G,x\times X}$,
$\Fl_{G,\Delta_X}$) denote the preimage in $\Fl_{G,X^2}$ of the
corresponding subvariety in $X^2$.  By considering the iterated
version of $\Fl_{G,X^2}$, the following description of the D-module
\eqref{on square} was obtained in \cite{Ga1}:

\begin{itemize}

\item The restriction of the D-module \eqref{on square} to
$\Fl_{G,X\times x}\simeq \Fl_{G,X}$ identifies with
$$\on{ker}
\biggl(j_x{}_*j_x^*\Bigl((\CF_{V,X}\otimes \CE_n)\boxtimes \CZ_W\Bigr)\to 
i_x{}_!(\CZ_V\underset{I}\star \CZ_W)\biggr).$$

\item The restriction of \eqref{on square} to $\Fl_{G,x\times X}\simeq
\Fl_{G,X}$ identifies with
$$\on{ker}
\biggl(j_x{}_*j_x^*\Bigl((\CF_{W,X}\otimes \CE_m)\boxtimes \CZ_V\Bigr)\to 
i_x{}_!(\CZ_W\underset{I}\star \CZ_V)\biggr).$$

\item The restriction of \eqref{on square} to $\Fl_{G,\Delta_X}\simeq
\Fl_{G,X}$ identifies with
$$\on{ker} \biggl(j_x{}_*j_x^*\Bigl((\CF_{V\otimes W,X}\otimes
\CE_n\otimes \CE_m)\boxtimes \delta_{1,\Fl_G}\Bigr)
\overset{\text{\eqref{n & m}}}\longrightarrow i_x{}_!
(\CZ_V\underset{I}\star \CZ_W)\biggr).$$

\end{itemize}

Hence, we obtain three maps 
\begin{equation} \label{three maps downstairs, Iwahori}
j_*j^*\Bigl((\CF_{V,X}\otimes \CE_n)\boxtimes (\CF_{W,X}\otimes
\CE_m)\boxtimes \delta_{1,\Fl_G}\Bigr)\to
\Delta_!(\CZ_V\underset{I}\star \CZ_W),
\end{equation}
whose some equals to zero. By lifting the terms of \eqref{three maps
downstairs, Iwahori} and the corresponding maps by means of
$\pi_{\Fl}:\on{Jets}^\mer(G)_{X^2\times x}\to \Fl_{G,X^2}$, we obtain
three maps
$$j_*j^*\Bigl((\Gamma(\Gr_{G,X},\CF_{V,X})\otimes \CE_n)\boxtimes
(\Gamma(\Gr_{G,X},\CF_{W,X})\otimes \CE_m)\boxtimes
\delta_{I,G\ppart}\Bigr)\to \Delta_!\Bigl(\CZ_{V\otimes W}\underset{I}\star
\delta_{I,G\ppart}\Bigr).$$ We claim that these maps coincide with
those of of \propref{assoc, Iwahori} for
$\CM=\delta_{I,G\ppart}$. This follows from \lemref{action as fusion}
in the same way as in the proof of \propref{comp with ten products}.

The case of a general $\CM$ follows from that of $\delta_{I,G\ppart}$
by the construction of the maps \eqref{can pairing, Iwahori}.

\end{proof}

\ssec{Proof of \thmref{main}}

The proof will be parallel to that of \thmref{main, spherical}. It
suffices to construct the isomorphisms
\begin{equation}  \label{req iso}
\fs_{V}^{-1}:
\CV_{\wt\fZ^{\int,\nilp}_{\fg,x}}\underset{\wt\fZ^{\int,\nilp}_{\fg,x}}
\otimes \CM\simeq \CZ_V\underset{I}\star \CM,
\end{equation}
for $\CM=\delta_{I,G\ppart}$, which commute with the right action of
$\hg_\crit$, which intertwine the actions of
$N_{\CV_{\wt\fZ^{\int,\nilp}_{\fg,x}}}$ and $N_V$, and which are
compatible with tensor products of representations.  Given $V\in
\Rep(\cG)$ and $\CM\in \hg_\crit\mod^I$ will construct the map
$\fs_V^{-1}$ functorially, provided that $\CZ_V\underset{I}\star \CM$
is acyclic away from cohomological degree $0$ for all $V\in
\Rep(\cG)$.

\medskip

Let $\CN_1$ and $\CN_2$ be two $\hg_\crit$-modules, such that the
support of $\CN_1$ over $\Spec(\fZ_\fg)$ is in
$\Spec(\wt\fZ^{\int,\nilp}_{\fg,x})$.  The next assertion is deduced
from \lemref{pairings over center} and \propref{commutative fusion} as
was the case of \thmref{repr of pair from ten product}:

\begin{lem}  \label{com pairing, Iwahori}
Chiral pairings
$$\{(\CV_{\fz_\fg}\underset{\fz_\fg}\otimes
\CA_{\fg,\crit})\underset{\CO_X}\otimes \CE_n,\CN_1\} \to \CN_2$$ are
in bijection with maps of $\hg_\crit$-modules
$$\Bigl(\CV_{\wt\fZ^{\int,\nilp}_{\fg,x}}
\underset{\wt\fZ^{\int,\nilp}_{\fg,x}} \otimes \CN_1\Bigr)
_{N_{\CV_{\wt\fZ^{\int,\nilp}_{\fg,x}}}^n}\to \CN_2.$$
\end{lem}

More generally, for a local system $\CE'$ on $X-x$ with a nilpotent
monodromy around $x$, chiral pairings
$$\{(\CV_{\fz_\fg}\underset{\fz_\fg}\otimes
\CA_{\fg,\crit})\underset{\CO_X}\otimes \CE',\CN_1\} \to \CN_2$$ are
in bijection with maps of $\hg_\crit$-modules
$\Bigl(\Bigl(\CV_{\wt\fZ^{\int,\nilp}_{\fg,x}}
\underset{\wt\fZ^{\int,\nilp}_{\fg,x}} \otimes \CN_1\Bigr) \otimes
\CE'_x\Bigr)_{N_{\CV_{\wt\fZ^{\int,\nilp}_{\fg,x}}}+N_{\CE'}}\to
\CN_2$.

\medskip

Setting $\CN_1=\CM$ and $\CN_2=\CZ_V\underset{I}\star \CM$ and taking
$n>>0$ so that $N_{\CV_{\wt\fZ^{\int,\nilp}_{\fg,x}}}^n=0$, from the
map \eqref{can pairing, Iwahori}, we produce the desired map
$$\fs^{-1}_V:\CV_{\wt\fZ^{\int,\nilp}_{\fg,x}}
\underset{\wt\fZ^{\int,\nilp}_{\fg,x}}\otimes \CM \to
\CZ_V\underset{I}\star \CM.$$

\medskip

We claim that this map is compatible with tensor products of objects
of $\Rep(\cG)$, i.e., that the following three maps coincide:
$$\CV_{\wt\fZ^{\int,\nilp}_{\fg,x}}
\underset{\wt\fZ^{\int,\nilp}_{\fg,x}}\otimes
\CW_{\wt\fZ^{\int,\nilp}_{\fg,x}}
\underset{\wt\fZ^{\int,\nilp}_{\fg,x}}\otimes \CM
\overset{\fs^{-1}_V}\to \CZ_V \underset{I}\star
(\CW_{\wt\fZ^{\int,\nilp}_{\fg,x}}
\underset{\wt\fZ^{\int,\nilp}_{\fg,x}}\otimes \CM)
\overset{\on{id}_{\CZ_V} \underset{I}\star \fs^{-1}_W}\to
\CZ_V\underset{I}\star \CZ_W \underset{I}\star \CM,$$
$$\CV_{\wt\fZ^{\int,\nilp}_{\fg,x}}
\underset{\wt\fZ^{\int,\nilp}_{\fg,x}}\otimes
\CW_{\wt\fZ^{\int,\nilp}_{\fg,x}}\underset{\wt\fZ^{\int,\nilp}_{\fg,x}}
\otimes
\CM \overset{\on{id}_{\CV}\otimes \fs^{-1}_W}\to
\CV_{\wt\fZ^{\int,\nilp}_{\fg,x}}\underset{\wt\fZ^{\int,\nilp}_{\fg,x}}
\otimes
(\CZ_W\underset{I}\star \CM)\overset{\fs^{-1}_V}\longrightarrow
\CZ_V\underset{I}\star \CZ_W\underset{I}\star \CM$$ and
$$\CV_{\wt\fZ^{\int,\nilp}_{\fg,x}}\underset{\wt\fZ^{\int,\nilp}_{\fg,x}}
\otimes
\CW_{\wt\fZ^{\int,\nilp}_{\fg,x}}\underset{\wt\fZ^{\int,\nilp}_{\fg,x}}
\otimes \CM \overset{\fs^{-1}_{V\otimes W}}\to \CZ_{V\otimes
W}\underset{I}\star \CM \simeq \CZ_V\underset{I}\star
\CZ_W\underset{I}\star \CM.$$ This follows from \propref{assoc,
Iwahori} and \lemref{add assoc} as in the proof of \propref{map is
assoc}.

\medskip

Finally, the fact that the maps $\fs^{-1}_V$ are isomorphisms follows
from the above compatibility with tensor products as in the proof of
\thmref{main, spherical}.

\medskip

\noindent{\it Remark.} Let us notice that the assertion of \thmref{representability, Iwahori}
is stronger than that of its spherical counterpart, \thmref{representability, spherical}.
Indeed, the former theorem, combined with \lemref{com pairing, Iwahori}, implies
the existence of an isomorphism \eqref{req iso} directly. In other words, the additional
argument, involving duailization is not necessary.

\section{Convolution and twisting for general chiral algebras, continued}   
  \label{sect 7}

\ssec{}

The goal of this section is to prove \thmref{representability,
Iwahori} and give another interpretation of the map \eqref{can
pairing, Iwahori}, parallel to what was done in \secref{proof of rep}
in the spherical case.

We shall first discuss a generalization of the construction of
\propref{global}.  Let $\CA$ be a chiral algebra as in \secref{general
alg}, and let $\{\CM_1,\CM_2\}\to \CM_3$ be a chiral pairing between
$\CA$-modules.

Let also $\CF'_1,\CF'_2,\CF'_3$ be chiral modules over
$\fD_{\kappa,G}$, and $\{\CF'_1,\CF'_2\}\to \CF'_3$ be a chiral
pairing.

\begin{propconstr} \label{very general}
In the above situation there exists a naturally defined chiral pairing
$$\{\CF'_1\tboxtimes \CM_1,\CF'_2\tboxtimes \CM_2\} \to
\CF'_3\tboxtimes \CM_3,$$ compatible with the actions of $\CA$ and
$\CA_{\fg,2\kappa_\crit}$.
\end{propconstr}

\begin{proof}

Each $\CF'_i\tboxtimes \CM_i$ is isomorphic to $\CF'_i\otimes \CM_i$
for $i=1,2,3$ as a D-module, and we define the desired pairing as the 
one coming from \eqref{mult chiral}. The fact that this pairing is compatible 
with the action of chiral algebras follows from the construction of
$\CF'_i\tboxtimes \CM_i$ in \secref{global case}.

\end{proof}

By multiplying the terms of the pairing of \propconstrref{very
general} by the chiral Clifford algebra we obtain a chiral pairing of
complexes of $\CA$-modules:
\begin{equation} \label{fund pairing}
\{\fC^\semiinf(L_\fg,\CF'_1\otimes
\CM_1),\fC^\semiinf(L_\fg,\CF'_2\otimes \CM_2)\} \to
\fC^\semiinf(L_\fg,\CF'_3\otimes \CM_3).
\end{equation}

Let us consider a particular case when $\CF'_1\simeq \fD_{G,\kappa}$,
$\CF'_2\simeq \CF'_3=:\CF'$ and the pairing is given by the chiral
action. Assume also that $\CM_1$ is $\on{Jets}(G)_X$-equivariant. In
this case the map $\CM_1\to \CO_{\on{Jets}(G)_X}\otimes \CM_1$ induces
a map of complexes of $\CA$-modules
$$\CM_1\to \fC^\semiinf(L_\fg,\CF'_1\otimes \CM_1).$$

\begin{lem}
The resulting chiral pairing
$$\{\CM_1, \fC^\semiinf(L_\fg,\CF'\otimes \CM_2)\}
\to \fC^\semiinf(L_\fg,\CF'\otimes \CM_3)$$
coincides with that of \eqref{global pairing}.
\end{lem}

\ssec{}

Let us assume now that in the previous set-up the chiral
$\fD_{G,\kappa}$-module $\CF'_1$ is $\on{Jets}(G)_X$-equivariant on
the right. Let $\CF'_2$ and $\CF'_3$ be both supported at $x\in X$,
and assume that as twisted D-modules on $G\ppart$ they are both
$I$-equivariant on the right.

Let $\CM_1$ be $\on{Jets}(G)_X$-equivariant, and let $\CM_2,\CM_3$ be
supported at $x\in X$ and $I$-equivariant. In this case the pairing
\eqref{fund pairing} gives rise to a chiral pairing
\begin{equation} \label{fund pairing, Iwahori}
\{\fC^\semiinf\Bigl(L_\fg;\fg,\CF'_1\otimes \CM_1\Bigr),
\fC^\semiinf\Bigl(L_\fg;\fh,\CF'_2\otimes \CM_2\Bigr)\}
\to \fC^\semiinf\Bigl(L_\fg;\fh,\CF'_3\otimes \CM_3\Bigr).
\end{equation}

Let us specialize further the case when $\CM_1\simeq \CA$,
$\CM_2\simeq \CM_3:=\CM$, with the chiral pairing being given by the
action. Let the terms of $\{\CF'_1,\CF'_2\}\to \CF'_3$ be those of
\eqref{map upstairs, Iwahori} for $n>>0$.  We obtain a chiral pairing
$$\{\fC^\semiinf\Bigl(L_\fg;\fg,\pi^*(\CF_{V,X-x})\otimes
\CA\Bigr)\otimes \CE_n,
\fC^\semiinf\Bigl(L_\fg;\fh,\delta_{I,G\ppart}\otimes \CM\Bigr)\}\to
\fC^\semiinf\Bigl(L_\fg;\fh,\pi_{\Fl}^*(\CZ_V)\otimes \CM\Bigr)$$

Let us assume now that $\kappa$ is such that $\CF_{V,X}\star \CA$ is
acyclic away from cohomological degree $0$. Passing to the $0$-th
cohomology in the previous expression, we obtain a functorial chiral
pairing
\begin{equation}  \label{main pairing, Iwahori}
\{\Bigl(\CF_{V,X-x}\star \CA)\otimes \CE_n\Bigr),\CM\}\to
h^0(\CZ_V\underset{I}\star \CM),
\end{equation}
such that the action of $N_{\CE}$ on the LHS goes over to the action of
$N_V$ on the RHS.

Parallel to the spherical case, we have the following generalization 
of \thmref{representability, Iwahori}:

\begin{thm}   \label{gen repr, Iwahori}  \hfill

\smallskip

\noindent{\em (1)} For any $\CM\in \CA\mod_x^I$ the convolution
$\CZ_V\underset{I}\star \CM$ is acyclic away from cohomological degree
$0$.

\smallskip

\noindent{\em (2)} For all $n$ that are large enough, the covariant
functor on $\CA\mod_x$ that sends an object $\CN$ to the set of chiral
pairings
$$\{\Bigl((\CF_{V,X-x}\star \CA)\otimes \CE_n\Bigr),\CM\}\to \CN$$
is representable by $\CZ_V\underset{I}\star \CM$.

\end{thm}

The proof of this theorem is based on the following construction. Let
$$\{\Bigl((\CF_{V,X-x}\star \CA)\otimes \CE_n\Bigr),\CM\}\to
\CN^\bullet$$ be a chiral pairing of complexes of $\CA$-modules, where
$\CN^\bullet$ is a strongly $\on{Jets}(G)_X$-equivariant complex.

Then for $m>>0$, from \eqref{fund pairing, Iwahori} we obtain a chiral
pairing of complexes
$$\{\fC^\semiinf\Bigl(L_\fg;\fg, (\pi^*(\CF_{V^*,X})\otimes \CE_m)
\otimes ((\CF_{V,X}\star \CA)\otimes \CE_n)\Bigr), \CM\}\to
\fC^\semiinf(L_\fg;\fh,\pi^*_{\Fl}(\CZ_{V^*})\otimes \CN^\bullet),$$
and by passing to the $0$-th cohomology and composing
with $\delta_{1,\Gr_{G,X}}\to \CF_{V,X}\star \CF_{V^*,X}$ and $\CE_{m+n}\to
\CE_m\otimes \CE_n$, we obtain a chiral pairing
$$\{(\CA\otimes \CE_{m+n}),\CM\}\to h^0(\CZ_{V^*}\underset{I}\star
\CN^\bullet).$$ Any such pairing factors through
$\CE_{m+n}\twoheadrightarrow \CO_X$, and therefore corresponds to a
map $\CM\to h^0(\CF_{V^*}\underset{I}\star \CN^\bullet)$.  The rest of
the argument repeats that of the the proof of \thmref{repr, general}.

\section{Appendix: proof of \propref{action as fusion}}   \label{app A}

\ssec{}

Let $\CF_X$ be a $\kappa$-twisted D-module on $\Gr_{G,X}$ and
let $\CF'_X$ denote its pull-back to $\on{Jets}^\mer(G)_X$.
Consider the D-module $({\bf 1}_{1,1})_!(\omega_X\boxtimes \CF_X)$
on $\Gr_{G,X^2}$. Corresponding to it there is a map
$$j_*j^*(\delta_{1,\Gr_{G,X}}\boxtimes \CF_X)\to \Delta_!(\CF_X)$$
of D-modules on $\Gr_{G,X^2}$. Lifting this map by means of
$\pi$ we obtain a map
\begin{equation}  \label{eq act as fusion}
j_*j^*\Bigl(\fD_{G,\kappa}\boxtimes
\Gamma(\on{Jets}^\mer(G)_X,\CF'_X)\Bigr)\to
\Delta_!\Bigl(\Gamma(\on{Jets}^\mer(G)_X,\CF'_X)\Bigr).
\end{equation}

Along with \propref{action as fusion} we will prove the following:

\begin{prop} \label{action as fusion moving} 
The map of \eqref{eq act as fusion} equals the map corresponding to
the chiral action of $\fD_{G,\kappa}$ on
$\Gamma(\on{Jets}^\mer(G)_X,\CF'_X)$.
\end{prop}

Applying this to $\CF_X=\delta_{1,\Gr_{G,X}}$, we obtain a description
of the chiral bracket on $\fD_{G,\kappa}$ in terms of distributions on
$\on{Jets}^\mer(G)_{X^2}$.

\ssec{}

Since the chiral algebra $\fD_{G,\kappa}$ is generated by
$\CO_{\on{Jets}(G)_X}$ and $L_{\fg,\kappa}$, to prove both
\propref{action as fusion} and \propref{action as fusion moving}, it
suffices to show that the two maps in question coincide when instead
of $\fD_{G,\kappa}$, as one of the multiples in the LHS, we take
$\CO_{\on{Jets}(G)_X}$ or $L_{\fg,\kappa}$.

The assertion concerning $\CO_{\on{Jets}(G)_X}$ follows tautologically
from the definition of the group ind-scheme
$\on{Jets}^\mer(G)_{X^2}$. In the case of $L_{\fg,\kappa}$ we shall
discuss the set-up of \propref{action as fusion}, while the case of
\propref{action as fusion moving} is similar. We need to establish
the following:

\medskip

Let $\bg$ be any point of $G\ppart/K$ and let ${\bf 1}_{1,1}(\bg)$ be
the corresponding section $X\to \Gr_{G,X;K}$. The normal to this
section, considered as a right D-module on $X-x$, has the property
that
\begin{equation} \label{ident normal}
\CN_{{\bf 1}_{1,1}(\bg)}|_{X-x}\simeq L_\fg|_{X-x}\oplus
T_\bg(G\ppart/K)\otimes \omega_{X-x} \text{ and } \CN_{{\bf
1}_{1,1}(\bg)}|_x\simeq T_\bg(G\ppart/K).
\end{equation}
In particular, we obtain a map
$$j_x{}_*j_x^*(L_\fg)\to i_x{}_!\Bigl(T_\bg(G\ppart/K)\Bigr),$$
which is equivalent to a map
$$\fg\ppart\simeq H^0_{DR}(\D^\times_x,L_\fg)\to T_\bg(G\ppart/K).$$

We need to show that the latter map equals the natural projection
$$\fg\ppart\twoheadrightarrow \fg\ppart/T_\bg(G\ppart/K),$$
corresponding to the left action of $G\ppart$ on $G\ppart/K$. Using
the above $G\ppart$-action, we reduce the assertion to the case when
$\bg$ is the unit point of $G\ppart$, in which case 
$T_\bg(G\ppart/K)\simeq \fg\ppart/\sk$, where $\sk$ is the Lie algebra of $K$.

Thus, we need to show that the map $L_\fg|_{X-x}\to \CN_{{\bf
1}_{1,1}(\bg)}|_{X-x}$ extends to a map of D-modules $L_{\fg.\sk}\to
\CN_{{\bf 1}_{1,1}(\bg)}$, where
$$L_{\fg.\sk}:=\on{ker}\Bigl(j_x{}_*j_x^*(L_\fg)\to
i_x{}_!(\fg\ppart/\sk)\Bigr).$$ This is a direct calculation performed
below.

\ssec{}

Without restriction of generality, we can assume that $X$ is affine,
and let $\xi(x_1,x_2)$ be a $\fg$-valued map on 
$X\times X-(\Delta_X\sqcup X\times x)$. 
We will show that $\xi(x_1,x_2)$ gives rise to a vertical
vector field on $\Gr_{G,X;K}$, relative to its projection onto $X$,
and that the arising normal vector field to the unit section vanishes
if and only if $\xi(x_1,x_2)$ has no poles at the diagonal, and for
any fixed $x_1$ the Laurent expansion of $\xi(x_1,\cdot)$ at $x_2=x$,
viewed as an element of $\fg\ppart$, belongs to $\sk$. This will imply
the required assertion.

\medskip

It will be more convenient to use a group-theoretic notation. I.e., we
will think of $\xi(x_1,x_2)$ as a map $X\times X-(\Delta_X\sqcup
X\times x)\to G$, parameterized by the scheme of dual numbers. More
generally, we will work with a map $\bg(x_1,x_2):\bigl(X\times
X-(\Delta_X\sqcup X\times x)\bigr)\times S\to G$ for an arbitrary
scheme $S$.

\medskip

To $\bg(x_1,x_2)$ as above we attach an $S$-valued automorphism of
$\Gr_{G,X;K}$ as follows. Given a point $(x_1,\CP_G,\beta,\alpha)$
of $\Gr_{G,X;K}$ we define a new point by leaving $(x_1,\CP_G,\alpha)$
the same, but multiplying $\beta$ by $\bg(x_1,\cdot)$, thought of
as a $G$-valued function on $X-(x_1\sqcup x)$.

\medskip

Applying this automorphism to the unit section of $\Gr_{G,X;K}$, the
resulting new $S\times X$-valued point of $\Gr_{G,X;K}$ will be
isomorphic to the initial one if and only if there exists an $S\times
X$-valued automorphism of $\CP^0_G$, preserving $\alpha$, and whose
value at any $x_1\in X$ equals that of $\bg(x_1,\cdot)$. But this
precisely means that $\bg(x_1,x_2)$ extends regularly to the diagonal
and $X\times x$, and the Taylor expansion of $\bg(x_1,\cdot)$ around
$x$ belongs to $K$.

 \end{document}